\documentclass[11pt,a4paper,reqno]{amsart}
\usepackage{amssymb,amsmath,graphics,verbatim}
\usepackage{latexsym}
\usepackage{eucal}
\usepackage{a4wide}
\usepackage[colorlinks]{hyperref}
\usepackage[normalem]{ulem}
\usepackage{soul}
\usepackage{cancel}
\usepackage{todonotes}

\usepackage{geometry}
\usepackage{marginnote}
\newcommand{\stkout}[1]{\ifmmode\text{\sout{\ensuremath{#1}}}\else\sout{#1}\fi}

\newcommand{\FLO}{{\rm FLO}}

\newcommand \ccl{{\rm ccl}}

\newcommand \Char{{\rm Char}}
\newcommand \B{\mathcal B}

\newcommand\R{\rr}

\newcommand{\rr}{\mathbb{R}}
\newcommand{\nn}{\mathbb{N}}

\newcommand{\0}{{\bf 0}}
\newcommand{\singsupp}{{\rm singsupp}}
\newcommand{\U} {{\mathcal U}}
\newcommand{\M}{{\mathcal M}}

\newcommand{\oo}{{\mathcal O}}
\newtheorem{theorem}{Theorem}
\newtheorem{lemma}[theorem]{Lemma}
\newtheorem{proposition}[theorem]{Proposition}
\newtheorem{corollary}[theorem]{Corollary}

\theoremstyle{definition}

\newtheorem{remark}[theorem]{Remark}

\numberwithin{equation}{section}
\numberwithin{theorem}{section}

\DeclareMathOperator \spn {span}

\DeclareMathOperator \WF {WF}

\DeclareMathOperator \supp {supp}

\title[Determining Riemannian Manifolds From Measurements at a Single Point]{Determining Riemannian Manifolds From Nonlinear Wave Observations at a Single Point}

\author{Leo Tzou}
\address{Leo Tzou\newline\indent School of Mathematics and Statistics, University of Sydney, NSW 2006, Australia}
\email{leo.tzou@gmail.com}

\date{Mar 18, 2021}

\begin{document}
\begin{abstract}
We show that on an a-priori unknown Riemannian manifold $(M,g)$, measuring the source-to-solution map for the semilinear wave equation at a single point determines its topological, differential, and geometric structure.


\end{abstract}

\maketitle

\section{Introduction}

The work of \cite{Kurylev:2014aa, Kurylev:2018aa, Kurylev:2014ab} developed the new idea of exploiting nonlinearity to solve inverse problems involving nonlinear hyperbolic equations. Roughly speaking, these results reconstruct the topological, differential, and conformal structure of a Lorentzian manifold from the source-to-solution data of some nonlinear hyperbolic equation measured in an open space-time cylinder of a timelike path. Other works involving the use of this idea to solve inverse problems for nonlinear hyperbolic equations can be found in \cite{Uhlmann:2020aa, Uhlmann:2018aa, Uhlmann:2019aa, Lassas:2018aa, Lassas:2017aa,Hoop:2019aa, Hoop:2020aa, Feizmohammadi:2019aa, Oksanen:2020aa, chen2021inverse}

Formally, the source-to-solution data considered in all of the above works are over-determined so it is natural to ask whether we really need so much data. Furthermore, these works model tomography where data is collected by {\em infinitely} many receivers, when in reality many imaging methods in seismology and astronomy are best modeled by only {\em one} receiver. Another motivation for reducing the amount of data necessary is that it can be seen as a nonlinear hyperbolic analogue to the ``partial data'' elliptic Calder\'on problem \cite{MR2299741}.

This paper provides a proof-of-concept model for a {\em single observer} inverse problems for nonlinear hyperbolic equations.
We show that when the Lorentzian metric is the direct product of $\R$ with a time-independent Riemannian manifold, the Riemannian manifold is determined by data collected at a {\em single spacial point}. Of course, the technique we develop here is far more robust than just the ultra-static Lorentzian case. In a forthcoming paper \cite{NOT} we will adapt these methods to recover conformal structure of general hyperbolic Lorentzian metrics with just a single observer.



Consider two Riemannian manifolds $(M_1, g_1)$ and $(M_2,g_2)$, and the corresponding Lorentzian manifold $(\M_j , -dt^2 + g_j)$ where $\M_j = \R \times M_j$. For each $j =1,2$, fix $p_j\in M_j$ and consider small open neighborhoods $\U_j\subset M_j$ containing $p_j$, and a diffeomorphism identifying $\U_1 \equiv \U_2$ and $p_1 \equiv p_2$. We will denote the common point to be $\0$ and the common neighborhood to be $\U$.

For a Riemannian manifold $(M,g)$, denote by $\Box_{g}=-\partial_t^2+\Delta_{g}$ the corresponding wave operator. The choice of nonlinearity we make here is cubic, and this is motivated by the cubic nonlinear interaction in Yang-Mills equations (see introduction of \cite{Chen:2019aa} for this motivation), 
\begin{eqnarray}
\label{cubic wave}
\Box_{g} u + u^3 = f {\rm\ in\ } \M   \\\nonumber
u = 0,\ t<-1.
\end{eqnarray}
This Cauchy problem has a unique solution (see e.g. \cite{Kurylev:2018aa} or \cite{Chen:2019aa}) if we assume $f$ is compactly supported in the cylinder $(-1,\infty)\times \mathcal U$ and has small $C^4_0((-1,\infty)\times \mathcal U)$ norm.

The data we wish to collect is the {\em single point} source-to-solution map
\begin{eqnarray}
\label{source to sol}
L_{M,g, 2T+1} f := u(t,\0) \mid_{t \in (0,2T+1)},\ \ \ f\in C^4_0((-1,2T+1)\times \mathcal U)
\end{eqnarray}
for a fixed $T>>1$. Note that while we are allowed to place sources in a cylindrical neighborhood of the timelike segment $\{(t,\0)\in \M\mid t\in (-1,2T+1)\}$, we are only collecting data on the time-like curve $(0, 2T+1)\times \0$. We denote by $B_{g_j}(x_j,r_j)$ the Riemannian open ball of $(M_j,g_j)$ centered at $x_j\in M_j$ with radius $r_j$.
Our result is the following:

\begin{theorem}
\label{metric recovery}
Let $(M_1,g_1)$ and $(M_2,g_2)$ be three-dimensional Riemannian manifolds. If their the source-to-solution map coincide $L_{M_1,g_1, 2T+1} = L_{M_2,g_2, 2T+1}$, then the Riemannian balls $(B_{g_1}(\0,T), g_1)$ and $(B_{g_2}(\0, T), g_2)$ are isometric. 
\end{theorem}
To the best of the author's knowledge, this is the first inverse problems result where measurements are done at a single point. Note also that a-priori the two Riemannian balls need not to have the same topological structure.
\begin{remark}
The intervals $[-1, 2T+1]$ in \eqref{cubic wave} and $(0,2T+1)$ in the definition of source-to-solution can be replaced by $[-\delta, 2T+\delta]$ and $(0, 2T+\delta)$ for any fixed $\delta >0$. As this article already uses an abundance of small parameters in proving Theorem \ref{metric recovery}, we consider the constant $1$ in place of a parameter $\delta>0$ to avoid having to keep track of yet another small parameter. 
\end{remark}

In this article we will assume that 
\begin{eqnarray}
\label{same as eucl}
g_1\mid_{\U} = g_2\mid_\U := g_0\mid_{\U}
\end{eqnarray}
 since deducing this fact from the source-to-solution map is trivial. We may assume without loss of generality that
\begin{eqnarray}
\label{mathcal U is small}
\mathcal U = B_{g_0}(\0; \delta_0)
\end{eqnarray}
where 
\begin{eqnarray}
\label{radius of inj delta0}
\delta_0 <\frac{1}{2}{\rm min}\left( {\rm inj}(B_{g_1}(\0, 2T), g_1), {\rm inj}(B_{g_2}(\0, 2T), g_2)\right)
\end{eqnarray}
is chosen small enough for $\U$ to be convex.

Our method uses the nonlinear interactions idea of \cite{Kurylev:2014aa, Kurylev:2018aa, Kurylev:2014ab}. We motivate the idea of our approach by considering the special case when $g_1$ and $g_2$ are both Riemannian metrics on $\R^3$ and assume the absence of conjugate points. Let $y_1\in \R^3$ and $x_1\in \U$. At time $t=0$ we create  a point source at $\0$ which produces a wave whose singularity propagates along the unit speed $g_1$-geodesic ray $[\0,y_1]_{g_1}$ joining $\0$ and $y_1$ (see \cite{Melrose:1979aa}). At another time $t_1$ (which may be greater or less than $0$), we send another singularity along the unit speed $g_1$-geodesic ray $[x_1, y_1]_{g_1}$ so that the two rays collide at $y_1$ at time $t =  R_0 := d_{g_1}(\0,y_1)$. Arrange for an auxiliary singularity to go along unit speed $g_1$-geodesic ray $[x_2, y_1]_{g_1}$ from $x_2\in\U$ to $y_1$ so that all three ray collide at the same time. When arranged judiciously (see Section \ref{sect nonlinear wf}), the nonlinear interaction given by \eqref{cubic wave} will result in a singularity coming back from $y_1$ to $\0$ along the $g_1$-geodesic ray $[\0,y_1]_{g_1}$. The single observer sitting at $\0$ will then see this singularity at time $t = 2R_0$. We then deduce, using the fact that source-to-solution maps are equal, that the same collisions must have occurred for the $g_2$ metric at time $R_0$ at some point $y_2$ which is distance $R_0$ away from $\0$. This allows us to conclude that $d_{g_1}(y_1, x) = d_{g_2}(y_2, x) = R_0 - t_1$ for all $x\in \U$ (see Prop \ref{distances are equal} and its corollary). This, together with a simple Riemannian argument (Prop. \ref{Riemannian lemma}), implies that the Riemannian balls $(B_{g_1}(\0,T), g_1)$ and $(B_{g_2}(\0, T), g_2)$ are isometric.

The difficulty in carrying out this program is that sources do not produce waves which propagate along a single light ray but rather traveling planes. The interaction of three traveling plane waves generically produce a spacelike curve as an artificial source (see Sect \ref{3fold}) and light emanating from such a curve can produce focal points even if we assume the absence of conjugate points. A single observer along a timelike path may not be able to ``see'' such singularity. See Remark \ref{counter eg}.

To overcome this issue we place small spherical ``mirrors'' near $\0$. When the location and timing of these mirrors are chosen judiciously, further nonlinear interaction between the returning light wave and these mirrors will reflect a visible {\em conormal} wave back to $\0$ (see Sect \ref{5fold}) which is visible to the lone observer, allowing us to circumvent the difficulty of focal points.

We note that it may be possible to avoid focal points by replacing the cubic nonlinearity in \eqref{cubic wave} by a quadratic nonlinearity and work with artificial point sources resulting from a four wave interaction. In addition to its relationship to Yang-Mills equations, we have chosen to consider the cubic nonlinearity because it motivates us to develop a robust method which has the potential to be applied to a broader range of problems. For example, this approach of using mirrors can potentially allow us to recover the ``earliest light observation set'' in the general Lorentzian setting (see \cite{Kurylev:2018aa}) from artificial sources to a narrow cylinder in space-time by only looking at singularity along a timelike curve contained in the cylinder.
\section*{Acknowledgement}
The author wishes to thank Marco Mazzucchelli, Alan Greenleaf, Suresh Eswarathasan, Mikko Salo, and Lauri Oksanen for the useful discussions. The author is supported by ARC DP190103302 and ARC DP190103451.

\section{Preliminaries}
\subsection{Basic Definitions}
Let $(M,g)$ be a Riemannian manifold. Associated to this Riemannian structure is a stationary Lorentzian metric $-dt^2 + g$ on the space-time 
\[\M = \{(t,x)\mid t\in \R, x\in M\}.\]
The time orientation on $\M$ is given by $\partial_t$.
It is sometimes convenient to have a Riemannian structure on $\M$ and the obvious choice is $G := dt^2 + g$.

Throughout this article, covectors in $\M$ are denoted by $\xi=\tau dt + \xi'$, where $\tau\in\R$ and $\xi'\in T^*M$. We denote by $\Char_g\subset T^*\M$ the subspace of lightlike covectors, which are those covectors that are multiple of $\mp dt + \xi'$, with $\xi'\in S^*M$. Here, the covectors $-dt + \xi'$ and $+ dt + \xi'$ are future and past pointing respectively. Notice that $\Char_g$ is the characteristic set of the homogeneous principal symbol of $\Box_g$.  

If $\Lambda_1, \Lambda_2\subset T^*M$ we use the notation
$$\Lambda_1 + \Lambda_2 := \{(t,x,\xi)\mid (t,x) \in \pi(\Lambda_1)\cap \pi(\Lambda_2),\ \xi = \xi_1 + \xi_2,\ {\rm with}\ (t,x,\xi_j) \in \Lambda_j\}.$$
When $\Lambda_1\cap T^*_{(t,x)}\M$ and $\Lambda_2\cap T^*_{(t,x)}\M$ are linearly independent for all $(t,x) \in \pi(\Lambda_1)\cap \pi(\Lambda_2)$, we will emphasize this fact with $\oplus$ in place of just $+$.

A smooth curve $s\mapsto(t(s),\gamma(s))\in\M$ is causal when its tangent vectors are timelike (i.e.\ $-\dot t(s)^2+\|\dot\gamma(s)\|_g^2<0$) or lightlike (i.e.\ $-\dot t(s)^2+\|\dot\gamma(s)\|_g^2=0$). The curve is future-directed if $\dot t(s)>0$, and past-directed if $\dot t(s)<0$. The future and past causal cones at $(t_0,x_0)\in \M$ are defined respectively by
\begin{align*}
I_g^+ (t_0,x_0) & = \{(t,x)\mid \mbox{ there exists a future-directed causal curve from}\ (t_0, x_0)\ \mbox {to}\ (t,x)  \},\\
I_g^- (t_0,x_0) & = \{(t,x)\mid \mbox{ there exists a past-directed causal curve from}\ (t_0, x_0)\ \mbox {to}\ (t,x)  \}.
\end{align*}

\subsection{Geodesics and Conjugate Points}

Let $H$ denote the Hamiltonian vector field on $T^*\M$ associated to the principal symbol of $\Box_g$, that is, to the Hamiltonian function
\begin{align}
\label{e:principal_symbol}
 (t,x,\tau dt+\xi')\mapsto \tfrac12\big(\|\xi'\|_g^2-\tau^2\big).
\end{align} 
The flow of $H$ preserves $\Char_g$, and its orbits on $\Char_g$ are precisely the lifts of the lightlike geodesics of $\M$. We define the future oriented flow on $\Char_g$ by 
\begin{eqnarray}
\label{def of e+}
e_+(s,t,x,\xi) = e_+^s(t,x,\xi) :=
\left\{
    \begin{array}{@{}ll}
       e^{sH}(t,x,\xi), & \mbox{if }\xi\ \mbox{ is future pointing,}  \vspace{5pt}\\
       e^{-sH}(t,x,\xi), & \mbox{if }\xi\ \mbox{ is past pointing}
    \end{array}
\right.
\end{eqnarray}
for $s\geq 0$. If $(t,x,\xi)\in \Char_g$ we denote its future oriented flowout by 
$$\FLO_g^+(t,x,\xi) :=\{ e_+^s(t,x,\xi)\mid s\geq 0\}.$$
Since the Lorentzian metric is a direct sum of $-dt$ with a Riemannian metric $g$, the projection to $M$ of the lightlike geodesics in $\M$ are Riemannian geodesics. If $(t,x,\xi)\in \Char_g$ and $\xi = -dt + \xi'$ with $\xi'\in S^*_xM$, then
\begin{eqnarray}
\label{explicit expression for LL geodesics}
\pi\circ e_+^s(t,x,\xi)= (t+s, \gamma_g(s, \xi'))
\end{eqnarray}
where $\gamma_g(s,\xi')$ is the unit speed geodesic in $M$ with initial condition $(x,\xi')$. We say that a Riemannian geodesic segment is length minimizing when it is the shortest curve joining its endpoints; we further say that it is unique length minimizing if every other curve joining the same endpoints is strictly longer.

%
%

\begin{lemma}
\label{uniform unique minimizer}
Assume that, for some $(y_0, \eta_0') \in S^*M$ and $R>0$, the geodesic segment $\gamma_g([0,R], \eta_0')$ is unique length minimizing and does not contain conjugate points.
Then, there exist $\delta>0$ and an open neighborhood $U'\subset S^*M$ of $(y_0, \eta'_0)$ such that, for each $y\in \pi (U')$, the map
\begin{align}
\label{lemma_diffeo_with_exp}
(0, R+\delta]\times (S^*_yM \cap U') \to M,
\qquad
(t, \eta') \mapsto \exp_y(t\eta')
\end{align}
is a diffeomorphism onto its image.
Up to reducing both $U'$ and $\delta$, for any $(y,\eta')\in U'$ and $t \in (0, R+\delta)$, the geodesics segment $\gamma_g([0,t], \eta')$ is unique length minimizing; in fact, there exists $c_0>0$ such that, for all $(y,\eta')\in U'$, $t\in  (R-\delta, R+\delta)$, $\hat \eta'\in S^*M\setminus\{\eta'\}$, and $\hat t >0$, if $\exp_y(t \eta') = \exp_y(\hat t \hat \eta')$ then $\hat t > t + c_0$.
\end{lemma}

\begin{proof}
The fact that~\eqref{lemma_diffeo_with_exp} is a diffeomorphism onto its image is well known due to the fact that the endpoints of $\gamma([0,R], \eta'_0)$ are not conjugate along the geodesic segment. As for the second part of the statement, assume by contradiction that there exist  sequences $S^*M\ni(y_j, \eta'_j)\to(y_0,\eta'_0)$, $t_j\to t_0\in(0, R]$, $\hat \eta_j' \in S^*_{y_j}M\setminus\{\eta_j'\}$, and $\hat t_j \leq t_j+ j^{-1}$ such that 
\begin{eqnarray}
\label{the other geodesic}
\gamma_g(\hat t_j,\hat\eta'_j) = \gamma_g(t_j, \eta'_j).
\end{eqnarray}
Since the map~\eqref{lemma_diffeo_with_exp} is a diffeomorphism onto its image, we must have that $\hat \eta_j' \notin U'$. Therefore, up to extracting a subsequence, $\hat \eta_j' \to \hat \eta' \in S^*_{y_0}M\backslash U'$ and $\hat t_j\to \hat t_0$. In particular $\hat \eta'\neq \eta'_0$. This implies $\gamma_g(\hat t_0,\hat\eta') = \gamma_g(t_0, \eta'_0)$, contradicting the fact that $\gamma_g([0,t_0], \eta'_0)$ is unique length minimizing.
\end{proof}

We recall that a Lagrangian submanifold $\Lambda \subset T^*\M$ is conic when $\lambda\xi\in\Lambda$ for all $\lambda>0$ and $\xi\in\Lambda$. The following Lemmas are classical, and we will not provide a proof.

\begin{lemma}
\label{N*K cap Char is smooth}
Let $K\subset \M$ be a spacelike curve then the intersection $N^*K \cap \Char_g\subset T^*\M$ is a three-dimensional submanifold.\hfill\qed
\end{lemma}
\begin{lemma}
\label{FLO is Lagrangian}
Let $\Lambda \subset T^*\M$ be a conic Lagrangian submanifold such that the intersection $\Lambda\cap \Char_g$ is a smooth three-dimensional submanifold nowhere tangent to $H$. Then its flowout
$\FLO_g^+(\Lambda\cap\Char_g)$
is Lagrangian. \hfill\qed
\end{lemma}

Let $\xi_0\in T_{(t_0,x_0)}^*\M$. For any $h>0$ small, we set 
\[ {\mathcal B}_h(\xi_0) := \big\{\xi \in T_{(t_0,x_0)}^*\M\ \big|\  \|\xi-\xi_0\|_{G} <h \big\} \]
Denote by $\ccl \B_h(\xi_0)$ its cone closure.
We recall that a hypersurface $S\subset\M$ is called lightlike when the Lorentzian metric $-dt^2+g$ restricts to a degenerate bilinear form on $TS$, or equivalently its conormal bundle consists of lightlike covectors.

\begin{lemma}
\label{conormal when not conjugate}
Let $\xi'_0\in S_{x_0}^*M$ and $t_0>0$ be such that the end points of the geodesic $\gamma([0,t_0],\xi'_0)$ are not conjugate. We set $\xi_0 := -dt + \xi'_0\in T^*_{(0,x_0)}\M\cap \Char_g$, and consider for $h>0$  small enough the Lagrangian submanifolds
\begin{align*}
\Lambda^\pm:=\FLO_g^+({\rm ccl}{\mathcal B}_h(\pm\xi_0)\cap\Char_g).
\end{align*}
Then, for any sufficiently small open neighborhood ${\mathcal V}\subset\M$ of $(t_0, \gamma(t_0,\xi'_0))$, there exists a lightlike codimension-one hypersurface $S\subset\mathcal V$ such that 
\[
(\Lambda^+\cup\Lambda^-)\cap T^*{\mathcal V} = N^*S. 
\]
\end{lemma}
\begin{proof}
Let $U\subset S_{x_0}^*M$ be a small open neighborhood of $\xi_0'$, and $\epsilon>0$ a small positive quantity. We define the maps
\begin{gather*}
 \iota_\pm:(t_0-\epsilon,t_0+\epsilon)\times U\to \Char_g\subset T^*\M,\\
 \iota_\pm(t,\xi')=e_+(t,0,x_0,\mp dt\pm\xi')=(t,\gamma(t,\xi'),\mp dt\pm\dot\gamma(t,\xi')^\flat),
\end{gather*}
where $\dot\gamma(t,\xi'):=\partial_t\gamma(t,\xi')$ and $\dot\gamma(t,\xi')^\flat:=g(\dot\gamma(t,\xi'),\cdot)$. Since the end points of the geodesic $\gamma([0,t_0],\xi'_0)$ are not conjugate, $\pi_2\circ\iota$ is a diffeomorphism onto its image, where $\pi_2:T\M\to M$ is the projection onto the factor $M$ of the base $\M = \R\times M$. 
If ${\mathcal V}\subset\M$ is a sufficiently small open neighborhood of $(t_0, \gamma(t_0,\xi'_0))$, and $\tilde U:=\iota^{-1}_+(T^*\mathcal V)=\iota^{-1}_-(T^*\mathcal V)$, then clearly
\begin{align*}
 \Lambda^\pm\cap T^*{\mathcal V}=\ccl(\iota_\pm(\tilde U)).
\end{align*}
The codimension-one hypersurface $S\subset\mathcal V$ of the statement is given by
\begin{align*}
 S:=\{(t,\gamma(t,\xi'))\ |\ (t,\xi')\in\tilde U\}.
\end{align*}
If $J_1,J_2$ are two linearly independent Jacobi fields along the geodesic $\gamma(\cdot,\xi')$ that are orthogonal to it and satisfy $J_1(0)=J_2(0)=0$, then the tangent spaces of $S$ are given by
\begin{align*}
T_{(t,\gamma(t,\xi'))}S=\mathrm{span}\{\partial_t+\dot\gamma(t,\xi'), J_1(t),J_2(t)\}.
\end{align*}
We claim that $\iota_\pm(\tilde U)\subset N^*S$. Indeed
\begin{align*}
(\mp dt\pm\dot\gamma(t,\xi')^\flat)(\partial_t+\dot\gamma(t,\xi')) &= \mp1 \pm \|\dot\gamma(t,\xi')\|_g^2 =0,\\
(\mp dt\pm\dot\gamma(t,\xi')^\flat)(J_i(t)) & = g(\pm\dot\gamma(t,\xi'),J_i(t)) =0.
\end{align*}
Since the conormal bundle $N^*S$ has rank one, we conclude $\ccl(\iota_+(\tilde U)\cup\iota_-(\tilde U))= N^*S$. Finally, since $\iota_\pm$ has image inside the characteristic set $\Char_g$, $S$ is a lightlike hypersurface.
\end{proof}

\begin{lemma}
\label{conormal intersection}
Let $S_j\subset\M$, $j=1,2,3$, be lightlike hypersurfaces with a triple transverse intersection at $(0,\0)$, i.e.
\begin{eqnarray}
\label{intersection of three is transverse}
{\rm Dim}\left( N^*_{(0,\0)} S_1 + N^*_{(0,\0)} S_2 + N^*_{(0,\0)} S_3\right) = 3.
\end{eqnarray}
Then there is an open set $\mathcal O\subset \M$ containing $(0,\0)$ so that $ K := S_1\cap S_2\cap S_3\cap \mathcal O $ is a spacelike segment. Furthermore,
\begin{eqnarray}
\label{NstarK}
N_{(t,x)}^*K = N^*_{(t,x)} S_1+N^*_{(t,x)} S_2+N^*_{(t,x)} S_3,\qquad\forall (t,x) \in K.
\end{eqnarray}
\end{lemma}
\begin{proof}
The transversality condition \eqref{intersection of three is transverse} ensures that near $(0,\0)$ the triple intersection $S_1\cap S_2\cap S_3$ is an open curve segment $K$ and that \eqref{NstarK} holds.

It remains to check that $K$ is spacelike. Observe that due to \eqref{NstarK} if $V\in T_{(0,\0)}K\subset T_{(0,\0)}\M$ then there are three independent covectors $\xi_j \in T^*_{(0,\0)}\M\cap \Char_g$, $j=1,2,3$, so that $\xi_j(V) = 0$. Without loss of generality $\xi_j = -dt + \xi'_j$ with $\xi'_j \in S^*M$. By linear independence we have that
\begin{eqnarray}
\label{xi'1 neq xi'2}
\xi_1'\neq \xi_2'
\end{eqnarray}

Without loss of generality we can write $V = a \partial_t + V'$ for some $V'\in S_{\0}M$. Using $\xi_1(V) = \xi_2(V) = 0$, we see that $\xi'_1(V') = \xi'_2 (V') = a$. Due to \eqref{xi'1 neq xi'2} we have that $|a|<1$. Therefore $V$ is spacelike and this implies that the segment $K$ is spacelike near $(0,\0)$.
\end{proof}

\begin{lemma}
\label{unique lightlike vector of lightlike surface}
Let $S\subset \M$ be a lightlike hypersurface and let $\xi\in N^*_{(t,x)}S$. If $V \in T_{(t,x)}S$ is a lightlike vector then $V \in \spn\{ \xi^\sharp\}$ where $\sharp$ is taken with respect to the Lorentzian metric $-dt^2+g$.
\end{lemma}
\begin{proof}
Write $\xi = -dt + \xi'$ with $\xi'\in S^*_xM$ then $\xi^\sharp = \partial_t + \xi'^\sharp$ where $\xi'^\sharp$ is the raising of index with respect to the Riemannian metric $g$ on $M$. Since $V\in T_{(t,x)} S$ is lightlike it is of the form (modulo multiplication by $\R$), $V= \partial_t + V'$ for some $V'\in S_x M$. By the fact that $\xi\in N^*_{(t,x)}S$, $\xi(V) = 0$ meaning that
$$1 = \xi'(V') = \langle \xi'^\sharp, V'\rangle_g.$$
Since both $\xi'^\sharp$ and $V'$ are unit vectors, we must have that $V' = \xi'^\sharp$. That is $V = \xi^\sharp$. 
\end{proof}

\begin{lemma}
\label{full span lem}
For all linearly independent covectors $\eta'_0, \eta'_1\in S^*_yM$, there exists $\eta'_{2} \in \spn\{ \eta_0', \eta_1'\}\cap S^*_yM$ arbitrarily close to $\eta'_{0}$ such that 
\begin{eqnarray}
{\rm Dim}\left(\spn\{-dt - \eta'_{0}, -dt - \eta'_{1}, -dt - \eta'_{2}\}\right) = 3.
\end{eqnarray}
\end{lemma}
\begin{proof}
We choose orthonormal coordinates on $T^*_yM$ so that 
\begin{align*}
 -\eta'_{1} = (1,0,0),
 \qquad
 - \eta'_{0} = (s, \sqrt{1-s^2},0),
\end{align*}
for some $s\in[-1,1]$.
Since $\eta'_{1}$ and $\eta'_{0}$ are linearly independent, we have $s \neq \pm1$. We choose $\eta'_{2}$ of the form $$-\eta'_{2} = (r ,\sqrt{1-r^2},0)$$ for $r\in(-1,1)\setminus\{s\}$ near $s$, so that $\eta'_2\in \spn\{\eta_1', \eta_0'\}$. The three vectors $-dt-\eta_i'$, $i=1,2,3$, are linearly independent if and only if the matrix
$$\begin{pmatrix}
-1 & 1 & 0 \\
-1 & s & \sqrt{1-s^2} \\
-1& r& \sqrt{1-r^2}
\end{pmatrix}$$
is non-singular, that is, if and only if $\frac{1-s}{1+s} \neq \frac{1-r}{1+r}$.
Since $t\mapsto \frac{1-t}{1+t}$ is monotonic for $t\in (-1,1)$, this inequality is always satisfied provided $r\in(-1,1)\setminus\{s\}$.
\end{proof}

\subsection{Analysis of PDE and Wave Front}

We first construct sources $f$ which are compactly supported near a point in $\M$ with wavefront set microlocalized near a single direction.

Let  $\xi_0' \in S^*_\0M$ and write $\xi_0 = - dt + \xi_0'\in T^*_{(0,\0)}\M$ as an element of $\Char_g$. Using local coordinates we can construct a homogeneous degree zero symbol $\omega_{0,h}(t,x,\xi) \in S^0(\M)$ such that $\omega_{0,h}(t,x,\xi) = 1$ for $\xi$ near $\pm\xi_0$ and  
$$\omega_{0,h}(t,x,D) \delta_{(0,\0)}\in I(\M; {\rm ccl}\mathcal B_h(\xi_0)\cup\ccl \B_h(-\xi_0)).$$
For a Lagrangian submanifold $\Lambda\subset T^*\M$ we define 
$$I(\M,\Lambda) := \bigcup_{m\in \nn} I_1^m(\M, \Lambda)$$
to be the space of Lagrangian distributions associated to $\Lambda$. We refer the reader to Definition 4.2.1 of \cite{Duistermaat:1996aa} for definition of $I_1^m(\M, \Lambda)$.

Set $\chi_{1,h} \in C_0^\infty(\M)$ such that $\chi_{1,h}(0,\0) = 1$ but $\supp(\chi_{1,h})\subset B_{G}((0,\0);h)$. 
Define the compactly supported conormal distribution $f_h \in I(\M; {\rm ccl}\mathcal B_h(\xi_0)\cup\ccl \B_h(-\xi_0))$ by 
\begin{eqnarray}
\label{source}
f_h (t,x):=\chi_{1,h}(t,x) \omega_{0,h}(t,x,D)\langle D\rangle^{-N} \delta_{(0,\0)}
\end{eqnarray}
for large fixed $N>0$ and elliptic classical pseudodifferential operator $\langle D\rangle$ of order $1$. We will often drop the $h$ subscript to simplify notation. We choose $N>0$ large enough so that $f_h(t,x) \in C^7_c(\M)$.

Let $u\in I(\M, \Lambda)$ be a Lagrangian distribution for some conic Lagrangian $\Lambda$ and 
$$\sigma(u) \in S^\mu(\Lambda; \Omega^{1/2}\otimes L)/ S^{\mu-1}(\Lambda;\Omega^{1/2}\otimes L)$$ 
be its principal symbol where $\Omega^{1/2}$ is the half-density bundle and $L \to \Lambda$ is the Maslov complex line bundle with local trivializations defined on open conic neighborhoods $U_j\subset \Lambda$. The transition functions of two local trivialization of $L$ on two intersecting open conic sets $U_j$ and $U_k$ are given by $i^{m_{j,k}}$ where $m_{j,k} \in\mathbb Z$. This observation allows us to define the notion of homogeneity of sections of $L$. Let $\rho >0$ and define $m_\rho : (t,x,\xi) \mapsto (t,x,\rho\xi)$ mapping $T^*\M$ to itself. If $\ell$ is a section of $L$ whose trivialization in the open conic neighborhood $U_j$ is given by the $\mathbb C$ valued function $\ell_j$, we can define $m_\rho^*\ell$ by the section whose local trivialization on $U_j$ is $m_\rho^*\ell_j$. Note that since the transition functions are constants (and thus remain unchanged under pullback by $m_\rho$), this definition is coordinate independent.

The fact that transition functions are constants also allows us to define the notion of a Lie derivative and pullback by a flow. Indeed if $X$ is a section of $T\Lambda$ and $\Phi^X$ denotes its flow, we define $(\Phi^X_s)^* \ell$ for $s$ small by $(\Phi^X_s)^*\ell_j(t,x,\xi) := \ell_j(\Phi^X_s(t,x, \xi))$ provided that both $\Phi^X_s(t,x,\xi)$ and $(t,x,\xi)$ are in $U_j$. Using the group property of $s\mapsto \Phi^x_s$ we can extend this definition for any $s\in \R$ by observing that transition functions between local trivializations are constants. The Lie derivative can then be defined by differentiating the pullback with respect to the parameter $s$.

It will also be useful to have local representations for sections of the half-density bundle $\Omega^{1/2}$ on $\Lambda$. Locally, a conic neighborhood $U\subset \Lambda$ can be written as 
$$U = \{(t,x, d_{t,x}\phi(t,x,\theta)\mid (t,x,\theta) \in C_\phi\}$$
for some $\mathcal V\subset \M$, $\Theta\subset \R^N$ a conic subset, $\phi: {\mathcal V} \times \Theta \to \R$ a homogeneous degree $1$ nondegenerate phase function, and 
$$C_\phi := \{(t,x,\theta) \in {\mathcal V}\times\Theta\mid d_\theta \phi(t,x,\theta) = 0\}.$$
Denote by $T_\phi: C_\phi \to U$ the diffeomorphism $(t,x,\theta)\in C_\phi \mapsto (t,x, d_{t,x}\phi(t,x,\theta))$. Define $d_{C_\phi} := (d_\theta \phi)^*\delta_{0,\R^N}$ to be a smooth section of $\Omega^1$ on $C_\phi$. Here $\delta_{0,\R^N}$ is the delta distribution on $\R^N$ at the origin.

To see that $d_{C_\phi}$ is degree $N$, it is useful to write down an explicit expression for $d_{C_\phi}$. Let $\lambda_1,\dots, \lambda_n$ be homogeneous degree $1$ local coordinates on $C_\phi$ which we extend smoothly to a neighourbhood of $C_\phi$ in $\mathcal V \times \R^N$. In these coordinates
\begin{eqnarray}
\label{dc in coordinates}
d_{C_\phi} = \left|\frac{D(\lambda_1,\dots, \lambda_n,\phi_{\theta_1},\dots,\phi_{\theta_N})}{D(t,x,\theta)} \right|^{-1} d\lambda_1\dots d\lambda_n.
\end{eqnarray}
The determinant is homogeneous of order $N-n$ and $d\lambda_1\dots d\lambda_n$ is order $n$. So $d_{C_\phi}$ is homogeneous of degree $N$.


If $[a] \in S^\mu(\Lambda, \Omega^{1/2} \otimes L)/ S^{\mu-1}(\Lambda, \Omega^{1/2} \otimes L)$ has homogeneous principal symbol and $(t,x,\xi)\in \Lambda$, we say that 
$$[a](t,x,\xi)\neq 0$$ 
if there is a an element $a\in [a]$ homogeneous of degree $\mu$ such that
$$a(t,x,\xi) \neq 0.$$

If $\ell_\phi$ is a non-vanishing homogeneous degree $0$ section of the pullback bundle $T^*_\phi L$ and $[a](t_0,x_0,\xi_0) \neq 0$ then under pullback by $T_\phi$,
$$T_\phi^* a(t_0,x_0,\theta_0) = a_\phi(t_0,x_0,\theta_0) \sqrt{d_{C_\phi}}\otimes \ell_\phi$$
for some $a_\phi(t,x,\theta) \in S^{\mu-N/2}({\mathcal V}\times \Theta)$ which satisfies
$$|a_\phi(t_0,x_0,\rho\theta_0) |\geq C_{t_0,x_0,\theta_0} |\rho|^{\mu-N/2}$$
for $\rho>0$ sufficiently large. 
\begin{lemma}
\label{transport equation}
Suppose $K\subset \M$ is a smooth submanifold whose conormal bundle $N^*K$ is transverse to $H$. Let $\Lambda = \FLO_g^+(N^*K\cap\Char_g)$ and let $a\in S^\mu(\Lambda, \Omega^{1/2}\otimes L)$ solve the transport equation 
$${\mathcal L}_H a = 0\ {\rm on}\ \Lambda,\ a\mid_{\partial \Lambda} = a_{\partial \Lambda}$$
for some $a_{\partial \Lambda} \in C^\infty(\partial\Lambda, \Omega^{1/2}\otimes L)$.

Suppose  $(m_\rho^*a)(t_0,x_0,\xi_0) = \rho^\mu a(t_0,x_0,\xi_0)$ and $a(t_0,x_0,\xi_0)\neq 0$ for some $(t_0,x_0,\xi_0)\in \partial \Lambda$ then $(m_\rho^*a)(t,x,\xi) = \rho^\mu a(t,x,\xi)$ and $a(t,x,\xi) \neq 0$ for all $(t,x,\xi)\in \FLO_{g}^+(t_0,x_0,\xi_0)$.
\end{lemma}
\begin{proof}
Assume that $(t_0, x_0,\xi_0), (t,x,\xi) \in \Lambda\cup\partial\Lambda$ can be represented by $(t_0,x_0,\theta_0), (t,x,\theta)\in C_\phi$ via the mapping $T_\phi$ and $(t,x,\xi) = e_+(s_0,t_0,x_0,\xi_0)$. In this proof we adopt the notation
$$e_s(t,x,\xi) := e_+(s,t,x,\xi).$$
To prove the lemma it suffices to show that if $a(t_0,x_0,\xi_0)\neq 0$ and $(m_\rho^*a)(t_0,x_0,\rho \xi_0) = \rho^\mu a(t_0,x_0,\xi_0)$ then
\begin{eqnarray}
\label{nonvanishing condition}
m_\rho^*a(t,x,\xi) = \rho^\mu a(t,x,\xi),\quad a(t,x,\xi)\neq 0.
\end{eqnarray}
To this end let $(t(s), x(s), \xi(s)) := e_s(t_0,x_0,\xi_0)$ for $s\in (0,s_0)$. Note that
\begin{eqnarray}
\label{rescale dynamic}
(t(s), x(s), \rho\xi(s)) = e_{s/\rho}( t_0,x_0, \rho\xi_0)
\end{eqnarray}
for $s\in (0, s_0/\rho)$.
 Since ${\mathcal L}_H a = 0$ we have that
\begin{eqnarray}
\label{pullback}
 a(t(s), x(s),\xi(s)) = (e_s^* a)(t(s), x(s), \xi(s))
 \end{eqnarray}
for all $s$. So using the rescaling formula \eqref{rescale dynamic} together with the invariance formula \eqref{pullback}, we see that \eqref{nonvanishing condition} is equivalent to 
 \begin{eqnarray}
\label{rescale and pullback}
m_\rho^* e_{-s_0/\rho}^* a(t, x, \xi) = \rho^\mu e_{-s_0}^* a(t,x,\xi),\quad e_{-s_0}^*a(t,x,\xi)\neq 0
\end{eqnarray}
 
Define $\tilde e_s := T_\phi^{-1} \circ e_s \circ T_\phi$ as the pullback flow defined on $C_\phi$ and write $(t(s),x(s), \theta(s)) = \tilde e_s(t_0,x_0,\theta_0)$. Using these coordinates we see that \eqref{rescale and pullback} is equivalent to


\begin{eqnarray}
\label{nonvanishing condition 2}
&&m_\rho^*\tilde e_{-s_0/\rho}^*\left((a_\phi \sqrt{ d_{C_\phi}}\otimes \ell_\phi) (t,x, \theta)\right) = \rho^\mu \tilde e_{-s_0}^*\left(a_\phi \sqrt{ d_{C_\phi}}\otimes \ell_\phi\right) (t,x,\theta).
\end{eqnarray}

We compute the LHS of \eqref{nonvanishing condition 2}. By definition 
\begin{eqnarray}
\label{scale the function}
m_\rho^*\tilde e_{-s_0/\rho}^*(a_\phi)(t,x,\xi) = \tilde e_{-s_0/\rho}^*a_\phi(t,x,\rho\xi) = a_\phi(t_0,x_0,\rho\xi_0) = \rho^{\mu-N/2}a_\phi(t_0,x_0,\xi_0).
\end{eqnarray}
The last equality comes from assumption that the symbol is homogeneous at $(t_0,x_0,\xi_0)$. The second last equality comes from the rescaling formula \eqref{rescale dynamic} and its analogous relation for the $\tilde e_s$ flow. 

The half-density $\sqrt{d_{C_\phi}}$ is 
given by \eqref{dc in coordinates} where the coordinate functions $\lambda_j$ are homogeneous of degree $1$. So using \eqref{rescale dynamic} we have that 
\begin{eqnarray}
\label{rescale density}
m_\rho^* \tilde e^*_{-s_0/\rho}\sqrt{d_{C_\phi}} = \rho^{N/2} \tilde e_{-s_0}^*\sqrt{d_{C_\phi}}
\end{eqnarray}
The section of the Maslov bundle $\ell_\phi$ is homogeneous of degree $0$ so 
\begin{eqnarray}
\label{rescale maslov}
m_\rho^*\tilde e_{-s_0/\rho}^*\ell_\phi = \tilde e_{-s_0}^*\ell_\phi.
\end{eqnarray}
Substituting \eqref{scale the function}, \eqref{rescale density}, and \eqref{rescale maslov} into the LHS of \eqref{nonvanishing condition 2} we have that \eqref{nonvanishing condition 2} is equivalent to
\begin{eqnarray}
\label{nonvanishing condition 3}
 \rho^{\mu}a_\phi(t_0,x_0,\xi_0)\tilde e_{-s_0}^*\left(\sqrt{d_{C_\phi}}\otimes  \ell_\phi \right)(t,x,\theta)= \rho^\mu \tilde e_{-s_0}^*\left(a_\phi \sqrt{ d_{C_\phi}}\otimes \ell_\phi\right) (t,x,\theta).
\end{eqnarray}
Observing that $\tilde e_{-s_0}^*a_\phi(t,x,\theta) = a_\phi(t_0,x_0,\theta_0)$ verifies \eqref{nonvanishing condition 3}.
\end{proof}

Two Lagrangians $\Lambda_0, \Lambda_1\subset T^*\M$ are said to be an intersecting pair if $\Lambda_0\cap \Lambda_1 = \partial \Lambda_1$ and
$$T_\lambda \Lambda_1 \cap T_\lambda \Lambda_0 = T_\lambda \partial \Lambda_1,\ \forall \lambda\in \partial\Lambda_1$$
In this case for all $m\in {\mathbb Z}$, Definition 3.1 of \cite{Melrose:1979aa} defined the class of paired Lagrangian distributions $I^m(\M,\Lambda_0,\Lambda_1)$. We denote by $I(\M, \Lambda_0, \Lambda_) := \bigcup_{m\in {\mathbb Z}} I^m(\M, \Lambda_0,\Lambda_1)$.

\begin{proposition}
\label{propagation}
Let $K\subset \M$ be a smooth submanifold and $f \in I(\M,N^*K)$ where $N^*K$ is transverse to $H$ and assume that $f$ has homogeneous principal symbol. Let $u$ be a solution of 
\begin{eqnarray}
\label{linear wave source}
\Box_g u = f,\ \ u\mid_{t<-1} = 0.
\end{eqnarray}
Then $N^*K$ and $\Lambda := \FLO_g^+(N^*K\cap \Char_g)$ form an intersecting pair of Lagrangians. The solution $u$ belongs to $u \in I(\M, N^*K, \Lambda )$ with homogeneous principal symbol on $\Lambda $. Furthermore if $(t_0,x_0,\xi_0)\in \WF(f)\cap\Char_g$ with $\sigma(f)(t_0,x_0,\xi_0)\neq 0$, then $\FLO_g^+(t_0,x_0,\xi_0)\subset \WF(u)$ with $\sigma(u)(t,x,\xi)\neq 0$ for all $(t,x,\xi) \in \FLO_g^+((t_0,x_0,\xi_0))$.
\end{proposition}
\begin{proof}
 It suffices to construct a parametrix for \eqref{linear wave source} satisfying $\Box_g u -f \in C^\infty$. This was done in Proposition 6.6 of \cite{Melrose:1979aa} where the homogeneous principal symbol $\sigma(u)$ satisfies
 $${\mathcal  L}_H \sigma(u) = 0\ {\mbox on}\ \Lambda, \  \sigma(u)\mid_{N^*K \cap \Lambda}= {\mathcal R} \left( (\tau^2 - \|\xi'\|_g^2)^{-1}\sigma(f) \right) \mid_{N^*K \cap \Lambda}.$$
 Here $\mathcal R$ is the operator defined via Definition (4.7) in \cite{Melrose:1979aa}.
By Lemma \ref{transport equation} it suffices to show that if $(t_0,x_0,\xi_0)\in N^*K \cap\Char_g$ satisfies $\sigma(f)(t_0,x_0,\xi_0)\neq 0$ then 
$${\mathcal R} \left( (\tau^2 - \|\xi'\|_g^2)^{-1}\sigma(f) \right)(t_0,x_0,\xi_0) \neq 0.$$
To this end observe that $H$ is the Hamiltonian vector field of the Hamiltonian $h_4(t,x,\tau dt+\xi') := \tau^2 -\|\xi'\|_g^2$ whose zero energy set is $\Char_g$. So since the Hamiltonian vector field $H$ is transverse to the Lagrangian $N^*K$, $N^*K\cap \Char_g$ is a smooth three-dimensional submanifold of $T^*\M$. So near $(t_0,x_0,\xi_0)$ we can choose $h_1,\dots h_{3}$ to be a coordinate system for $N^*K \cap \Char_g$. We then have that $\{h_1,\dots, h_4\}$ is a coordinate system for $N^*K$ in a neighbourhood of $(t_0,x_0,\xi_0)$. Note that since $H$ is the Hamiltonian flow of $h_4$, $h_4 = 0$ on $\FLO_g^+(N^*K\cap \Char_g)$ and
$$N^*K\cap \Char_g = N^*K \cap \Lambda = \partial \Lambda.$$


The Lagrangian $\Lambda$ is diffeomorphic to $\R^+ \times \partial \Lambda$ via the map $(s,t,x,\xi) \mapsto e_+(s, t,x,\xi)$. Define $h_5\in C^\infty(\Lambda)$ by $h_5(e_+(s,t,x,\xi)) = s$ with $(t,x,\xi)\in \partial \Lambda$. Clearly
$$\{h_4,h_5\} = H h_5 = 1>0$$
and $\{h_1, h_2, h_3, h_5\}$ forms a local coordinate system for $\Lambda$ near $\FLO_g^+(t_0,x_0,\xi_0)$. By (4.7) of \cite{Melrose:1979aa}, in these coordinates if 
$$\sigma(f) = \hat f |dh_1\wedge\dots\wedge dh_4| ^{1/2}\otimes \ell$$
then
$${\mathcal R} \left( h_4^{-1}\sigma(f) \right) = \hat f |dh_1\wedge dh_2\wedge dh_3\wedge dh_5| ^{1/2}\otimes \ell.$$
So by our assumption that $\sigma(f)(t_0,x_0,\xi_0)\neq 0$, we can conclude that 
$${\mathcal R} \left(h_4^{-1}\sigma(f) \right)(t_0,x_0,\xi_0)\neq 0$$ thus completing the proof. 

\end{proof}
\begin{lemma}
\label{WF of product}
Let $K_j \subset \M$, $j=1,2,3$ be submanifolds such that $K_1$ is transverse to $K_2$ and $K_3$ is transverse to $K_1\cap K_2$. Let $u_j\in I(\M; N^*K_j)$ for $j=1,2,3$.\\
i) Let $\omega(t, x, D)$ be a pseudodifferential operator whose wavefront is disjoint from $N^*K_1 \cup N^*K_2$. Then 
$$\omega(t,x,D)(u_1 u_2) \in I(\M, N^*(K_1\cap K_2))$$
with principal symbol given by
$$\sigma(u_1 u_2)(t,x,\xi) = \omega(t,x,\xi) \sigma(u_1)(t,x,\xi^{(1)})\sigma(u_2)(t,x,\xi^{(2)})$$
where $\xi = \xi^{(1)} + \xi^{(2)}$.\\
ii) Let $\omega(t,x,D)$ be a pseudodifferential operator whose wavefront is disjoint from 
$$\bigcup_{j=1}^3 N^*K_j \cup \bigcup_{j,k=1,j\neq k}^3 N^*\left(K_j\cap K_k \right)$$
then 
$$\omega(x,D)(u_1 u_2u_3) \in I(\M, N^*(K_1\cap K_2\cap K_3))$$
with principal symbol given by
$$\sigma(u_1 u_2u_3)(\xi) = \omega(x,\xi) \sigma(u_1)(\xi^{(1)})\sigma(u_2)(\xi^{(2)})\sigma(u_3)(\xi^{(3)})$$
where $\xi = \xi^{(1)} + \xi^{(2)} + \xi^{(3)}$.\\
\end{lemma}
\begin{proof}
Part i) is Lemma 2 of \cite{Chen:2019aa}. Part ii) is proof of equation (56) in \cite{Chen:2019aa}. Note that both are based on the work of \cite{Greenleaf:1993aa}.
\end{proof}

\begin{lemma}
\label{only threefold propagate}
Suppose that $S_1$ and $S_2$ are lightlike hypersurfaces of $\M$ of codimension one and intersects transversely,  
then
$$ N^*(S_1 \cap S_2) \cap \Char_g \subset (N^*S_1 \cup N^*S_2).$$
\end{lemma}
\begin{proof}
This is a trivial linear algebra calculation. See e.g. (30) in \cite{Chen:2019aa}.
\end{proof}

\begin{lemma}
\label{restriction has singularity}
If $u\in I(\M, N^*S)$ is a conormal distribution for some lightlike hypersurface $S$ containing the point $(0,\0)$. Suppose $\sigma(u)(0,\0,\xi_0)\neq0$, then for all timelike curve $\alpha(\cdot)$ which intersects $S$ transversely at $(0,\0)$ the distribution $\alpha^* u \in {\mathcal D}'(\R)$ is singular at the point $0\in \R$. 
\end{lemma}
\begin{proof} 
Since everything is local we may assume that $\M =\R^{1+3}$ endowed with Lorentzian metric $-dt^2+g$ where $g$ is a time-independent Riemannian metric on $\R^3$.
We choose local coordinates $z= (z_0,z_1,z_2,z_3)\in \R^{1+3}$ so that $S = \{z_0 = 0\}$, $\xi_0 = dz_0$, and $\alpha(t) = (\alpha_0(t), \alpha_1(t), \alpha_2(t), \alpha_3(t))$ with $\alpha(0) = (0,0,0,0)$. Since the intersection at the origin is transverse,
\begin{eqnarray}
\label{dot alpha}
\dot \alpha_0(0) \neq 0
\end{eqnarray}
In these coordinates, since $\sigma(u)(0,\0,\xi)\neq0$,
$$u(z) = \int_\R e^{i\theta z_0} a(z, \theta) d\theta$$
where for some $N\in \R$, the symbol $a(z,\theta)$ satisfies
\begin{eqnarray}
\label{a is positive}
|a(0,\theta)|\geq C |\theta|^{N}
\end{eqnarray}
for all $\theta >1$. The pullback $\alpha^*u(t)$ is then 
\begin{eqnarray}
\label{pullback distribution}
\alpha^*u(t) = \int_\R e^{i\theta \alpha_0(t)} \tilde a(t,\theta) d\theta
\end{eqnarray}
where $\tilde a(t,\theta) = a(\alpha(t),\theta)$. By \eqref{dot alpha} the map $(t,\theta)\mapsto \alpha_0(t)\theta$ is a non-degenerate phase function. So \eqref{pullback distribution} is a oscillatory distribution with nondegenerate phase. For any $\chi\in C^\infty_c(\R)$ supported in a neighborhood of $0$, we can take the Fourier Transform of $\chi \alpha^*u$. Apply Thm 2.3.1 of \cite{Duistermaat:1996aa} we have that due to \eqref{a is positive}
$$|{\mathcal F}(\chi\alpha^*u)(\tau)| \geq  C|\tau|^{N'}$$
for all $\tau >1$. Therefore $0\in \singsupp(\alpha^*u)$.
\end{proof}
\begin{remark}
\label{counter eg}
Lemma \ref{pullback distribution} states that an observer along a timelike world line is able to see a the singularity of distribution whose wavefront is the conormal bundle of a lightlike hypersurface. This is not true in general for lightlike Lagrangians.

Consider the the Minknowski metric $-dt^2+{\rm eucl}$ on $\R^{1+3}$ and the corresponding fundamental solution 
$$\Box_{\R^{1+3}} G(t,x) = 0,\ \ G(0,x) = 0,\ \partial_t G(0,x) = \delta_\0(x)$$
where $\delta_\0(x)$ is the delta distribution at the origin $\0\in \R^3$. The distribution $G(\cdot,\cdot)$ acts on $\varphi\in C^\infty_c(\R^{1+3})$ by
$$\langle \varphi, G(\cdot, \cdot)\rangle_{\R^{1+3}}  = \int_{-\infty}^\infty t^{-1} \int_{\partial B_{\rm eucl}(\0, |t|)} \varphi(x) {\rm dVol}_t(x) dt$$
where ${\rm dVol}_t$ is the volume form of the round sphere of radius $|t|$. This is a Langrangian distribution on $\R^{1+3}$ which is not conormal at $(0,\0)\in \R^{1+3}$. Take the curve 
$$\alpha(s) = (s,\0),\ s \in \R.$$
Then $(\alpha^*G)(s) = 0$ for all $s\in \R$. 
\end{remark}

We conclude the section with a lemma about the linearization of solutions of \eqref{cubic wave}. Henceforth if $\beta \in \nn^5$ is a multi-index, we denote by $\beta_j\in \nn^5$ to be the element with $1$ in the $j$th slot and zero elsewhere. If $j,k\in \{0,1,2,3,4\}$ are distinct, we denote by $\beta_{j,k}\in \nn^5$ to be the element with $1$ in the $j$th and $k$th slot and zero elsewhere. The multi-indices $\beta_{j,k,l}$, $\beta_{j,k,l,m}$ are defined analogously. We also denote by $\hat \beta\in \nn^5$ to be the element $(1,1,1,1,1)$.
\begin{lemma}
\label{expansion of sol}
Let $u$ solve \eqref{cubic wave} with source give by 
$\sum_{l=0}^4 \epsilon_l f_l$ where 
$$f_l\in C_c^N((-1,2T+1) \times \U)$$
 for large but finite $N$. Then $u = u_{\rm ans} + u_{\rm r}$ where
$$u_{\rm ans} := \sum_{\beta\in \nn^5, |\beta|\leq 5}\epsilon^\beta u_{\beta}$$ 
and $\|u_{\rm r}\|_{C([-1,2T+1]\times M)} \leq C \sum\limits_{6\leq |\beta|\leq 15}\epsilon^\beta$.\\
i) We have
$$\Box_g u_{\beta_j} = f_j,\ u_{\beta_j}\mid_{t<-1} = 0.$$
ii) If $|\beta | = 2$ then $u_\beta = 0$. \\
iii) We have 
$$\Box_g u_{\beta_{j,k,l}} = -6 u_{\beta_j}u_{\beta_k} u_{\beta_l} ,\ u_{\beta_{j,k,l}}\mid_{t<-1} = 0.$$
iv) If $|\beta| = 4$ then $u_\beta = 0$.\\
v) We have that $u_{\hat \beta}$ solves
$$\Box_g u_{\hat\beta} = -\sum\limits_{\substack{j,k,l,m,n=0\\ j\neq k\neq l\neq m\neq n}}^4 u_{\beta_{j,k,l}} u_{\beta_m} u_{\beta_n},\ u_{\hat \beta}\mid_{t<-1} = 0.$$
\end{lemma}
\begin{remark}
If $u$ is as in Lemma \ref{expansion of sol}, we denote by $u_j := u_{\beta_j}$ as in i) of the lemma, $u_{jk} := u_{\beta_{j,k}}$ as in ii) of lemma,  $u_{jkl} := u_{\beta_{j,k,l}}$ as in iii) of the lemma, $u_{jklm} := u_{\beta_{j,k,l,m}}$ as in iv) of lemma, and $u_{01234} := u_{\hat \beta}$ as in v) of the Lemma. 
\end{remark}
\begin{proof}
For each multiindex $\beta \in \nn^5$ satisfying $|\beta|\leq 5$ choose $u_\beta$ with $u_\beta \mid_{t<-1} = 0$ so that the ansatz
$$u_{\rm ans} := \sum_{\beta\in \nn^5, |\beta|\leq 5}\epsilon^\beta u_{\beta}$$ solve the approximate \eqref{cubic wave}:
$$\Box_g u_{\rm ans} + u^3_{\rm ans} = \sum_{l=0}^4 \epsilon_l f_l + \sum_{6\leq|\beta|\leq 15}\epsilon^\beta R_{\beta}.$$
where each $R_{\beta} \in C^{N'}([-1,2T+1]\times M)$ for some large $N'$ and $R_{\beta}\mid_{t<-1} = 0$. Direct (but cumbersome) calculation shows that i) - v) are satisfied.

Using the result for solving (172) in Appendix C of \cite{Kurylev:2013aa} we can find $u_{\rm r} \in C([-1, 2T+1]\times M)$ solving
$$\Box_g u_{\rm r} + 3 u_{\rm r} u_{\rm ans}^2 + u_{\rm r}^3   + 3 u_{\rm r}^2 u_{\rm ans} = \sum_{6\leq|\beta|\leq 15}\epsilon^\beta R_{\beta}.$$
Furthermore we have the well-posedness estimate
\begin{eqnarray}
\label{remainder estimate}
\|u_{\rm r} \|_{ C([-1,2T+1]\times M)} \leq \sum_{6\leq|\beta|\leq 15}\epsilon^\beta \|R_{\beta}\|_{C^{N'}([-1,2T+1]\times M))}.
\end{eqnarray}
Simple calculation shows that $u = u_{\rm ans} + u_{\rm r}$ solves \eqref{cubic wave}.
\end{proof}

It is useful to have an explicit formula for $u_{01234}$ in terms of the solution $u$ of \eqref{cubic wave} with source $\sum_{l=0}^4\epsilon_l f_l$. To this end we write $\epsilon := (\epsilon_0,\dots, \epsilon_4)$ and adopt the notation $K_0(\epsilon) := (0,\epsilon_1,\dots, \epsilon_4)$. That is, we replace the first entry of $\epsilon$ by $0$. Define $K_j(\epsilon)$ accordingly for $j = 0,\cdots, 4$. Similarly for distinct $j,k\in\{0,\dots,4\}$, $K_{j,k}(\epsilon)$ replaces the $j$th and $k$th entry of $\epsilon$ with $0$ and leaves the rest unchanged. Define the operators $K_{j,k,l}$ and $K_{j,k,l,m}$ accordingly for distinct $j,k,l,m\in \{0,\dots, 4\}$.

To highlight the dependence of $u$ on the parameter $\epsilon$ we write $u = u(\epsilon)$ in the following discussion. Using the estimate for $u_{\rm r}$ in Lemma \ref{expansion of sol} and the formal expansion of $u_{\rm ans}$, we see that
\begin{lemma}
\label{differentiate then restrict}
If $\alpha : [-1,2T+1]\to\M$ is a smooth timelike curve then 
$$\alpha^*u_{0} =  \lim\limits_{\substack{\epsilon_0\to 0\\ \epsilon_1=\dots=\epsilon_4 = 0}} \frac{\alpha^*u(\epsilon)  }{\epsilon_0}.$$

$$\alpha^*u_{012} =  \lim\limits_{\substack{\epsilon\to 0,\\ \epsilon_0 = \epsilon_1 = \epsilon_2,\\ \epsilon_3 = \epsilon_4 = 0}} \frac{\alpha^*u(\epsilon) - \sum\limits_{j=0}^2 \alpha^*u(K_j\epsilon)  +  \sum\limits_{\substack{j,k=0,\\ j \neq k}}^2 \alpha^*u(K_{j,k}\epsilon) }{\epsilon_0\dots \epsilon_{2}}.$$ 

{\small$$\alpha^* u_{01234} = \lim\limits_{\substack{\epsilon\to 0,\\ \epsilon_0 = \dots = \epsilon_4}} \frac{\alpha^*u(\epsilon) - \sum\limits_{j=0}^4 \alpha^*u(K_j\epsilon)  +  \sum\limits_{\substack{j,k=0,\\ j \neq k}}^4 \alpha^*u(K_{j,k}\epsilon)- \sum\limits_{\substack{j,k,l=0,\\ j \neq k\neq l}}^4 \alpha^*u(K_{j,k,l}\epsilon)+ \sum\limits_{\substack{j,k,l,m=0,\\ j \neq k\neq l\neq m}}^4 \alpha^*u(K_{j,k,l,m}\epsilon) }{\epsilon_0\dots \epsilon_{4}}.$$}


We formally write the limit on the right hand side as 
$$\partial_{\epsilon_0} (\alpha^*u) \mid_{\epsilon_0 = \dots = \epsilon_4=0},$$ 
$$\partial_{\epsilon_0\dots\epsilon_2}^3 (\alpha^*u) \mid_{\epsilon_0 = \dots = \epsilon_4=0},$$
$$\partial_{\epsilon_0\dots\epsilon_4}^5 (\alpha^*u) \mid_{\epsilon_0 = \dots = \epsilon_4=0}$$ 
etc.

\end{lemma}
Formally, we can view the result of Lemma \ref{differentiate then restrict} as the commutation of the restriction operator with differentiation with respect to $\epsilon$.
\section{Analysis of Multilinear Wave Interaction}
\label{sect nonlinear wf}

\subsection{Lightlike Flowout from Spacelike Curves}
\label{Lightlike Flowout from Spacelike Curves}
In this subsection we deduce some conormal properties for the Langrangian generated by flowing out by lightlike covectors on the normal bundle of a spacelike curve. The projection of these Lagrangians are in general very irregular objects. However, we deduce sufficient conormal properties for us to carry on our analysis later. 

Let $K\subset \M$ be a smooth and short spacelike segment containing the point $(0,y_0)$. Suppose $-dt + \eta'_0 \in N^*_{(0,y_0)}K$ is a lightlike covector. Let $\gamma_g([0,R],\eta'_0)$ be the unique length minimizing geodesic between the initial point  $y_0$ and the end point $x_0\in \mathcal U$ which satisfies $d_g(x_0,\0)<\delta_0/4$ where $\delta_0$ is defined in \eqref{radius of inj delta0}.  Assume the segment does not contain conjugate points. Note though that there may be other (longer) geodesics joining the two end points. Set  
$$\Lambda := \FLO_g^+(N^*K\cap \Char_g)$$
to be the flowout and define the point $\lambda_0 \in \Lambda$ to be 
$$\lambda_0 := (R, x_0, -dt + \xi'_0)$$
where
 $\xi'_0= \dot\gamma_g(R, \eta'_0)^\flat$.
The set $N^*K\cap \Char_g$ is a smooth 3-manifold by Lemma \ref{N*K cap Char is smooth} and so $\Lambda$ is a smooth Lagrangian submanifold of $T^*M$ by Lemma \ref{FLO is Lagrangian}. Observe that the future oriented flow \eqref{def of e+} gives a diffeomorphism
\begin{align}
\label{e:restriction_e+}
\R^+_s \times (N^*K\cap\Char_g)   \to \Lambda,\qquad
 (s, t,x,\xi)  \mapsto e_+(s,t,x,\xi).
\end{align}
Due to their conic nature, $N^*K\cap \Char_g$ and $\Lambda$ are sometimes inconvenient to work with, as the differential of the base projection acting on these manifolds has at least a one dimensional kernel. For this reason, we will sometimes work with the fiberwise quotients by $\R^+$ of the complements of the zero-sections of $N^*K \cap \Char_g$ and $\Lambda$; in order to simplify the notation, we will not indicate the removal of the zero-section in the notation, and simply write $(N^*K \cap \Char_g)/\R^+$ and $\Lambda/\R^+$ for such quotients. We will see them by means of the identifications
\begin{align}
\label{quotient conormal}
(N^*K \cap \Char_g)/\R^+ &\cong \{(t,x,\xi)\in N^*K \mid \xi = \pm dt + \xi',\ \xi'\in S_x^*M\},\\
\label{quotient lagrangian}
\Lambda/\R^+ &\cong \{(t,x,\xi)\in \Lambda \mid \xi = \pm dt + \xi',\ \xi'\in S_x^*M \}.
\end{align}
The future oriented Hamiltonian flow~\eqref{e:restriction_e+} induces a diffeomorphism of the quotients
\begin{equation}
\label{PhiH explicitly}
\begin{gathered}
\R^+ \times (N^*K\cap \Char_g)/\R^+  \to \Lambda/\R^+ \\
e_+(s,r,x,\mp dt + \eta')  = (r+ s, \gamma_g(s,\pm \eta'), \mp dt \pm \dot\gamma_g(s,\pm\eta')^\flat). 
\end{gathered}
\end{equation}

\begin{proposition}
\label{flowout lagrangian is almost conormal}
There exists sequences
\begin{align*}
\Lambda \ni \lambda_j = (R_j, x_j, -dt + \xi'_j) & \to \lambda_0 = (R, x_0, -dt + \xi'_0), \\
(N^*K\cap \Char_g)/\R^+ \ni \eta_j = (r_j, y_j, -dt + \eta'_j) & \to (0,y_0, -dt+\eta'_0),\\
s_j & \to R,
\end{align*}
with $x_j\neq x_0$, $R_j> R$, $e_+(s_j, \eta_j) = \lambda_j$, 
open neighborhoods $\mathcal O_j\subset \M$ of $(R_j, x_j)$, open neighborhoods $U_j\subset (N^*K\cap\Char_g)/\R^+$ of $\eta_j$, quantities $\delta_j>0$, and open neighborhoods $\Gamma_j\subset K$ of $(r_j, y_j)$ such that 
\begin{align*}
S_j:= \pi\circ e_+((s_j-\delta_j, s_j+\delta_j)\times U_j)
\end{align*}
is a lightlike codimension-one hypersurface in $\M$ satisfying
$$T^* \mathcal O_j \cap \FLO_g^+(N^*\Gamma_j\cap \Char_g) = N^*(S_j\cap\mathcal O_j).$$
Finally we can arrange the sequences $\lambda_j$, $S_j$, and $\mathcal O_j$ so that 
\begin{eqnarray}
\label{flowout doesn't hit (t,x0)}
\pi\circ\FLO_g^+(N^*(S_j \cap \mathcal O_j)) \cap\left( [R, R+\delta_0/2] \times x_0 \right)= \emptyset.
\end{eqnarray}
\end{proposition}


\begin{proof}

Since $K$ is spacelike and short, there is a smooth function $t_K :M\to \R$ such that $K = \{(t_K(x),x)\mid x\in K'\}$ for some smooth curve $K'\subset M$. Furthermore we can assume without loss of generality that $K$ is short enough so that
\begin{eqnarray}
\label{K contained in U}
K' \subset \pi(U')
\end{eqnarray}
where $U'$ is as in Lemma \ref{uniform unique minimizer}. 

By assumption $\Lambda \subset T^*\M$ is a conic Lagrangian which contains the element $\lambda_0 = (R, x_0, -dt + \xi_0')$. The set of points $\lambda\in\Lambda$ for which $\pi_{\Lambda} := \pi\mid_{\Lambda}$ has constant rank in a neighborhood of $\lambda$ is open and dense. So there is a sequence 
\begin{eqnarray}
\label{laj}
\lambda_j = (R_j, x_j, -dt + \xi'_j)\in \Lambda,\ \xi'_j\in S^*_{x_j}M
\end{eqnarray}
converging to $\lambda_0$ with $R_j >R$ and $x_j\neq x_0$
such that $\pi_{\Lambda}$ has constant rank in a neighborhood of $\lambda_j$. 
If $j\in \nn$ is sufficiently large, 
\begin{eqnarray}
\label{future flowout of lambdaj doesn't intersect}
\pi\circ\FLO_g^+(\lambda_j)\cap [R,R+\delta_0/2] \times \{x_0\} = \emptyset
\end{eqnarray}
where $\delta_0$ is the radius of injectivity defined in \eqref{radius of inj delta0}.

As $\Lambda$ is the flowout of a submanifold transverse to the vector field $H$, for each $\lambda_j$ there is a unique $(s_j, r_j, y_j, -dt + \eta_j')$, which we denote by $(s_j,\eta_j)$, so that 
$$e_+(s_j, \eta_j) = \lambda_j.$$
Using the constant rank theorem there are conic open sets $V_j \subset T^*\M$ containing $\lambda_j$ such that $\pi_{\Lambda}(\Lambda \cap V_j) = S_j$ is a smooth codimension $k$ submanifold of $\M$ and that $\Lambda \cap V_j$ is a conic open subset of $N^*S_j$ (see remark on p.83 of \cite{Duistermaat:1996aa}). 

By observing the the fibers of $\Lambda$ are lightlike covectors, we can deduce the codimension $k$. Indeed, suppose $k\geq 2$ and let $(t,x)\in S_j$. By the fact that $\Lambda \cap V_j$ is a conic open subset of $N^*S_j$, $\pi^{-1}((t,x)) \cap \Lambda$ is an open subset of the vector space $N_{(t,x)}^*S_j$ and therefore contains an open convex subset. Since 
the vector space is assumed to have dimension greater than two, this open convex subset contains two linearly independent covectors $\xi, \tilde \xi \in \pi^{-1}((t,x))\cap \Lambda$ such that $(1-a) \xi + a\tilde\xi \in\pi^{-1}((t,x))\cap \Lambda$ for all $a\in (0,1)$. But nontrivial convex combination of two linearly independent lightlike covectors cannot be lightlike. So our assumption that $k\geq 2$ leads to a contradiction. Since $k\geq1$ by the fact that $\Lambda$ is conic, we conclude that $k = 1$. In other words,
\begin{eqnarray}
\label{dim of Sj}
{\rm Dim}(S_j)  = {\rm rank}(D\pi_{\Lambda}(\lambda)) = 3,  \forall \lambda \in V_j\cap \Lambda.
\end{eqnarray}
Since $S_j$ is lightlike we may assume that it is contained in a sufficiently small neighborhood of $(R_j, x_j)$ so that
\begin{eqnarray}
\label{Sj is a graph}
S_j = \{(t_{S_j}(x),x)\mid x \in S_j'\}
\end{eqnarray}
for some open subset $S_j'\subset M$ and smooth function $t_{S_j}(\cdot) : S_j' \to \R$.

Recall that the quotient of $\Lambda$ and $N^*K\cap \Char_g$ by $\R^+$ is identified with \eqref{quotient conormal} and \eqref{quotient lagrangian}. The projection $\pi_{\Lambda}(\lambda) = \pi_{\Lambda/\R^+}(\lambda)$ so we may assume that $V_j$ is chosen sufficiently small so that 
\begin{eqnarray}
\label{proj is diffeo}
\pi_{\Lambda/\R^+} : (\Lambda\cap V_j)/\R^+ \to S_j\ \mbox{is a diffeomorphism}.
\end{eqnarray}
So for each $\lambda_j$ of the form \eqref{laj} there is $\delta_j>0$ and $U_j \subset (N^*K \cap \Char_g)/\R^+$ open, 
$$ (s_j, r_j, y_j, -dt+\eta'_j)\in (s_j-\delta_j, s_j+\delta_j)\times U_j $$ 
such that the map $(s,\eta)\mapsto \pi\circ e_+(s,\eta)$ is a diffeomorphism from
$$(s_j-\delta_j, s_j+\delta_j)\times U_j  \to S_j.$$
By \eqref{K contained in U}, we may assume that $U_j$ is chosen sufficiently small so that
\begin{eqnarray}
\label{Uj contained in U}
U_j \subset \{(t_K(y),y,-dt+ \eta')\mid (y, \eta')\in U'\}.
\end{eqnarray}
By shrinking $S_j$ if necessary we may assume that it is the image of $(s_j-\delta_j, s_j+\delta_j)\times U_j$ under $\pi\circ e_+$.

We now need to show that there exists an open set $\mathcal O_j \subset \M$ containing $(R_j,x_j)$ and open segment $\Gamma_j\subset K$ containing $(r_j,y_j)$ such that 
$$T^*\mathcal O_j\cap  \FLO_g^+(N^*\Gamma_j\cap\Char_g) \subset N^*S_j.$$ We will do this by first showing that

\begin{lemma}
\label{above and below}
For each fixed $j\in\nn$ sufficiently large, there exists an open set $\oo_j\subset \M$ containing $(R_j, x_j)$ and an open segment $\Gamma_j\subset K$ containing $(r_j, y_j)$ such that 
\begin{eqnarray}
\oo_j \cap  \pi(\Lambda_j)=  \oo_j\cap S_j.
\end{eqnarray}
Here, $\Lambda_j := \FLO_g^+(N^*\Gamma_j\cap\Char_g) $.

\end{lemma}
 
\begin{proof}[Proof of Lemma \ref{above and below}]
The map $(s,\eta)\mapsto e_+(s,\eta)$ is a diffeomorphism from 
$$\R^+\times ( N^*K\cap \Char_g )\to \Lambda.$$
So for each $\lambda_j$ of the form \eqref{laj} there exists a unique $s_j\in\R$ and $\eta_j =(r_j, y_j, -dt+ \eta'_j) \in \left(N^*K\cap\Char_g\right)/\R^+$ so that $\lambda_j = e_+(s_j, \eta_j)$. We also have that $s_j\to R$ and $\eta_j \to (0,y_0, -dt + \eta'_0)$. In particular, $x_j = \gamma_g(s_j, \eta'_j)$.
Furthermore, for $j\in\nn$ sufficiently large, $(y_j, \eta'_j)\in U'\subset S^*M$ where $U'$ is as in Lemma \ref{uniform unique minimizer}. 

The geodesic segment $\gamma_g([0,R], \eta'_0)$ is the unique minimizer between end points without conjugate points. Therefore, since $(y_j, \eta'_j)\in U'\subset S^*M$ where $U'$ is as in Lemma \ref{uniform unique minimizer}, the unit speed Riemannian geodesic segment $\gamma_g([0,s_j], \eta'_j)$ satisfies the following condition
\begin{eqnarray}
\label{minimizer on the conormal flowout}
\gamma_g([0,s_j], \eta'_j)\ \mbox{is the unique minimizer between its endpoints }\ x_j\ {\rm and}\ y_j.
\end{eqnarray}
We will show that \eqref{minimizer on the conormal flowout} is contradicted if the lemma fails to hold. To this end, suppose the lemma fails for some $j\in\nn$ large enough so that $(y_j, \eta'_j)\in U'$. Then with $j\in \nn$ fixed there exists a sequence of points $\{(r_{j,l}, y_{j,l})\}_{l\in \nn}\subset K$ converging to $(r_j, y_j)$, a sequence of future-pointing lightlike covectors $\eta_{j,l} \in (N_{(r_{j,l}, y_{j,l})}^*K\cap\Char_g)/\R^+$, and $s_{j,l}$ so that $\pi\circ e_+(s_{j,l}, \eta_{j,l}) \to (R_j, x_j)$ but $\pi\circ e_+(s_{j,l}, \eta_{j,l})\notin S_j$. Furthermore since 
$$\pi\circ e_+(\cdot,\cdot) : (s_j-\delta_j, s_j+\delta_j)\times U_j  \to S_j$$
is a diffeomorphism, we have that 
\begin{eqnarray}
\label{sequence not in Uj}
(s_{j,l}, \eta_{j,l})\notin (s_j-\delta_j, s_j+\delta_j)\times U_j.
\end{eqnarray}

Using the explicit expression for the flow \eqref{PhiH explicitly} we see that $s_{j,l} \to s_j$. This and \eqref{sequence not in Uj} forces $\eta_{j,l} \notin U_j$ for all $l$ sufficiently large. Using \eqref{PhiH explicitly} again we have that 
$$\gamma_g(s_{j,l}, \eta'_{j,l}) \to \gamma_g(s_j, \hat \eta'_j) = x_j$$
 for some $g$-unit covector $\hat \eta'_j \neq \eta'_j$. This contradicts \eqref{minimizer on the conormal flowout}. Thus the lemma is established.
\end{proof}
Now that Lemma \ref{above and below} is established, we only need to show that, by taking $\Gamma_j$ shorter if necessary, for all $(t,x)\in\oo_j \cap \pi(\Lambda_j)$, the fibers of 
$$\Lambda_j = \FLO_g^+(N^*\Gamma_j \cap\Char_g)$$ 
over $(t,x)$ satisfies 
\begin{eqnarray}
\label{lambdaj is subset}
T^*_{(t,x)}\M \cap \Lambda_j \subset N^*_{(t,x)} S_j.
\end{eqnarray}
Suppose \eqref{lambdaj is subset} is false. Then there exists a $(\hat t,\hat x)\in\oo_j \cap \pi(\Lambda_j)$ and $$(\hat r, \hat y, \hat \eta) \in( N^*\Gamma_j \cap \Char_g)/\R^+$$ such that 
$$(\hat t,\hat x) \in \pi\circ\FLO_g^+(\hat r, \hat y, \hat \eta)$$
but 
$$\FLO_g^+(\hat r, \hat y, \hat \eta) \cap T^*_{(\hat t,\hat x)}\M \not\subset  N^*_{(\hat t,\hat x)} S_j .$$

Lemma \ref{unique lightlike vector of lightlike surface} then forces $\pi\circ\FLO_g^+(\hat r, \hat y, \hat \eta)$ to intersect $S_j$ transversely at $(\hat t, \hat x)$. This would contradict Lemma \ref{above and below}

Finally, \eqref{flowout doesn't hit (t,x0)} is a consequence \eqref{future flowout of lambdaj doesn't intersect} if we choose $S_j$ and $\mathcal O_j$ sufficiently small.
\end{proof}

\subsection{Threefold Interaction Producing Conormal Waves} \label{3fold}
In this section we will use Proposition \ref{flowout lagrangian is almost conormal} to deduce certain conormal properties produced by interacting waves. 

To this end consider distinct points $x_0, x_1, x_2\in \mathcal U$ and $y_0\in M \backslash\U$ and assume that each $x_l$ is joined to $y_0$ by the unit speed geodesic segment $\gamma_g([0,R_l], \xi'_l)$ which is the unique minimizer between the end points with no conjugate points along the segment. Observe that as a consequence, 
\begin{eqnarray}
\label{gets further after R0}
d_g(x_0, \gamma_g(R_0+t, \xi'_0)) > R_0-t
\end{eqnarray}
for all $t>0$.

 Furthermore we assume that
\begin{eqnarray}
\label{linear independent at the target}
\dot \gamma_g(R_l,\xi_l') \neq \pm \dot \gamma_g(R_k,\xi'_k)
\end{eqnarray}
when $l\neq k$.
We label the opposite of the arrival direction at $y_0$ of the geodesic segments $\gamma_g([0, R_l],\xi'_l)$ to be 
\begin{eqnarray}
\label{def of etal}
\eta'_l := -\dot \gamma_g(R_l, \xi'_l)^\flat.
\end{eqnarray}
Due to \eqref{linear independent at the target}, any two of $\{\eta'_k, \eta'_l\}$ is linearly independent. We assume in addition that
\begin{eqnarray}
\label{span is only 2}
{\rm Dim}\spn\{\eta'_0, \eta'_1, \eta'_2\} = 2
\end{eqnarray}
Furthermore we assume that \begin{eqnarray}
\label{full span}
{\rm Dim}\left(\spn\{-dt - \eta'_0, -dt - \eta'_1, -dt - \eta'_2\} \right)= 3.
\end{eqnarray}
Note that as a consequence of \eqref{span is only 2} and \eqref{full span} we have that 
$$dt \in \spn\{-dt - \eta'_0, -dt - \eta'_1, -dt - \eta'_2\}$$ 
which then implies
\begin{eqnarray}
\label{-dt+eta'0 is in span}
-dt + \eta'_0 \in \spn\{-dt - \eta'_0, -dt - \eta'_1, -dt - \eta'_2\}.
\end{eqnarray}

Let $t_0 = 0$ and $t_1, t_2\in\R$ satisfy $R_0 = R_1+t_1 = R_2 + t_2$. Denote by 
$$\xi_l = -dt + \xi'_l \in T^*_{(t_l, x_l)}\M$$ 
and observe that
\begin{eqnarray}
\label{R0y0 in flowout}
(R_0,y_0) \in \bigcap_{l=0}^2 \pi\circ\FLO_g^+(t_l, x_l, \xi_l).
\end{eqnarray} 
Also, by assumption \eqref{linear independent at the target} the point $x_l$ is not on the minimizing geodesic segment $\gamma_g([0,R_l], \xi'_l)$ if $l\neq k$. So 
\begin{eqnarray}
\label{doesn't get in the way}
(t_k,x_k)\notin I_g^+(t_l, x_l).
\end{eqnarray}
Combining the fact that $\gamma_g([0,R_0], \xi'_0)$ is the unique minimizer between the end points with \eqref{linear independent at the target} we get that
\begin{eqnarray}
\label{doesn't circle back, threefold interaction 1}
(2R_0, x_0)\notin \pi\circ\FLO_g^+(t_l,x_l,\xi_l)
\end{eqnarray}
for $l = 0,1,2$.

\begin{proposition}
\label{main 3fold interaction thm}
Assume that \eqref{R0y0 in flowout}, \eqref{linear independent at the target}, \eqref{doesn't get in the way}, \eqref{span is only 2}, and \eqref{full span} are satisfied. For each $l=0,1,2$ there is a sequence of $\xi_{l;j} \in T^*_{(t_l, x_l)}\M \cap \Char_g$ converging to $\xi_l$ satisfying the following properties:

For each element $ \xi_{l;j}$ we may choose an $h_j>0$ so that if $h\in (0,h_j)$ and $f_{l;j}$ is a ($h$-dependent) distribution of the form \eqref{source} with 
$$\xi_{l;j}\in \WF(f_{l;j})\subset \ccl \B_h( \xi_{l;j})\cup \ccl\B_h(-\xi_{l;j}) \subset\subset \ccl \B_{h_0}( \xi_{l})\cup \ccl\B_{h_0}(- \xi_{l}),$$ 
$$\sigma(f_{l;j})(t_l, x_l, \xi_{l;j})\neq 0,$$
and $v^j$ are solutions of \eqref{cubic wave} with source $f^j :=\sum_{l=0}^2 \epsilon_l f_{l;j}$ then $\singsupp (v_{012}^j)$ contains a point $(\tilde T_j, \tilde x_j)$ with $\tilde T_j > 2R_0$ and $\tilde x_j\neq x_0$. The sequence $\{(\tilde T_j, \tilde x_j)\}_{j\in \nn}$ converges to $(2R_0, x_0)$. 

Furthermore for each $j\in\nn$ sufficiently large, there are open sets $\tilde{\mathcal O}_j$ containing $(\tilde T_j, \tilde x_j)$ with $\tilde{\mathcal O}_j \cap\left(\R\times \{x_0\}\right) = \emptyset$ such that for all $h\in (0,h_j)$,
\begin{eqnarray}
\label{threefold interaction symbol}
v_{012}^j\mid_{\tilde{\mathcal O}_j}\in I(\M, N^*\tilde S_j), \ \sigma(v_{012}^j)(\tilde T_j,\tilde x_j,\tilde \xi)\neq 0,\ {\rm for}\ \tilde\xi\in N_{(\tilde T_j, \tilde x_j)}^*\tilde S_j
\end{eqnarray}
for some lightlike hypersurfaces $\tilde S_j$ which satisfy \eqref{flowout doesn't hit (t,x0)} with $2R_0$ in place of $R$.
\end{proposition}

%
%
%
%
%
%
%
%
%
%
%
%
Due to the absence of conjugate points along $\gamma_{g}([0, R_l],\xi'_l)$, Lemma \ref{conormal when not conjugate} asserts that there is a $\delta_{y_0}>0$ such that for all $h_0>0$ sufficiently small 
$$\FLO_{g}^+(\ccl(\B_{h_0}(\xi_l)\cup\B_{h_0}(-\xi_l))\cap\Char_g)\cap T^*B_G(R_0, y_0; \delta_{y_0}) = N^*\hat S_l$$
for some lightlike hypersurface $\hat S_l$. Due to \eqref{full span} and Lemma \ref{conormal intersection}, the triple intersection 
\begin{eqnarray}\label{triple intersection in Omega}
(R_0,y_0)\in K := B_G(R_0, y_0; \delta_{y_0}) \cap\bigcap_{l=0}^2 \pi\circ\FLO_g^+(\ccl\B_{h_0}(\xi_l)\cap\Char_g) 
\end{eqnarray}
is a spacelike curve.
Lemma \ref{conormal intersection} also states that 
$$N_{(R_{0}, y_0)}^*K = \spn\{-dt - \eta'_{0}, -dt - \eta'_{1}, -dt - \eta'_{2}\}.$$
Due to \eqref{-dt+eta'0 is in span}, $-dt + \eta'_{0}\in N_{(R_{0}, y_0)}^*K\cap\Char_g$. Therefore inserting the definition of $\eta'_{0}$ (see \eqref{def of etal}) into formula \eqref{PhiH explicitly} we see that 
$$(2R_0, x_0, \pm dt \pm \xi'_{0}) \in \Lambda := \FLO_g^+(N^*K\cap\Char_g).$$
Due to \eqref{doesn't circle back, threefold interaction 1} there is a $\delta>0$ such that for all $h_0>0$ sufficiently small,
\begin{eqnarray}
\label{doesn't circle back, threefold interaction}
B_G(2R_0, x_0;\delta)\cap \pi\circ\FLO_g^+(\ccl\B_{h_0}(\xi_l)\cap\Char_g) = \emptyset
\end{eqnarray}
for $l = 0,1,2$.

The geodesic segments $\gamma_g([0,R_{l}], \eta'_{l})$ are the unique minimizers from $y_0$ to $x_l$ and contains no conjugate points. Therefore, $y_0$ is the first point along $\gamma_g([0,R_{0}], \xi'_{0})$ that intersects $\gamma_g([0,R_{l}], \xi'_{l})$ for $l = 1,2$. Furthermore, uniqueness of the minimizers $\gamma_g([0,R_{0}], \eta'_{0})$ and $\gamma_g([0,R_{l}], \eta'_{l})$ ensures that even if we can extend both $\gamma_g([0,R_{0}], \eta'_{0})$ and $\gamma_g([0,R_{l}], \xi'_{l})$ slightly, $y_0$ is still the only point of intersection. Observe also that due to \eqref{gets further after R0},
$$I_g^-(2R_0, x_0) \cap \pi\circ \FLO_g^+(t_0,x_0,\xi_0) \subset \{(t, \gamma_g(t,\xi'_0))\mid t\in [0, R_0]\}$$
(recall that $t_0  =0$).

Therefore, we can choose $\delta>0$ small enough so that \eqref{doesn't circle back, threefold interaction} is satisfied and in addition
\begin{eqnarray}
I^-_g(B_G(2R_0,x_0;\delta)) \cap\bigcap_{l=0}^2\pi\circ \FLO_g^+((t_l,x_l,\xi_{l})) = (R_{0}, y_0).
\end{eqnarray}
If $h_0>0$ is chosen small enough we can conclude from \eqref{triple intersection in Omega} that
{\begin{eqnarray}
\label{K is only intersection in past cone}
I^-_g(B_G(2R_0, x_0;\delta)) \cap\bigcap_{l=0}^2\pi\circ \FLO_g^+(\ccl\B_{h_0}(\xi_{l})\cap\Char_g)= K.
\end{eqnarray}}

We now evoke Proposition \ref{flowout lagrangian is almost conormal} to produce a sequence of elements $( \tilde R_j, \tilde y_j) \in K$ converging to $(R_{0}, y_0)$ and elements $(\tilde T_j, \tilde x_j, -dt+\tilde\xi'_j)\in T^*B_G(2R_0,x_0;\delta)\cap\Char_g$ converging to $(2R_0, x_0, -dt -\xi'_{0})$ with $\tilde T_j> 2R_0$ and $\tilde x_j\neq x_0$ such that
\begin{eqnarray}
\label{generic point is in flowout}
(\tilde T_j, \tilde x_j,  -dt+\tilde\xi'_j) \in \FLO_g^+(N^*_{(\tilde R_j, \tilde y_j)}K \cap \Char_g).
\end{eqnarray}
By Lemma \ref{uniform unique minimizer}, within $N^*_{(\tilde R_j, \tilde y_j)}K$ there is a unique lightlike covector 
$$-dt+ \tilde \eta'_j\in N^*_{(\tilde R_j, \tilde y_j)}K\cap\Char_g$$ 
and unique $\tilde s_j>0$ so that 
\begin{eqnarray}
\label{def of -dt + tildeeta'}
(\tilde T_j, \tilde x_j, \mp dt\pm\tilde\xi'_j)= e_+(\tilde s_j, \tilde R_j, \tilde y_j, \mp dt \pm \tilde\eta'_j)
\end{eqnarray}
with 
\begin{eqnarray}
\label{tildeetaj converges}
(\tilde s_j, \tilde R_j, \tilde y_j, - dt + \tilde\eta'_j) \to (R_0, R_0, y_0, -dt + \eta'_0).
\end{eqnarray}
Since $(\tilde R_j, \tilde y_j) \in K$ and for $l=0,1,2$ $\gamma_g([0,R_l],\xi'_l)$ are unique minimizing segments containing no conjugate pints, for each $j\in \nn$ there exists a unique (modulo multiplication by $\R$), $\xi_{l;j}\in \B_{h_0}(\xi_{l})\cap \Char_g$ such that $(\tilde R_j, \tilde y_j) \in \bigcap_{l=0}^2\pi\circ\FLO_g^+(t_l,x_l,\xi_{l;j})$. We may assume that $\xi_{l;j}$ is of the form
$$\xi_{l;j} = -dt + \xi'_{l;j},\ \xi'_{l;j}\in S^*_{x_l}M.$$
For each fixed $j\in \nn$ there is an $h_j<h_0$ such that all $h\in (0, h_j)$ we have that $\B_h(\pm\xi_{l;j}) \subset \B_{h_0}(\pm\xi_{l})$. The triple intersection
\begin{eqnarray}
\label{tilde triple intersection}
\Gamma_{h;j} := \bigcap_{l=0}^2\pi\circ\FLO_g^+(\ccl\B_h(\xi_{l;j})\cap\Char_g) \cap B_G(R_0, y_0; \delta_{y_0})\subset K
\end{eqnarray}
is transverse and that $\bigcap_{h>0} \Gamma_{h;j} = (\tilde R_j, \tilde y_j)$. 
Furthermore, if we define $ \tilde\eta_{l;j}\in T_{(\tilde R_j, \tilde y_j)}^*\M$ by
\begin{eqnarray}
\label{def of tilde etaj}
(\tilde R_j, \tilde y_j, \tilde\eta_{l;j}):= e_+(\tilde R_j - t_l,t_l, x_l, \xi_{l;j}) 
\end{eqnarray}
then by Lemma \ref{conormal intersection}
$$N^*_{(\tilde R_j, \tilde y_j)}\Gamma_h = \spn\{\tilde\eta_{0;j}, \tilde\eta_{1;j},\tilde\eta_{2;j}\}$$
because the triple intersection \eqref{tilde triple intersection} is transverse.
Observe that with $\eta'_l$, $l=0,1,2$, defined as in \eqref{def of etal}, we have that
\begin{eqnarray}
\label{tildeetalj converges to}
(\tilde R_j, \tilde y_j, \tilde \eta_{l;j}) \to (R_0, y_0, -dt- \eta'_l).
\end{eqnarray}
In particular, $-dt+\tilde\eta'_j\in N^*_{(\tilde R_j, \tilde y_j)}\Gamma_h$ defined by the relation \eqref{def of -dt + tildeeta'} obeys 
\begin{eqnarray}
\label{-dt+tildeeta in span}
-dt+\tilde\eta'_j\in \spn\{\tilde\eta_{0;j}, \tilde\eta_{1;j},\tilde\eta_{2;j}\}
\end{eqnarray}
and
\begin{eqnarray}
\label{-dt+tildeeta not in span of two}
 -dt+\tilde\eta'_j\notin\spn\{\tilde\eta_{l,j},\tilde\eta_{k,j}\}\ \  k,l=0,1,2.
\end{eqnarray}
To conclude \eqref{-dt+tildeeta not in span of two} we first observe that due to \eqref{linear independent at the target}, $-dt +\eta'_0 \notin \spn\{ -dt - \eta'_l\}$ for $l= 0,1,2$ (see definition \eqref{def of etal}). Therefore due to \eqref{tildeetaj converges} and \eqref{tildeetalj converges to}, for $j\in \nn$ sufficiently large, $-dt +\tilde \eta'_j \not\in\spn\{ \tilde \eta_{l,j}\}$ for $l=0,1,2$. Finally we use Lemma \ref{only threefold propagate}, while observing that $-dt+\tilde\eta'_j$, $\tilde\eta_{l;j}$, and $\tilde \eta_{k;j}$ are all lightlike covectors, to conclude \eqref{-dt+tildeeta not in span of two}.

\begin{lemma}
\label{-dt+tildeeta not in 2span}
For each $j\in \nn$ sufficiently large there exists an open conic neighourbhoods $P_j\subset T^*\M$ containing $\spn\{-dt+\tilde\eta'_j\} \subset N^*_{(\tilde R_j, \tilde y_j)}\Gamma_{h;j}$ such that
$$P_j \cap\FLO_g^+(\ccl\B_h(\pm\xi_{l;j})\cap\Char_g ) = \emptyset$$
and {\small
$$P_j\cap\left( \FLO_g^+(\ccl(\B_h(\xi_{l;j})\cup\B_h(-\xi_{l;j}) )\cap\Char_g ) + \FLO_g^+(\ccl(\B_h(\xi_{k;j})\cup\B_h(-\xi_{k;j}) )\cap\Char_g )\right)= \emptyset $$}
for any $k,l = 0,1,2$ and $h>0$ sufficiently small.
\end{lemma}
\begin{proof}
We identify $T^*\M/\R^+ \cong \M \times S^3$ locally and for each $\eta \in T^*\M$ denote by $\eta/\R^+$ the corresponding element in the quotient space. We first use \eqref{-dt+tildeeta in span} and \eqref{-dt+tildeeta not in span of two} to find open sets $N_j\subset S^3$ containing both $\pm (-dt +\tilde\eta'_j)/\R^+ \in S^3$ so that $\bar N_j\cap \spn\{\tilde\eta_{l;j},\tilde \eta_{k;j}\}/\R^+ = \emptyset$.

 Now we can choose $h>0$ sufficiently small so that if
 $$(R, y,\tilde\eta_l)\in \FLO_g^+(\ccl( \B_h(\xi_{l;j})\cup  \B_h(-\xi_{l;j}))) \cap T^* B_G(R_0, y_0; \delta_{y_0})$$ 
 and 
$$(R, y,\tilde\eta_k)\in \FLO_g^+(\ccl( \B_h(\xi_{k;j})\cup  \B_h(-\xi_{k;j}))) \cap T^* B_G(R_0, y_0; \delta_{y_0})$$ 
 then we also have $N_j\cap\spn\{\tilde \eta_k,\tilde\eta_l\}/\R^+ = \emptyset$. Recall that $ B_G(R_0, y_0; \delta_{y_0})\subset \M$ is the open set chosen so that \eqref{triple intersection in Omega} holds. We now choose $P_j =  B_G(R_0, y_0; \delta_{y_0})\times N_j\times \R^+$.
\end{proof}

Proposition \ref{flowout lagrangian is almost conormal} also states that for each $j\in \nn$ there is an open set $\mathcal O_j\subset \subset B_G(2R_0,x_0;\delta)$ containing $(\tilde T_j, \tilde x_j)$ such that if $h>0$ is sufficiently small, then
\begin{eqnarray}
\label{flowout from gamma is conormal}
\FLO_g^+(N^*\Gamma_{h;j}\cap \Char_g) \cap T^*\mathcal O_j \subset N^*(S_j\cap\mathcal O_j)
\end{eqnarray}
for some lightlike hypersurface 
$$S_j\cong (\tilde s_j -\delta_j, \tilde s_j+\delta_j)\times U_j$$ 
where the diffeomorphism is given by the map $\pi\circ e_+(\cdot,\cdot)$ and $U_j$ is an open subset of $(N^*K \cap\Char_g)/\R^+$ containing $(\tilde R_j, \tilde y_j,-dt + \tilde\eta'_j)$. Proposition \ref{flowout lagrangian is almost conormal} also states that $\oo_j$ and $S_j$ can be chosen to satisfy the flowout condition \eqref{flowout doesn't hit (t,x0)}. Without loss of generality we may assume that $U_j\subset\subset P_j/\R^+$ where $P_j\subset T^*\M$ is the open conic subset constructed in Lemma \ref{-dt+tildeeta not in 2span}.

Since $x_0\neq \tilde x_j$, we can require that $\mathcal O_j$ satisfies
\begin{eqnarray}
\label{O does not intersect vertical line}
\left(\R\times x_0 \right)\cap \mathcal O_j = \emptyset.
\end{eqnarray}

We also choose here $\tilde U_j\subset\subset U_j\subset\subset P_j/\R^+$ containing $(\tilde R_j, \tilde y_j,-dt + \tilde\eta'_j)$ and $\tilde\delta_j<\delta_j$ so that 
\begin{eqnarray}
\label{tilde S}
\tilde S_j := \pi\circ e_+((\tilde s_j -\tilde\delta_j, \tilde s_j+\tilde\delta_j)\times\tilde U_j)
\end{eqnarray}
is compactly contained in $S_j\cap \mathcal O_j$. Choose $\tilde{\mathcal O}_j\subset \subset \mathcal O_j\subset \subset B_G(2R_0,x_0;\delta)$ open and containing $(\tilde T_j, \tilde x_j)$ so that 
\begin{eqnarray}
\label{tilde O}
\tilde{\mathcal O}_j\cap S_j \subset\subset \tilde S_j
\end{eqnarray}
(see\eqref{doesn't circle back, threefold interaction} for criterion on choice of $B_G(2R_0,x_0;\delta)$).

In this setup we have that
\begin{eqnarray}
\label{tilde flowout from gamma is conormal}
\FLO_g^+(N^*\Gamma_{h;j}\cap \Char_g) \cap T^*\tilde{\mathcal O}_j \subset N^*\tilde S_j.
\end{eqnarray}

\begin{proof}[Proof of Proposition \ref{main 3fold interaction thm}]

Fix $j\in \nn$ large and for $l = 0,1,2$ construct sources of the form \eqref{source}
$$f_{l;j} \in I(\M,\ccl (\B_h( \xi_{l;j}) \cup \B_h(- \xi_{l;j})))$$ 
to be sufficiently smooth and supported in small neighborhoods of $(t_l,x_l)$ so that 
\begin{eqnarray}
\label{flj symbol}
\sigma(f_{l;j})(t_l,x_l,\pm \xi_{l;j})\neq 0.
\end{eqnarray}
%
%
%
%
%
%
%


For each $j$ we have that $v^j_{012}$ solves
$$\Box_g v^j_{012} = -6v^j_0 v^j_1 v^j_2,\ \ v^j_{012}\mid_{t<-1} = 0$$
where 
$$\Box_g v_l^j  = f_{l;j}, v_l^j\mid_{t<-1} =0.$$
So by Proposition \ref{propagation}
$$v_l^j \in I\left(\M, \ccl \left(\B_h(\xi_{l;j}) \cup  \B_h(-\xi_{l;j})\right), \FLO_g^+\left(\ccl\left(\B_h(\xi_{l;j})\cap\B_h(-\xi_{l;j}) \right)\cap\Char_g\right)\right)$$
with 
\begin{eqnarray}
\label{vl symbol}
\sigma(v^j_l)(t,x,\xi)\neq 0
\end{eqnarray}
for all $(t,x,\xi) \in \FLO_g^+(t_l, x_l, \pm\xi_{l;j})$. This means in particular
\begin{eqnarray}
\label{vl symbol at tildeetalj}
\sigma(v^j_l)(\tilde R_j,\tilde y_j,\pm\tilde \eta_{l;j})\neq 0
\end{eqnarray}
where $(\tilde R_j, \tilde y_j, \tilde \eta_{l;j})$ are as in \eqref{def of tilde etaj}.

Combining  Thm 23.2.9 of \cite{Hormander:2007aa} and Lemma \ref{only threefold propagate}, we can deduce from the nonhomogeneous linear equation for $v^j_{012}$ that
$$\WF(v^j_{012}) \cap \Char_g \subset \bigcup_{l=0}^2 \WF(v_l^j)
\cup \FLO_g^+\left(\left(\sum_{l=0}^2\WF(v_l^j)\right)\cap\Char_g\right).$$
If $j\in\nn$ is sufficiently large, $h\in (0,h_j)$ is sufficiently small, and $\tilde \oo_j \subset\subset B_G(2R_0,x_0;\delta)$ is chosen to be a sufficiently small neighborhood of $(\tilde T_j, \tilde x_j)$, for any $\chi_j\in C^\infty_c(\tilde{\mathcal O}_j)$ we can assert that
$$\WF(\chi_j v^j_{012})\subset \FLO_g^+\left(\left(\sum_{l=0}^2\WF(v_l^j)\right)\cap\Char_g\right)$$
due to \eqref{doesn't circle back, threefold interaction}. 

For $h>0$ small enough 
$$\pi^{-1} (I_g^-(B_G(2R_0,x_0;\delta))\cap\left(\sum_{l=0}^2\WF(v^j_l)\right) \subset N^*\Gamma_{h;j}$$ 
due to \eqref{tilde triple intersection} and \eqref{K is only intersection in past cone}. Furthermore due to \eqref{vl symbol at tildeetalj} we have that $\sigma(v_l^j)( \tilde R_j, \tilde y_j,\pm \tilde\eta_{l;j}) \neq 0$ where $\tilde\eta_{l;j}$ is defined by \eqref{def of tilde etaj}. Let $\omega_j(x,t,\xi)$ be a homogenous degree $0$ symbol whose support is contained in the open set $P_j$ constructed in Lemma \ref{-dt+tildeeta not in 2span}. Since $\tilde U_j\subset\subset U_j\subset\subset P_j/\R^+$, we can also ask that 
\begin{eqnarray}
\label{1-omegaj = 0 in tilde Uj}
1-\omega_j =0  
\end{eqnarray}
in an open conic set containing the closure of $\tilde U_j \times \R^+$ where $\tilde U_j$ was defined via \eqref{tilde S}. 
We have then by Lemma \ref{WF of product}
{\small \begin{eqnarray}
\label{source of vsing}
\omega_j(t,x,D)\left( v_0^jv_1^jv_2^j\right) \in I(\M,  N^*\Gamma_{h;j}),\ \sigma\left(\omega_j(t,x,D)\left( v_0^jv_1^jv_2^j\right)\right)\left(\tilde R_j, \tilde y_j, \mp dt \pm \tilde\eta'_j\right) \neq 0
\end{eqnarray}}
where $-dt + \tilde \eta'_j \in T^*_{(\tilde R_j, \tilde y_j)}\M\cap \Char_g$ is as in \eqref{-dt+tildeeta in span}.
At the same time
\begin{eqnarray}
\tilde U_j \cap \WF((1-\omega_j(t,x,D)) (v_0^jv_1^jv_2^j)) =\emptyset
\end{eqnarray}

We now write $v_{012}^j = v_{\rm reg} + v_{\rm sing}$ where
$$\Box_g v_{\rm sing} = -6\omega_j(t,x,D)\left( v^j_0v^j_1v^j_2\right),\ \ v_{\rm sing}\mid_{t<-1} = 0$$
and
$$\Box_g v_{\rm reg} = -6(1-\omega_j(t,x,D))\left( v^j_0v^j_1v^j_2\right),\ \ v_{\rm reg}\mid_{t<-1} = 0.$$

\begin{lemma}
\label{vreg is regular}
We have that for any $\chi_j \in C^\infty_c(\tilde\oo_j)$, $\chi_j v_{\rm reg} \in C^\infty_c(\M)$.
\end{lemma}

We assuming Lemma \ref{vreg is regular} for the time being and focus on the microlocal analysis of $v_{\rm sing}$. Combining the flow condition \eqref{def of -dt + tildeeta'} and the wavefront property \eqref{source of vsing} allows us to use Proposition \ref{propagation} to deduce that 
$$\sigma(v_{\rm sing})(\tilde T_j,\tilde x_j, \mp dt \pm \tilde\xi'_j) \neq 0.$$
Furthermore, $\WF(v_{\rm sing}) \cap\Char_g\subset \FLO_g^+\left(\WF\left(\omega_j(t,x,D)\left( v^j_0v_1^jv_2^j\right)\right)\cap \Char_g\right)$ so by \eqref{source of vsing} we can use \eqref{flowout from gamma is conormal} to deduce that $\WF(v_{\rm sing})\cap T^*\tilde{\mathcal O}_j \subset N^* \tilde S_j$ where $\tilde S_j$ is given in \eqref{tilde S} and $\tilde{ \mathcal O}_j\subset\subset B_G(2R_0,0;\delta)$ is chosen to satisfy \eqref{tilde O}. Therefore if $\tilde{ \mathcal O}_j\subset\subset B_G(2R_0,0;\delta)$ in \eqref{tilde O} is chosen small enough, $v^j_{012}\mid_{\tilde{\mathcal O}_j}\in I(\M, N^*\tilde S_j)$ and $\sigma(v_{012}^j)(\tilde T_j, \tilde x_j,\tilde \xi) \neq 0$ for all $\tilde\xi\in N^*_{(\tilde T_j, \tilde x_j)}\tilde S_j$. So we have verified \eqref{threefold interaction symbol} for $\tilde S_j$ satisfying \eqref{flowout doesn't hit (t,x0)}.
\end{proof}
It remains to give a 
\begin{proof}[Proof of Lemma \ref{vreg is regular}]
Let $(t,x,\pm dt + \eta') \in T^*\tilde {\mathcal  O}_j\cap \Char_g$ and assume that it is an element of $\WF(\chi_j v_{\rm reg})$. By Thm 23.2.9 \cite{Hormander:2007aa} 
\begin{eqnarray}
\label{backward flo intersects source}
\FLO_g^-((t,x,\pm dt + \eta'))\cap \WF\left(\left(1-\omega_j(t,x,D)\right)\left( v_0^jv_1^jv_2^j\right)\right)\neq \emptyset.
\end{eqnarray}
We have that by Lemma \ref{only threefold propagate}
$$\WF(v_0^jv_1^jv_2^j) \cap\Char_g \subset\left( \bigcup_{l=0}^2\WF(v_l^j) \right)
\cup \left(\left(\sum_{l=0}^2\WF(v_l^j)\right )\cap \Char_g\right). $$
Since 
$$\WF(v_l^j) \cap T^*\tilde\oo_j \subset \FLO_g^+\left(\ccl(\B_h( \xi_{l;j})\cup \B_h(-\xi_{l;j}))\cap\Char_g\right),$$ an element in $\FLO_g^-((t,x,\pm dt +\eta'))\cap\WF(v_l^j)$ would violate  \eqref{doesn't circle back, threefold interaction} and the fact that $\tilde \oo_j \subset \subset B_G(2R_0, x_0; \delta)$.

So in order for \eqref{backward flo intersects source} to hold, we must have that
$$\FLO_g^-((t,x,\pm dt + \eta'))\cap \left(\sum_{l=0}^2\WF(v_l^j)\right)\cap \supp(1-\omega_j)  \neq \emptyset.$$

However, $\tilde{\mathcal O}_j \subset\subset B_G(2R_0,x_0;\delta)$ and 
$$I^{-}_g(B_G(2R_0,x_0;\delta) )\cap \bigcap_{l=0}^2\pi( \WF(v_l^j) )= \Gamma_{h;j}$$ 
by \eqref{K is only intersection in past cone}, so we must have 
\begin{eqnarray}
\label{in supp 1-omegaj}
\FLO_g^-((t,x,\pm dt + \eta'))\cap N^*\Gamma_{h;j} \cap\Char_g\cap \supp(1-\omega_j(t,x,\xi))  \neq \emptyset.
\end{eqnarray}
This means that $(t,x, \pm dt + \eta') \in \FLO_g^+\left(N^*\Gamma_{h;j} \cap\Char_g\right)$ so by the fact that $(t,x,\pm dt + \eta')\in T^*\tilde{\mathcal O}_j$, we can apply \eqref{tilde flowout from gamma is conormal} to conclude that $(t,x, \pm dt + \eta')\in N^*\tilde S_j$. By \eqref{tilde S} we can conclude that $\FLO_g^-((t,x,\pm dt + \eta'))\cap N^*\Gamma_{h;j} \cap\Char_g \subset \tilde U_j$. Meanwhile by \eqref{1-omegaj = 0 in tilde Uj} $1-\omega_j = 0$ in $\tilde U_j$. This contradicts \eqref{in supp 1-omegaj}. So we conclude that \eqref{backward flo intersects source} must be false and therefore $\WF(\chi_j v_{\rm reg}) = \emptyset$.
\end{proof}

\subsection{Fivefold Interaction Producing Singularity along $\R\times \0$}\label{5fold}
We will analyze five-fold interaction of solutions of \eqref{cubic wave} with source of the form $\sum_{l=0}^4 \epsilon_l f_{l;j}$ where $f_{l;j}$ for $l =0,1,2$ are as in Proposition \ref{main 3fold interaction thm}. We continue the notation adopted in Proposition \ref{main 3fold interaction thm} except we now set $x_0 = \0$.

Before we construct sources $f_{3;j}$ and $f_{4;j}$ we need to analyze further the geometry of $N^*\tilde S_j$ produced by Proposition \ref{main 3fold interaction thm}. We denote by 
\begin{eqnarray}
\label{bf Tj}
{\bf T}_j := \R\times\{\tilde x_j\}
\end{eqnarray}
 and by the fact that $\tilde x_j\neq \0$, ${\bf T}_j\cap \left(\R\times\{\0\}\right) = \emptyset$. 


Since $(\tilde T_j, \tilde x_j, -dt + \tilde \xi'_j)\in N^*(\tilde S_j \cap \tilde{\mathcal O}_j)$ converges to $(2R_0, \0, -dt - \xi'_0)$ as $j\to\infty$ and \eqref{flowout doesn't hit (t,x0)} holds, there is a unique unit covector $\xi'_{\0;j}\in S^*_{\tilde x_j}M$, $\xi'_{\0;j} \neq \tilde \xi'_j$ such that \begin{eqnarray}
\label{short trip to origin}
\gamma_g(d_g(\0,\tilde x_j), \xi'_{\0;j}) = \0.
\end{eqnarray}

The elements of $N^*_{(\tilde T_j, \tilde x_j)}{\bf T}_j$ belong to ${\rm Ker}(\partial_t)$ and $N^*_{(\tilde T_j, \tilde x_j)}\tilde S_j$ is spanned by a single lightlike covector. Therefore it is easily checked by a dimension count that
\begin{eqnarray}
\label{only way going to origin}
\pm(-dt + \xi'_{\0;j}) \in T_{(\tilde T_j, \tilde x_j)}^*\M = N^*_{(\tilde T_j, \tilde x_j)}{\bf T}_j \oplus N_{(\tilde T_j, \tilde x_j)}^*\tilde S_j.
\end{eqnarray}


Denote by $\delta_{{\bf T}_j}\in{\mathcal D}'(\M)$ the distribution obtained by integrating along ${\bf T}_j$ and for $0<a<<h$ small choose $\chi_{a;j}\in C^\infty_c(\tilde{\mathcal O}_j)$
such that $\chi_{a;j} = 1$ in a neighborhood of $(\tilde T_j, \tilde x_j)$ and $\supp(\chi_{a;j})\subset\subset B_{G}(\tilde T_j, \tilde x_j; a)$. 

Clearly $\WF(\chi_{a;j}\delta_{{\bf T}_j}) \subset N^*{\bf T}_j$ is spacelike.
Choose $N$ a large positive number and define the conormal distributions $f_{3;j}$ and $f_{4;j}$ by
\begin{eqnarray}
\label{def of f4 f5}
f_{3;j} := \Box_g\left(\chi_{a;j}\langle D\rangle^{-N} \delta_{{\bf T}_j}\right)\in I(\M, N^*{\bf T}_j), 
\ f_{4;j} =\Box_g \chi_{a;j}(t,x)\in C^\infty_c(\tilde{\mathcal O}_j)
\end{eqnarray}
Note that while both depend on the parameter $a>0$, we suppress its dependence in the notation until we start manipulating this parameter in later sections. The distribution $\left(\chi_{a;j}\langle D\rangle^{-N} \delta_{{\bf T}_j}\right)$ is clearly conormal and 
\begin{eqnarray}
\label{symbol of v3 1}
\sigma\left( \chi_{a;j}\langle D\rangle^{-N} \delta_{{\bf T}_j}\right) (\tilde T_j, \tilde x_j, \xi')\neq 0
\end{eqnarray}
for all $\xi'\in N^*_{(\tilde T_j, \tilde x_j)}  {\bf T}_j \subset {\rm Ker}(\partial_t)$.
Due to the definition of \eqref{def of f4 f5}, uniqueness of the linear wave equation states that if $l=3,4$ and $v_l^j$ solves
$$\Box_g v_l^j = f_{l;j},\ v_l^j \mid_{t<-1} = 0$$
then
{\small \begin{eqnarray}
\label{v3 and v4}
v^j_4 = \chi_{a;j} \in C^\infty_c(\tilde{\mathcal O}_j),\ v^j_3 =  \chi_{a;j}\langle D\rangle^{-N} \delta_{{\bf T}_j}\in I(\M, N^*{\bf T}_j).
\end{eqnarray}}
Note that this holds for both $g= g_1$ and $g= g_2$ since $g_1 = g_2$ in $\U$ by assumption \eqref{same as eucl}.

Due to \eqref{symbol of v3 1}
\begin{eqnarray}
\label{symbol of v3}
\sigma( v_3^j) (\tilde T_j, \tilde x_j, \xi')\neq 0
\end{eqnarray}
for all $\xi'\in N^*_{(\tilde T_j, \tilde x_j)}  {\bf T}_j \subset {\rm Ker}(\partial_t)$
\begin{lemma}
\label{microlocal cutoff at the nearby source}
For each $j\in \nn$ one can find $\omega_j(t,x,\xi)$ a homogenous degree zero symbol which is identically $1$ in small conic neighborhoods of $\pm(-dt + \xi'_{\0;j} )\in T^*_{(\tilde T_j, \tilde x_j)}\M$ such that

$$\chi_{a;j}\omega_j(t,x,D) \left( v^j_3v^j_{012}\right)\in I(\M, T^*_{(\tilde T_j, \tilde x_j)}\M).$$
Furthermore, 
$$ \sigma(\omega_j(t,x,D)( v^j_3 v^j_{012}))(\tilde T_j, \tilde x_j, \pm(-dt + \xi'_{\0;j})) \neq 0$$
and 
$$(\tilde T_j, \tilde x_j,\pm( -dt + \xi'_{\0;j})) \notin\WF\left(1- \omega_j(t,x,D)\right).$$
\end{lemma}
\begin{proof}
By \eqref{threefold interaction symbol}, for $a\in (0,h)$ sufficiently small, $\chi_{a;j} v_{012}^j\in I(\M, N^*\tilde S_j)$ for some lightlike hypersurface $\tilde S_j$ satisfying \eqref{flowout doesn't hit (t,x0)} and that 
\begin{eqnarray}
\label{symbol v012j 2}
\sigma( v_{012}^j)(\tilde T_j, \tilde x_j, -dt + \tilde \xi'_j)\neq 0
\end{eqnarray}
if $-dt +\tilde \xi'_j \in N^*_{(\tilde T_j, \tilde x_j)} \tilde S_j$.
By \eqref{short trip to origin} the lightlike covector $-dt + \xi'_{\0;j}$ satisfies 
$$(\tilde T_j + d_g(\tilde x_j, \0), \0) \in \pi\circ \FLO_g^+(\tilde T_j, \tilde x_j, -dt + \xi'_{\0;j})$$
so by \eqref{flowout doesn't hit (t,x0)},
$$-dt + \xi'_{\0;j} \neq -dt + \tilde \xi'_j.$$
So we can find a homogenous degree zero symbol $\omega_j(t,x,\xi)$ which is identically $1$ in a conic neighbourhood containing $(\tilde T_j, \tilde x_j, -dt +  \xi'_{\0;j})$ but 
$$(\tilde T_j, \tilde x_j,  \pm(-dt + \tilde \xi'_j)) \notin \supp(\omega_j).$$

The set ${\bf T}_j$ defined in \eqref{bf Tj} intersects $\tilde S_j$ transversely at the point $(\tilde T_j, \tilde x_j)\in\M$ and the conormal bundle $N^*{\bf T}_j$ is spacelike. So by \eqref{v3 and v4} we can choose $\omega_j(t,x,\xi)$  so that
$$\WF(v^j_3) \cap \supp(\omega_j) = \emptyset.$$
We are now in a position to apply Lemma \ref{WF of product} to conclude that 
$$\chi_{a;j}\omega_j(t,x,D)(v_{012}^j v^j_3)\in I(\M, T^*_{(\tilde T_j, \tilde x_j)}\M).$$
Furthermore, due to \eqref{only way going to origin}, there is a unique $b_j\in\R$ and $\hat \xi'_j\in S^*_{\tilde x_j}M$ such that 
$$-dt + \xi'_{\0;j} = -dt + \tilde \xi'_j + b_j \hat \xi'_j \in T^*_{(\tilde T_j, \tilde x_j)}\M.$$
So Lemma \ref{WF of product} gives
$$\sigma(\omega_j(t,x,D)(v_{012}^j v^j_3)(\tilde T_j, \tilde x_j, -dt + \xi'_{\0;j}) = \sigma( v_{012}^j) (\tilde T_j, \tilde x_j, -dt + \tilde\xi'_j)\sigma(v^j_3)(\tilde T_j, \tilde x_j, b\hat \xi'_j) \neq 0$$
due to \eqref{symbol v012j 2} and \eqref{symbol of v3}. The same holds for the symbol evaluated at $-(-dt + \hat \xi'_{\0;j})$. \end{proof}

 The main result of this section is
\begin{proposition}
\label{fivefold interaction}
Let $j\in \nn$ be large, $h_j>0$ small. For all $h\in(0,h_j)$ and $a<<h$, let $v^j$ be the unique solution of \eqref{cubic wave} with source $\sum_{l = 0}^4 \epsilon_l f_{l;j}$. Then the distribution of one variable $ t\mapsto v^j_{01234}(\cdot; \0)$ has a singularity at $t = \tilde T_j + d_g(\tilde x_j,\0)$.
\end{proposition}
Observe that for distinct elements $k,l,m,n\in \{0,\dots,4\}$, $v_{k,l}^j = v^j_{klmn}= 0$ by Lemma \ref{expansion of sol}. So direct calculation yields that
\begin{eqnarray}
\label{five fold equation}
\Box_g v^j_{01234} = \sum_{\sigma\in S_5} v^j_{\sigma(0)\sigma(1)\sigma(2)} v^j_{\sigma(3)} v^j_{\sigma(4)},\ v^j_{01234}\mid_{t<-1} = 0
\end{eqnarray}
where $S_5$ is the symmetric group on the five letters $\{0,\dots,4\}$. Denote by $S_3\subset S_5$ to be the subgroup which maps $\{0,1,2\}$ to itself and define $v_{\rm reg}$ as the unique solution to
\begin{eqnarray}
\label{regular part of fivefold}
\Box_g v_{\rm reg} = \sum_{\sigma\in S_5\backslash S_3} v^j_{\sigma(0)\sigma(1)\sigma(2)} v^j_{\sigma(3)} v^j_{\sigma(4)},\ v_{\rm reg}\mid_{t<-1} = 0
\end{eqnarray}

\begin{lemma}
\label{not main term}
The solution of \eqref{regular part of fivefold} satisfies
$$(T_j + d_g(\tilde x_j,\0), \0) \notin \singsupp(v_{\rm reg}).$$
\end{lemma}
\begin{proof}
Consider $v^j_{3kl}$ with $k,l\in \{0,1,2\}$ which solves 
\begin{eqnarray}
\label{v3kl eq}
\Box_g v^j_{3kl} = -6 v^j_3 v^j_k v^j_l := f^j_{3kl},\ v^j_{3kl} \mid_{t<-1} = 0.
\end{eqnarray}
Since $v^j_3$ is given by \eqref{v3 and v4}, we have that 
$$\supp(f^j_{3kl})\subset \supp(v^j_3) \subset B_{G}(\tilde T_j, \tilde x_j;a) \subset\subset B_G(2R_0, \0;\delta).$$

 By Proposition \ref{propagation}, if $k \in \{0,1,2\}$,
\begin{eqnarray}
\label{vjk is flowout, fivefold interaction}
v^j_k \in I(\M, T^*_{(t_k, x_k)}\M, \FLO_g^+(\ccl\left(\B_{h_0}(\xi_{k})\cup \B_{h_0}(-\xi_{k})\right)\cap\Char_g) ).
\end{eqnarray}
So by \eqref{doesn't circle back, threefold interaction}, if $k\in \{0,1,2\}$
$\singsupp(v^j_k) \cap \supp(v^j_3) = \emptyset$ for $0<a<h$ all sufficiently small. Therefore $\WF(f^j_{3kl}) \subset \WF(v^j_3)\subset  N^*{\bf T}_j$. So by Thm 23.2.9 \cite{Hormander:2007aa},
\begin{eqnarray}
\label{wf of vj3kl} 
\WF(v^j_{3kl})  \subset N^*{\bf T}_j,\ \singsupp(v^j_{3kl}) \subset B_{G}(\tilde T_j, \tilde x_j;a)
\end{eqnarray}
for $k,l = 0,1,2$ and $0<a<h$ small. Similar argument yields 
\begin{eqnarray}
\label{wf of vj4kl}
v_{4kl}^j \in C^\infty(\M)
\end{eqnarray}
for $k, l=0,1,2$ and $0<a<h$ small. Also, observe that since 
$$\supp(v_3^j) \subset B_G(\tilde T_j, \tilde x_j;a) \subset\subset B_G(2R_0, \0;\delta)$$
we can conclude that for $k\in \{0,1,2\}$, $v_k^jv_4^j\in C^\infty(\M)$ due to \eqref{vjk is flowout, fivefold interaction} and \eqref{doesn't circle back, threefold interaction}. So $v_3^jv_4^j v_k^j \in I(\U, N^*{\bf T_j})$ with $\supp(v_3^jv_4^j v_k^j) \subset B_G(\tilde T_j, \tilde x_j; a)$. Therefore
\begin{eqnarray}
\label{wf of vj34l} 
\WF(v^j_{34l})  \subset N^*{\bf T}_j,\ \singsupp(v^j_{34l}) \subset B_{G}(\tilde T_j, \tilde x_j;a)
\end{eqnarray}
for $l = 0,1,2$ and $0<a<h$ small.
Therefore, combine \eqref{wf of vj3kl}, \eqref{wf of vj4kl}, \eqref{wf of vj34l}, and \eqref{doesn't circle back, threefold interaction}, for $\sigma'\in S_5\backslash S_3$
$$\WF(v^j_{\sigma'(0) \sigma'(1) \sigma'(2)} v^j_{\sigma'(3)} v^j_{\sigma'(4)}) \subset \WF(v^j_3) \cup \bigcup_{\sigma\in S_3}\left( \WF(v^j_{\sigma(0)}) + \WF(v^j_{\sigma(1)})\right) \cup \bigcup_{\sigma\in S_3} \WF(v^j_{\sigma(0)}).$$
So by Lemma \ref{only threefold propagate} and the fact that $\WF(v^j_3)$ is spacelike, we get that for $\sigma\in S_5\backslash S_3$
\begin{eqnarray}
\label{WF of all product'}
\WF(v^j_{\sigma(0) \sigma(1) \sigma(2)} v^j_{\sigma(3)} v^j_{\sigma(4)}) \cap \Char_g \subset \WF(v^j_0) \cup \WF(v^j_1)\cup\WF(v^j_2).
\end{eqnarray}
In view of the wavefront property given by \eqref{vjk is flowout, fivefold interaction}, \eqref{WF of all product'} becomes
{\small\begin{eqnarray}
\label{WF of all product five fold}
\WF(v^j_{\sigma(0) \sigma(1) \sigma(2)} v^j_{\sigma(3)} v^j_{\sigma(4)}) \cap \Char_g \subset \bigcup_{k=0}^2 \FLO_g^+(\ccl\left(\B_{h_0}(\xi_{k})\cup \B_{h_0}(-\xi_{k})\right)\cap\Char_g).
\end{eqnarray}}
So if $\xi\in T^*_{(\tilde T_j + d_g(\tilde x_j, \0), \0)}\M\cap \Char_g$, by the wavefront property \eqref{WF of all product five fold} and flowout property \eqref{doesn't circle back, threefold interaction}, we must have that for any $\sigma\in S_5\backslash S_3$,
$$\FLO_g^-(\tilde T_j + d_g(\tilde x_j, \0), \0,\xi)\cap \WF(v^j_{\sigma(0) \sigma(1) \sigma(2)} v^j_{\sigma(3)} v^j_{\sigma(4)}) = \emptyset.$$
By Thm 23.2.9 of \cite{Hormander:2007aa} this means $(\tilde T_j + d_g(\tilde x_j, \0), \0) \notin\singsupp(v_{\rm reg})$.\end{proof}

In view of Lemma \ref{not main term} we can write, $v_{01234}^j = v_{\rm reg} + w$ where $v_{\rm reg}$ is smooth in some neighborhood of $(T_j + d_g(\tilde x_j, \0),\0)$ and $w$ solves
\begin{eqnarray}
\label{equation for w}
\Box_g w = -6v^j_{012} v^j_3 v^j_4,\ \ w\mid_{t<-1} = 0
\end{eqnarray}

\begin{proof}[Proof of Proposition \ref{fivefold interaction}]
It suffices to show that $w$ solving \eqref{equation for w} satisfies
\begin{eqnarray}
\label{singsupp of w on line}
T_j + d_g(\tilde x_j,\0) \in \singsupp (w(\cdot,\0)).
\end{eqnarray}
Let $\omega_j(t,x,\xi)$ be the symbol constructed in Lemma \ref{microlocal cutoff at the nearby source}. We write $w = w_{\rm reg} + w_{\rm sing}$ where $w_{\rm sing}$ solves
\begin{eqnarray}
\label{w sing 1}
 \Box_g w_{\rm sing} = v_4^j \omega_j(t,x,D) \left( v_{012}^j v_3^j\right),\ w_{\rm sing}\mid_{t<-1} = 0
 \end{eqnarray}
and $w_{\rm reg}$ solves
\begin{eqnarray}
\label{w reg 1}
\Box_g w_{\rm reg} = v_4^j \left(1-\omega_j(t,x,D) \right)\left( v_{012}^j v_3^j\right),\ w_{\rm reg}\mid_{t<-1} = 0.
\end{eqnarray}
Inserting $v_3^j$ and $v_4^j$ from \eqref{v3 and v4} into the right side of \eqref{w sing 1} we get
\begin{eqnarray}
\label{w sing 2}\nonumber
 \Box_g w_{\rm sing} &=& f_{\rm source} := \chi_{a;j} \omega_j(D) \left( v_{012}^j  \left( \chi_{a;j}\langle D\rangle^{-N} \delta_{{\bf T}_j}\right)\right)\\ w_{\rm sing}\mid_{t<-1} &=& 0
 \end{eqnarray}
By Lemma \ref{microlocal cutoff at the nearby source}, $f_{\rm source}$ has wavefront set contained in the cotangent space of the point $(\tilde T_j, \tilde x_j)$ and that
$$\sigma(f_{\rm source})(\tilde T_j, \tilde x_j, \mp dt \pm \xi'_{\0;j}) \neq 0.$$
By \eqref{short trip to origin}, the element $(\tilde T_j +d_g(\tilde x_j,\0), \0, \mp dt \pm \dot \gamma_g(d_g(\tilde x_j, \0), \xi'_{\0;j})^\flat)$ belongs to the future flowout of $(\tilde T_j, \tilde x_j, \mp dt \pm \xi'_{\0;j})$. Therefore by Proposition \ref{propagation}
$$w_{\rm sing} \in I\left(\M, T^*_{(\tilde T_j, \tilde x_j)}\M, \FLO_g^+\left( T^*_{(\tilde T_j, \tilde x_j)}\M\cap\Char_g\right)\right)$$
and
\begin{eqnarray}
\label{wf at the line}
\sigma(w_{\rm sing})(\tilde T_j +d_g(\tilde x_j,\0), \0, \mp dt \pm \dot \gamma_g(d_g(\tilde x_j, \0), \xi'_{\0;j})^\flat)\neq0.
\end{eqnarray}

Since for $j\in\nn$ large $(\tilde T_j, \tilde x_j)$ is close to $(\tilde T_j +d_g(\tilde x_j,\0), \0)$, we have that there exists an open set $\hat {\mathcal O}_j$ containing $(\tilde T_j +d_g(\tilde x_j,\0), \0)$ such that 
$$\FLO_g^+\left( T^*_{(\tilde T_j, \tilde x_j)}\M\cap\Char_g\right)\cap T^*\hat{\mathcal O}_j = N^*(\hat S_j \cap \hat{\mathcal O}_j)$$ 
for some lightlike hypersurface $\hat S_j$. Therefore $w_{\rm sing}\mid_{\hat {\mathcal O}_j} \in I(\M, N^*\hat S_j)$. Due to \eqref{wf at the line} the hypothesis of Lemma \ref{restriction has singularity} are satisfied. Therefore we can conclude that $ \tilde T_j +d_g(\tilde x_j,\0)\in \singsupp(w_{\rm sing}(\cdot, \0))$.

To complete the proof we need to verify that $(\tilde T_j +d_g(\tilde x_j,\0),\0)\notin\singsupp(w_{\rm reg})$. The source term in \eqref{w reg 1} is supported in $B_{G}(\tilde T_j, \tilde x_j; a)$ so if there is an element 
\begin{eqnarray}
\label{something in the wf}
\xi\in \WF(w_{\rm reg}) \cap T^*_{(\tilde T_j + d_g(\tilde x_j, \0),\0)}\M
\end{eqnarray}
then by Thm 23.2.9 of \cite{Hormander:2007aa}
\begin{eqnarray}
\label{backward flow intersects}
\FLO_g^-((\tilde T_j + d_g(\tilde x_j, \0),\0,\xi)) \cap \WF( v_4^j \left(1-\omega_j(t,x,D) \right)\left( v_{012}^j v_3^j\right)) \neq \emptyset.
\end{eqnarray}
By \eqref{v3 and v4} we have that 
{\Small
$$ \WF\left( v_4^j \left(1-\omega_j(t,x,D) \right)\left( v_{012}^j v_3^j\right)\right)\subset \left(\WF(v^j_{012}) \cup \WF(v_3^j)\cup \left(\WF(v_3^j) + \WF(v_{012}^j)\right)\right)\cap \WF(1-\omega_j(t,x,D))$$}
and because $v_4^j\in C^\infty_c(B_G(\tilde T_j, \tilde x_j; a))$
$$\WF\left( v_4^j \left(1-\omega_j(t,x,D) \right)\left( v_{012}^j v_3^j\right)\right) \subset T^*B_{G}(\tilde T_j , \tilde x_j; a).$$
Combine this with the fact that $\WF(v_3^j)$ is spacelike, \eqref{backward flow intersects} becomes 
{\Small \begin{eqnarray}
\label{backward flow intersects 2'}\nonumber
T^*B_{G}(\tilde T_j , \tilde x_j; a)\cap \FLO_g^-((\tilde T_j + d_g(\tilde x_j, \0),\0,\xi))\cap\left(\WF(v^j_{012}) \cup \left(\WF(v_3^j) + \WF(v_{012}^j)\right)\right)\cap \WF(1-\omega_j(t,x,D)) \neq\emptyset\\
\end{eqnarray}}
Within $B_{G}(\tilde T_j , \tilde x_j; a) \subset\subset \tilde \oo_j$, $\WF\left(v_{012}^j \mid_{B_G(\tilde T_j,\tilde x_j; a)}\right) \subset N^*\tilde S_j$ by \eqref{threefold interaction symbol} and
$$\pi\circ \FLO_g^+(N^*\tilde S_j) \cap\left([2R_0, 2R_0 + \delta_0/2]\times \0)\right) =  \emptyset$$
as stated in Proposition \ref{main 3fold interaction thm}. Recall that $\tilde T_j>2R_0$ is a sequence converging to $2R_0$. Therefore, for $j\in \nn$ sufficiently large that $\tilde T_j + d_g(\tilde x_j, \0) \in[2R_0, 2R_0 + \delta_0/2]$,
$$\FLO_g^-((\tilde T_j + d_g(\tilde x_j, \0),\0,\xi)) \cap \WF(v_{012}^j) \cap T^*B_{G}(\tilde T_j , \tilde x_j; a)= \emptyset.$$
Therefore, \eqref{backward flow intersects 2'} becomes
\begin{eqnarray}
\label{backward flow intersects 2}
 \FLO_g^-((\tilde T_j + d_g(\tilde x_j, \0),\0,\xi))\cap \left(\WF(v_3^j) + \WF(v_{012}^j)\right)\cap \WF(1-\omega_j(t,x,D))\neq \emptyset.
\end{eqnarray}
By Proposition \ref{main 3fold interaction thm}, inside $B_{G}(\tilde T_j, \tilde x_j; a)$ the singular support of $v_{012}^j$ is contained in $\tilde S_j$. Since the singular support of $v_3^j$ is contained in $B_{G}(\tilde T_j, \tilde x_j;a)\cap {\bf T}_j$, we have that 
$$\left(\WF(v_3^j) + \WF(v_{012}^j) \right)\cap T^*B_{G}((\tilde T_j, \tilde x_j),a) \subset T^*_{(\tilde  T_j ,\tilde x_j)}\M.$$
Since $\tilde x_j$ is within the injectivity radius of $\0$, the only lightlike covectors in $T^*_{(\tilde  T_j ,\tilde x_j)}\M$ whose future flowout intersects $T^*_{(\tilde T_j + d_g(\tilde x_j,\0), \0)}\M$ are $\pm (-dt + \xi'_{\0;j})$ (see \eqref{short trip to origin} for definition of $\xi'_{\0;j}$). However, $\omega_j$ as constructed by Lemma \ref{microlocal cutoff at the nearby source} satisfies 
$$(\tilde T_j, \tilde x_j, \pm(-dt +  \xi'_{\0;j}))\notin\WF(1-\omega_j(t,x,D) ).$$ These combined to contradict \eqref{backward flow intersects 2}. Therefore \eqref{something in the wf} is false and 
$$\WF(w_{\rm reg}) \cap T^*_{(\tilde T_j + d_g(\tilde x_j, \0),\0)}\M = \emptyset.$$
\end{proof}

%
%
%


\section{Recovering Geometric Data}

Throughout this section we assume that $L_{M_1,g_1,2T+1} = L_{M_2, g_2,2T+1}$. 
\subsection{Geodesics Returning to $\mathcal U$}
In this subsection we take care of geodesics originating from $\U$ hitting $\0$ after exiting $\U$. For any $r\in (0,\delta_0]$, let $\partial_+ S^*B(\0,r)$, $\partial_0 S^*B(\0,r)$, and $\partial_- S^* B(\0,r)$ denote the outward, tangential, and inward pointing unit covector respectively. We shall prove that
\begin{proposition}
 \label{no boomerang}
For each $r\in (0,\delta_0/4)$, $\xi'\in \partial_+ S^*B(\0, r)\cup \partial_0S^*B(\0,r)$, and $R\in (0,T)$, we have $\0\notin \gamma_{g_2}([0,2R], \xi')$ if and only if $\0\notin \gamma_{g_1}([0,2R], \xi')$.
 \end{proposition}


 \begin{proof}
Suppose by contradiction that there exists a $\xi'_0\in \partial_+ S^*B(\0, r)\cup \partial_0S^*B(\0,r)$ based at a point $x_0\in \partial B(\0,r)$ such that $\gamma_{g_1}(T', \xi'_0) = \0$ for some $T'\in [0, 2R]$, but 
\begin{eqnarray}
\label{origin not in g2 geodesic}
\0\notin \gamma_{g_2}([0, 2R], \xi'_0). 
\end{eqnarray}
Choose a generic small $t>0$ so that 
$$x_1 := \gamma_{g_1}(t_1, \xi'_0) = \gamma_{g_2}(t_1, \xi'_0)\in \U$$
and that the end points of $\gamma_{g_1}([t_1,T'], \xi'_0)$ and $\gamma_{g_2}([t_1,T'],\xi'_0)$ are not conjugate to each other.
Define 
$$ \xi_1 := -dt + \dot\gamma_g( t_1, \xi'_0)^\flat \in T^*_{( t_1,  x_1)}\U\cap\Char_g,$$ 
and observe that, for all $h>0$ sufficiently small, there exists an open neighborhood $\oo\subset \R\times \U$ of $(T',\0)\in \R\times \U$ such that 
\begin{eqnarray}
\label{flowout is conormal no boomerang}
\FLO_{g_1}^+(\ccl(\B_h(\xi_1) \cup\B_h(-\xi_1))\cap \Char_{g_1}) \cap T^*\oo = N^*S,\ (T',\0)\in S
\end{eqnarray}
for some lightlike hypersurface $S$. However, due to \eqref{origin not in g2 geodesic}, for $h>0$ sufficiently small,
\begin{eqnarray}
\label{origin not in flowout of g2}
(T',\0)\notin\pi\circ\FLO_{g_2}^+(\ccl(\B_h(\xi_1) \cup\B_h(-\xi_1))\cap \Char_{g_2}) 
\end{eqnarray}
Let $f\in I(\U, \ccl(\B_h(\xi_0)\cup \B_h(-\xi_0)))$ be as constructed in \eqref{source} with $\sigma(f)(t_1, x_1, \xi_1)\neq 0$. Now let $u$ and $v$ be solutions of the non-linear wave equations
\begin{align*}
\label{cubic wave}
& \Box_{g_1} u + u^3 = \Box_{g_2} v + v^3 = \epsilon f,   \\
& u = v = 0,\mbox{ on } t<-1.
\end{align*}
The distributions $u'=\partial_\epsilon u|_{\epsilon=0}$ and $v'=\partial_\epsilon u|_{\epsilon=0}$ solves the linear wave equations $$\Box_{g_1} u = \Box_{g_2} v = f.$$
By \eqref{origin not in flowout of g2} and Thm 23.2.9 of \cite{Hormander:2007aa}, 
\begin{eqnarray}
\label{T'0 notin singsupp no boomerang}
(T',\0)\notin \singsupp(v').
\end{eqnarray}
However, due to Proposition \ref{propagation} and \eqref{flowout is conormal no boomerang}, $u'\mid_{\oo} \in I(\oo, N^*S)$ with $\sigma(u')(T',\0, \xi)\neq 0$ for $\xi\in N^*_{(T',\0)}S$. So Lemma \ref{restriction has singularity} asserts that $T'\in \singsupp(u'(\cdot,\0))$ for some $T'<2T$. This fact combined with \eqref{T'0 notin singsupp no boomerang} contradicts $L_{M_1, g_1, 2T+1} = L_{M_2, g_2, 2T+1}$.
\end{proof}

\subsection{Geodesics Starting From the Origin}
The goal of this section is to prove the following
\begin{proposition}
\label{unique minimizer from origin}
Let $R\in (0,T)$ and $\xi' \in S_\0^* \U$. The geodesic segment $\gamma_{g_1}([0,R];\xi') \subset M_1$ is a minimizer between its end points if and only if the geodesic segment $\gamma_{g_2}([0,R];\xi')\subset M_2$ is a minimizer between its end points.
\end{proposition}
Clearly it suffices to prove one direction of the ``if and only if'' implication. 
We will do so by contradiction. The idea is that suppose a segment $\gamma_{g_1}([0,R],\xi'_c)$ starting from $\0$ is minimizing but $\gamma_{g_2}([0,R],\xi'_c)$ is not minimizing. Then we can construct two $g_2$-lightlike geodesics, one of them with initial condition $(0,\0, -dt + \xi'_c)$ and the other one starting at a later time from $\0$ in a different direction, such that they collide at the same point in spacetime. This translates to singularities interacting in the semilinear wave equation $\Box_{g_2} u + u^3 = f$ for appropriately chosen source. This interaction will generate singularity coming back to $\0$. Such interaction will not occur for the $\Box_{g_1} u + u^3 = f$ equation since the two corresponding $g_1$-lightlike geodesics will not collide. This will contradict the fact that the source-to-solution maps are the same.

To realize this argument we construct a set of geodesics in case Proposition \ref{unique minimizer from origin} fails:

\begin{lemma}
\label{if not minimize}
Suppose Proposition \ref{unique minimizer from origin} fails to hold. That is, 
let $\xi'_c\in S^*_\0\U$ be a covector such that $\gamma_{g_1}([0,R],\xi'_c)$ is minimizing but $\gamma_{g_2}([0,R],\xi'_c)$ is not minimizing. Then, one can find $\xi'_{\rm long}, \xi'_0\in S^*_\0\U$ and $R_{\rm long} , R_0\in \R$ with $0<R_0<R_{\rm long}$ such that\\
\begin{enumerate}
\item \label{g1 is minimizing} The segment $\gamma_{g_1}([0, R_{\rm long}], \xi'_{\rm long})$ is the unique minimizer between end points and does not contain conjugate points.
\item \label{g2 geodesics meet} The segments $\gamma_{g_2}([0, R_{\rm long}], \xi'_{\rm long})$ and $\gamma_{g_2}([0, R_0], \xi'_0)$ meet at and only at the end points. In addition, $\dot \gamma_{g_2}(R_{\rm long},\xi'_{\rm long}) \neq -\dot\gamma_{g_2}(R_0,\xi'_0)$.
\item\label{R0 is minimizing} The segment $\gamma_{g_2}([0, R_0], \xi'_0)$ is the unique minimizer between the end points and does not contain conjugate points.
\item\label{minimizing if shorter} There is a $t_1>R_{\rm long} - R_0$ such that $\gamma_{g_2}(t_1,\xi'_{\rm long}) \in \U$ and $\gamma_{g_2}([t_1, R_{\rm long}],\xi'_{\rm long})$ is the unique minimizer between the end points and does not contain conjugate points.
\end{enumerate}
\end{lemma}
\begin{proof}
Let $\xi'_c\in S^*_\0\U$ be a covector such that $\gamma_{g_1}([0,R],\xi'_c)$ is minimizing but $\gamma_{g_2}([0,R],\xi'_c)$ is not minimizing. By taking $R$ slightly smaller we may assume without loss of generality that $\gamma_{g_1}([0,R],\xi'_c)$ is the unique minimizer and does not contain conjugate points. Note that this implies
\begin{eqnarray}
\label{gamma1 doesn't come back}
\0\notin \gamma_{g_1}((0,2R], \xi'_c)
\end{eqnarray}
Set
\begin{eqnarray}
\label{cut length}
R_c := \sup\{ t>0 \mid \gamma_{g_2}([0,t], \xi'_c) \subset M_2\ {\mbox{is a minimizing segment}}\}<R.
\end{eqnarray}

Let $t_1>0$ be chosen so that $\gamma_{g_2}([0,2t_1], \xi'_c) \subset \U$. Definition \eqref{cut length} implies that $\gamma_{g_2}([t_1, R_c])$ is the unique minimizing segment between end points and that it does not contain conjugate points. This means that for any $R_c'\in (R_c, R)$ sufficiently close to $R_c$, $\gamma_{g_2}([t_1, R_c'],\xi'_c)$ is a geodesic segment not containing conjugate points and is the unique minimizing segment between end points. We can also arrange that the end points of $\gamma_{g_2}([0,R_c'], \xi'_c)$ are not conjugate to each other (though the geodesic segment is not minimizing). 

For any $R_c'\in (R_c, R)$, the segment $\gamma_{g_2}([0,R_c'], \xi'_c)$ is not minimizing. So for any $R_c'\in (R_c, R)$ there is a unit covector $\xi'_0\in S^*_\0M_2$ and $R_0'< R_c'$ such that $\xi'_0\neq \xi'_c$ and $\gamma_{g_2}([0, R_0'],\xi'_0)$ is a minimizing geodesic between $\0$ and $\gamma_{g_2}(R_c', \xi'_c)$. We have the freedom to choose any $R_c'\in (R_c, R)$ so we choose it sufficiently close to $R_c$ so that the end points of $\gamma_{g_2}([0,R_c'],\xi'_c)$ are not conjugate and
\begin{eqnarray}
\label{R'c-R0'<t1}
R'_c - R_0' < t_1/2.
\end{eqnarray}
We may also assume that the only intersection of these two segments are at the end points. Furthermore, due to \eqref{gamma1 doesn't come back} and Proposition \ref{no boomerang}, 
\begin{eqnarray}
\label{original g2 doesn't come back}
\dot\gamma_{g_2}(R_c',\xi'_c) \neq - \dot\gamma_{g_2}(R_0', \xi'_0).
\end{eqnarray}

By the fact that the end points of $\gamma_{g_2}([0,R_c'], \xi'_c)$ are not conjugate to each other, there is a sequence $\xi'_j \in S^*_{\0}\U$ converging to $\xi'_c$ and $R_j$ converging to $R_c'$ such that
such that $\gamma_{g_2}(R_j, \xi'_j) = \gamma_{g_2}(R_0' - 1/j, \xi'_0)$. 
Choose $j\in \nn$ sufficiently large and set 
\begin{eqnarray}
\label{Rlong def}
R_{\rm long} := R_j,\ \xi'_{\rm long} := \xi'_j.
\end{eqnarray}
Due to \eqref{original g2 doesn't come back}, we can infer that $\dot\gamma_{g_2}(R_{\rm long},\xi'_j) \neq - \dot\gamma_{g_2}(R_1'-1/j, \xi'_0)$ for $j\in \nn$ large enough. So condition \eqref{g2 geodesics meet} is met. The geodesic segment $\gamma_{g_2}([0, R'_0- 1/j], \xi'_0)$ is the unique minimizer and does not contain conjugate point since it can be extended slightly and still be a minimizer. So condition \eqref{R0 is minimizing} is met by setting $R_0 := R'_0-1/j$. From \eqref{R'c-R0'<t1}, if $j\in \nn$ is chosen large enough
\begin{eqnarray}
\label{Rlong- R0 <t1}
R_{\rm long} - R_0 < t_1.
\end{eqnarray}

Recall that the geodesic segment $\gamma_{g_2}([t_1, R_c'], \xi'_c)$ does not contain conjugate points and is the unique minimizing segment between end points. So by Lemma \ref{uniform unique minimizer}, if $j\in\nn$ is chosen sufficiently large in definition \eqref{Rlong def}, $\gamma_{g_2}([t_1, R_{\rm long}], \xi'_{\rm long})$ is still the unique minimizer between its end points and does not contain conjugate points. So condition \eqref{minimizing if shorter} is met.
For $j\in\nn$ chosen large enough Lemma \ref{uniform unique minimizer} asserts that $\gamma_{g_1}([0,R], \xi'_{\rm long})$ is still the unique minimizer between end points and does not contain conjugate points. Since $R_{\rm long} <R$ if $j\in\nn$ is sufficiently large, $\gamma_{g_1}([0,R_{\rm long}], \xi'_{\rm long})$ is still the unique minimizer between end points and does not contain conjugate points. So condition \eqref{g1 is minimizing} is met.
\end{proof}
Observe that as a consequence of Condition \eqref{R0 is minimizing} of Lemma \ref{if not minimize} we have that $\0\notin \gamma_{g_2}((0,2R_0], \xi'_0)$. Using Proposition \ref{no boomerang} we have that 
\begin{eqnarray}
\label{0notingammag1}
\0\notin \gamma_{g_1}((0,2R_0], \xi'_0).
\end{eqnarray}

Set $y_0 := \gamma_{g_2}(R_{\rm long}, \xi'_{\rm long}) = \gamma_{g_2}(R_0, \xi'_0)$. Furthermore set $x_0 = \0$,
\begin{eqnarray}
\label{x1 def geodesic from origin}
x_1 := \gamma_{g_2}(t_1, \xi'_{\rm long}),
\end{eqnarray}
and 
\begin{eqnarray}
\label{xi'1 def geodesic from origin}
\xi'_1 := \dot \gamma_{g_2}(t_1, \xi'_{\rm long})^\flat \in S^*_{x_1} \U.
\end{eqnarray}

 Note that
$$\gamma_{g_2}([0,R_{\rm long}], \xi'_{\rm long}) =\gamma_{g_2}([0,t_1], \xi'_{\rm long})\cup \gamma_{g_2}([0,R_1], \xi'_1),$$
where
\begin{eqnarray}
\label{R1 = Rlong -t1}
R_1 := R_{\rm long} - t_1.
\end{eqnarray}
Condition \eqref{minimizing if shorter} means that the segment $\gamma_{g_2}([0,R_1], \xi'_1)$ is the unique minimizer and does not contain conjugate points.

By assumption \eqref{same as eucl}, $g_1=g_2$ in $\mathcal U$. Therefore the segment $\gamma_{g_1}([0,R_{\rm long}], \xi'_{\rm long})$ can be written as
\begin{equation}
\label{concatenation}
\begin{split}
\gamma_{g_1}([0,R_{\rm long}], \xi'_{\rm long}) 
& = \gamma_{g_1}([0,t_1], \xi'_{\rm long})\cup \gamma_{g_1}([0, R_1], \xi'_1) \\
& = \gamma_{g_2}([0,t_1], \xi'_{\rm long})\cup \gamma_{g_1}([0, R_1], \xi'_{1}).
\end{split}
\end{equation}

Set $\eta'_0 := -\dot \gamma_{g_2}(R_0, \xi'_0)^\flat$ and $\eta'_1 := -\dot \gamma_{g_2}(R_1, \xi'_1)^\flat$, both belonging to $S^*_{y_0}M_2$. By condition \eqref{g2 geodesics meet} of Lemma \ref{if not minimize}, $\eta_0' \neq -\eta_1'$. For any $\eta'_2\in S^*_{y_0}M_2 \cap\spn\{\eta'_0, \eta'_1\}\backslash\{\eta'_0\}$  sufficiently close to $\eta'_0$, Lemma \ref{uniform unique minimizer} guarantees $\gamma_{g_2}([0, R_0], \eta'_2)$ is the unique minimizing geodesic between the end points and does not contain conjugate points. We set 
\begin{eqnarray}
\label{x2 def geodesic from origin}
x_2 := \gamma_{g_2}(R_0, \eta'_2)
\end{eqnarray}
and 
\begin{eqnarray}
\label{xi'2 def geodesic from origin}
\xi'_2 := - \dot\gamma_{g_2}(R_0, \eta'_2)^\flat.
\end{eqnarray}
Observe that with this setup, if $\eta'_2$ is chosen sufficiently close to $\eta'_0$
\begin{eqnarray}
\label{no circling geod}
\dot\gamma_{g_2}(R_l, \xi'_l) \neq -\dot \gamma_{g_2}(R_k, \xi'_k),\ k,l\in\{0,1,2\}
\end{eqnarray}

For $\eta'_2$ chosen sufficiently close to $\eta'_0$ the segments $\gamma_{g_2}([0,2R_0], \xi'_0)$ and $\gamma_{g_2}([0,2R_0], \xi'_2)$ are close to each other. Condition \eqref{R0 is minimizing} of Lemma \ref{if not minimize} implies $\0\notin \gamma_{g_2}((0,2R_0], \xi'_0)$. So $\0\notin \gamma_{g_2}([0,2R_0], \xi'_2)$. Proposition \ref{no boomerang} then concludes that 
\begin{eqnarray}
\label{0 notin g1 geodesic xi2}
\0\notin \gamma_{g_1}([0,2R_0], \xi'_2).
\end{eqnarray}
The above conditions are satisfied for any $\eta'_2\in \spn\{\eta_0', \eta'_1\}$ chosen sufficiently close to $\eta'_0$. We now use Lemma \ref{full span lem} to give $\eta'_2 \in \spn\{\eta_0', \eta'_1\}$ the following additional linear algebra property:
\begin{eqnarray}
\label{e: span is 3}
{\rm Dim}\left(\spn\{-dt - \eta'_{0}, -dt - \eta'_{1}, -dt - \eta'_{2}\}\right) = 3,\ {\rm Dim}(\spn\{\eta_0', \eta_1', \eta_2'\} )= 2
\end{eqnarray}

Set $t_0$ and $t_2$ to be 
\begin{eqnarray}
\label{def of t0 and t2}
 0< t_0 = t_2 := R_{\rm long} - R_0< t_1
\end{eqnarray}
and observe that
$$(t_1, x_1)\notin I^+_{g_2}(t_0,\0)$$
by the strict inequality $t_1 = d_{g_2}(x_1,\0)> t_1-t_0$. And since $x_2$ is chosen sufficiently close to $\0$,
\begin{eqnarray}
\label{t1 is small}
(t_1, x_1)\notin I^+_{g_2}(t_0,\0)\cup I^+_{g_2}(t_2,x_2).
\end{eqnarray}

Define $\xi_l \in T^*_{(t_l, x_l)}\U$ to be $\xi_l := -dt + \xi'_l$.
In this setting 
\begin{eqnarray}
\label{three flowouts meet in M2}
(t_0+ R_0 ,y_0) \in \pi\circ\FLO^+_{g_2}(t_0,\0,\xi_0)\cap \pi\circ\FLO^+_{g_2}(t_1,x_1,\xi_1)\cap \pi\circ\FLO^+_{g_2}(t_2,x_2,\xi_2).
\end{eqnarray}
However, since $\gamma_{g_1}([0,R_{\rm long}], \xi'_{\rm long})$ is the unique length minimizing geodesic in $M_1$, we have the following
\begin{lemma}
\begin{eqnarray}
\label{does not intersect in M1}
I_{g_1}^-(t_0 + 2R_0, \0)\cap\pi\circ\FLO^+_{g_1}(t_0,\0,\xi_0)\cap \pi\circ\FLO^+_{g_1}(t_1,x_1,\xi_1) = \emptyset.
\end{eqnarray}
\end{lemma}
\begin{proof}
Note that $\gamma_{g_1}([0,R_1], \xi'_1)$ is contained in $\gamma_{g_1}([0,R_{\rm long}], \xi'_{\rm long})$ by \eqref{concatenation}. So it suffices to show that 
\begin{eqnarray}
\label{does not intersect in M1 1}
I_{g_1}^-(t_0 + 2R_0, \0)\cap\pi\circ\FLO^+_{g_1}(t_0,\0,\xi_0)\cap \pi\circ\FLO^+_{g_1}(0,\0,\xi_{\rm long}) = \emptyset
\end{eqnarray}
where $\xi_{\rm long} := -dt + \xi'_{\rm long} \in T^*_{(0,\0)}\U$.
We first observe that for $t\leq R_{\rm long}$, $\gamma_{g_1}(t, \xi'_{\rm long}) \in \partial B_{g_1}(\0,t)$ since the geodesic is uniquely minimizing up to $t = R_{\rm long}$ by condition \eqref{g1 is minimizing} of Lemma \ref{if not minimize}. Meanwhile $\gamma_{g_1}(t , \xi'_0) \in \overline{B_{g_1}(\0, t)}$. So by writing lightlike geodesics using expression \eqref{explicit expression for LL geodesics}, we see that since $t_0>0$,
{\small$$\left(\{t \leq R_{\rm long}\}\times M_1\right)\cap \pi\circ\FLO^+_{g_1}(t_0,\0,\xi_0)\cap \pi\circ\FLO^+_{g_1}(0,\0,\xi_{\rm long})  = \emptyset.$$}
So if there is an element in the LHS of \eqref{does not intersect in M1 1}, it must be in 
{\small
\begin{eqnarray}
\label{Ig1- is union of balls}
I_{g_1}^-(t_0 + 2R_0, \0)\cap\left( (R_{\rm long}, \infty)\times M_1\right) = \{(t_0 + 2R_0 -t, x) \mid d_{g_1}(\0, x) \leq t,\ t\in (0, R_0)\}.
\end{eqnarray}}
Note that by \eqref{def of t0 and t2}, $R_{\rm long}  = R_0 + t_0$.

However, since $\gamma_{g_1}([0, R_{\rm long}], \xi'_{\rm long})$ is the unique minimizer between end points and does not contain conjugate points,
\begin{eqnarray}
\label{g1 geod gets far away}
\gamma_{g_1}( R_{\rm long}+s, \xi'_{\rm long})\notin \overline{ B_{g_1}(\0,R_{\rm long}-s)}
\end{eqnarray}
for all $s>0$.

By \eqref{def of t0 and t2}, $R_{\rm long} = R_0 + t_0$ so for $\FLO_{g_1}^+(0,\0,\xi_{\rm long})$ to intersect 
$$ \{(t_0 + 2R_0 -t, x) \mid d_{g_1}(\0, x) \leq t,\ t\in(0,R_0)\},$$
there must be an $s>0$ such that 
$$\gamma_{g_1}(s + R_{\rm long}, \xi'_{\rm long}) \in \overline{B_{g_1}(\0,  R_0 -s)}.$$
This contradicts \eqref{g1 geod gets far away} since $R_0<R_{\rm long}$. Combine this with the expression \eqref{Ig1- is union of balls} we can conclude that 
$$I_{g_1}^-(t_0 + 2R_0, \0)\cap\left((R_{\rm long},\infty)\times M_1\right)\cap \FLO_{g_1}^+(0,\0,\xi_{\rm long}) = \emptyset.$$
So \eqref{does not intersect in M1 1} holds.
\end{proof}
We also have
\begin{lemma}
\label{not in g1 flo from xi1}
$$(t_0 + 2R_0, \0)\notin \pi\circ \FLO_{g_1}^+((t_1,x_1,\xi_1))$$
\end{lemma}
\begin{proof}
Condition \eqref{g1 is minimizing} states that $\gamma_{g_1}([0,R_{\rm long}], \xi'_{\rm long})$ is the unique minimizing segment between $\0$ and the end point. So we can conclude that $\0\notin \gamma_{g_1}((0,2R_{\rm long}], \xi'_{\rm long})$. By \eqref{concatenation} we have that 
\begin{eqnarray}
\label{0 notin R1+Rlong}
\0\notin \gamma_{g_1}([0, R_1 + R_{\rm long}],\xi_1').
\end{eqnarray}

Note that by definition \eqref{def of t0 and t2} and \eqref{R1 = Rlong -t1}
\begin{eqnarray}
\label{Rlong = t1+R1 = t0+R0}
R_{\rm long} = t_1 + R_1 = t_0 +R_0
\end{eqnarray}
Writing explicitly
$$\pi\circ \FLO_{g_1}^+((t_1,x_1,\xi_1)) = \{ (t_1 + s , \gamma_{g_1}(s, \xi'_1))\mid s>0\}$$
we see that if 
\begin{eqnarray}
\label{t0+2R0,0 in}
(t_0 + 2R_0, \0)\in \pi\circ \FLO_{g_1}^+((t_1,x_1,\xi_1))
\end{eqnarray}
then 
\begin{eqnarray}
\label{0=t0+2R0-t1}
\0 = \gamma_{g_1}(t_0 + 2R_0 -t_1, \xi'_1).
\end{eqnarray}
 By \eqref{Rlong = t1+R1 = t0+R0}, 
$$t_0 + 2R_0 -t_1 = R_{\rm long} + R_0 - t_1 <R_1+  R_{\rm long}.$$
So \eqref{0=t0+2R0-t1} would mean that $\0\in  \gamma_{g_1}([0, R_1 + R_{\rm long}],\xi_1')$, contradicting \eqref{0 notin R1+Rlong}. So \eqref{t0+2R0,0 in} must be false.
\end{proof}
With these facts established we are now ready to proceed with 

\begin{proof}[Proof of Proposition \ref{unique minimizer from origin}]
Let us suppose the contrary. Let $t_0$ be defined by \eqref{def of t0 and t2}, recall that $x_0 = \0$ was defined earlier, and $\xi'_0\in S_{\0}^*\U$ be as constructed in Lemma \ref{if not minimize}. Let $t_1$ be as in (\ref{minimizing if shorter}) of Lemma \ref{if not minimize}, $x_1$ be as in \ref{x1 def geodesic from origin}, $\xi'_1$  as in \eqref{xi'1 def geodesic from origin}. Let $t_2$ be defined by \eqref{def of t0 and t2}, $x_2$ be defined by \eqref{x2 def geodesic from origin}, and $\xi'_2$ be defined by \eqref{xi'2 def geodesic from origin}.

Recall the definition $\xi_l:= -dt + \xi'_l \in T^*_{(t_l, x_l)}\U$. Due to \eqref{does not intersect in M1} we can choose $h_0>0$ small so that 
\begin{eqnarray}
\label{flowout wedges do not intersect}\nonumber
I^-_{g_1}(t_0+2R_0+\delta, \0)\cap \pi\circ\FLO_{g_1}^+(\B_{h_0}(\xi_1)\cap\Char_{g_1}) \cap \pi\circ \FLO_{g_1}^+(\B_{h_0}(\xi_0)\cap \Char_{g_1}) = \emptyset\\
\end{eqnarray}
for all $\delta >0$ sufficiently small.

Due to \eqref{0notingammag1}, \eqref{0 notin g1 geodesic xi2}, and Lemma \ref{not in g1 flo from xi1}, $h_0>0$ can also be chosen so that
{\Small\begin{eqnarray}
\label{not if flo of g1 xi0 xi1 xi2}
B_G(t_0+ 2R_0, \0; \delta) \cap \left( \bigcup_{l=0}^2 \pi\circ\FLO_{g_1}^+(\B_{h_0}(\xi_l)\cap \Char_{g_1})\right) = \emptyset
\end{eqnarray}}
for all $\delta >0$ sufficiently small.

We now verify the hypothesis of Proposition \ref{main 3fold interaction thm}. Note that due to \eqref{t1 is small}, \eqref{doesn't get in the way} hold on $\M_2$. Furthermore due to \eqref{e: span is 3}, \eqref{span is only 2} and \eqref{full span} hold on $\M_2$. The condition \eqref{linear independent at the target} is satisfied by \eqref{no circling geod}. In the current setting, \eqref{three flowouts meet in M2} states that the three $g_2$-lightlike geodesics meet at $(R_0+t_0, y_0)$ which is analogous to \eqref{R0y0 in flowout} except that the time has shifted by $t_0$.

So we are allowed to evoke Proposition \ref{main 3fold interaction thm}. This means, for each $l=0,1,2$, there exists a sequence $\{ \xi_{l;j}\}_{j=1}^\infty \subset T^*_{(t_l,x_l)}\U$ converging to $\xi_l$ so that if $\{f_{l;j}\}_{j}$ are sequences of ($h$-dependent) distributions of the form \eqref{source} with 
$$\xi_{l;j}\in \WF(f_{l;j})\subset \B_h( \xi_{l;j})\subset \B_{h_0}(\xi_l),\ \sigma(f_{l;j})(t_l, x_l, \pm\xi_{l;j})\neq 0$$
and $w^j$ are solutions of \eqref{cubic wave} in $\M_2$ with source $f^j :=\sum_{l=0}^2 \epsilon_l f_{l;j}$ then each $\singsupp (w_{012}^j)$ contains a point $(\tilde T_j, \tilde x_j)$ with $\tilde T_j > 2R_0+t_0$ and $\tilde x_j\neq \0$ converging to $(2R_0+t_0, \0)$ as $j\to\infty$. Furthermore for each $j\in \nn$, if $h>0$ is small enough there are open sets $\tilde{\mathcal O}_j$ containing $(\tilde T_j, \tilde x_j)$ with $\tilde{\mathcal O}_j \cap\left(\R\times \{\0\}\right) = \emptyset$ such that for some lightlike hypersurfaces $\tilde S_j\subset \M_2$,
$$w_{012}^j\mid_{\tilde{\mathcal O}_j}\in I(\M_2, N^*\tilde S_j)$$
and
$$\sigma(w_{012}^j)(\tilde T_j, \tilde x_j, \tilde\xi_j) \neq 0$$
for all $\tilde\xi_j \in N^*_{(\tilde T_j, \tilde x_j)}\tilde S_j$.

We now choose $f_{3;j}$ and $f_{4;j}$ as in \eqref{def of f4 f5}. Recall that these two distributions depend on a parameter $a>0$ satisfying $a<<h$.
Let $v^j$ solve \eqref{cubic wave} in $\M_2$ with source $f^j := \sum_{l=0}^4f_{l;j}$ which depends now on the parameters $0<a<<h<<h_j$ where $h_j>0$ is a constant depending on $j\in\nn$. Then according to Proposition \ref{fivefold interaction}, $v^j_{01234}(\cdot,\0)$ has a singularity at $\tilde T_j + d_{g_2}(\tilde x_j,\0)$ for all $0<a<<h<<h_j$. By Lemma \ref{differentiate then restrict},
$$\tilde T_j + d_{g_2}(\tilde x_j,\0)\in \singsupp\left(\partial^5_{\epsilon_0\dots\epsilon_4}\left(v^j(\cdot, \0)\right)\mid_{\epsilon_0=\cdots= \epsilon_4=0}\right).$$

Let $u^j$ solve \eqref{cubic wave} in $\M_1$ and source $f^j = \sum_{l=0}^4f_{l;j}$. By the fact that source-to-solution maps are equal, $u^j(\cdot; \0) = v^j(\cdot;\0)$ so 
$$\tilde T_j + d_{g_2}(\tilde x_j,\0)\in \singsupp\left(\partial^5_{\epsilon_0\dots\epsilon_4}\left(u^j(\cdot, \0)\right)
\mid_{\epsilon_0=\cdots= \epsilon_4=0}\right).$$
By Lemma \ref{differentiate then restrict},
\begin{eqnarray}
\label{singsupp of uj}
\tilde T_j + d_g(\tilde x_j,\0)\in \singsupp(u^j_{01234}(\cdot,\0))
\end{eqnarray}
for all $0<a<<h<<h_j$ for some $h_j$ depending on $j\in\nn$. We will show that \eqref{singsupp of uj} is false and thus obtain a contradiction.

The distribution $u^j_l$ solves
$$\Box_{g_1}u^j_l = f_{l;j},\ u^j_l\mid_{t<-1} = 0.$$
Proposition \ref{propagation} states that 
\begin{eqnarray}
\label{singsuppujl no boomerang}
u^j_l \in I(\M_1, T_{(t_l,x_l)}^*\M_1 ,\FLO_{g_1}^+(\ccl(\B_{h_0}(\xi_l)\cup\B_{h_0}(-\xi_l))\cap\Char_{g_1}))
\end{eqnarray}
 for $l= 0,1,2$. So due to \eqref{flowout wedges do not intersect} and the fact that $(\tilde T_j, \tilde x_j) \to (2R_0 + t_0, \0)$, for $j\in \nn$ large enough,
\begin{eqnarray}
\label{triple intersection is empty for ujl}
\singsupp(u^j_0)\cap\singsupp(u^j_1)\cap\singsupp(u^j_2)\cap I_{g_1}^-((\tilde T_j+ d_{g_2}(\tilde x_j, \0),\0)) = \emptyset.
\end{eqnarray}
Note the $u^j_{012}$ solves 
$$\Box_{g_1} u^j_{012} = -6 u^j_0 u^j_1 u^j_2,\ u^j_{012} \mid_{t <-1} = 0.$$
As only lightlike wavefronts propagate, Lemma \ref{only threefold propagate} combined with \eqref{triple intersection is empty for ujl} forces 
{\small \begin{eqnarray}
\label{backward causal intersects three singsupp}
 \singsupp(u^j_{012}) \cap  I_{g_1}^-((\tilde T_j+ d_{g_2}(\tilde x_j, \0),\0))\subset \left(\bigcup_{l=0}^2\singsupp(u^j_l) \right)\cap I_{g_1}^-((\tilde T_j+ d_{g_2}(\tilde x_j, \0),\0)).
\end{eqnarray}}
Combining the above inclusion with \eqref{not if flo of g1 xi0 xi1 xi2} and \eqref{flowout wedges do not intersect} we have that there exists $\delta>0$ such that
\begin{eqnarray}
\label{singsupp uj012 is empty}
 \singsupp(u^j_{012}) \cap  I_{g_1}^-((\tilde T_j+ d_{g_2}(\tilde x_j, \0),\0))\cap B_G(2R_0 + t_0,\0; \delta) = \emptyset
\end{eqnarray}
for all $j\in\nn$ sufficiently large.

By \eqref{v3 and v4} we have
\begin{eqnarray}
\label{uj3}
u^j_3 = \left( \chi_{a;j}\langle D\rangle^{-N} \delta_{{\bf T}_j}\right)
\end{eqnarray}
and 
\begin{eqnarray}
\label{uj4}
u^j_4 = \chi_{a;j}(t,x).
\end{eqnarray}
 We recall for convenience of reader that the definitions of $\chi_{a;j}$ and $\delta_{{\bf T}_j}$ are given in \eqref{def of f4 f5}. Note that if $0<a<<h$, 
 $$\supp(\chi_{a;j}) \subset B_G(\tilde x_j, \tilde T_j;a)\subset\subset B_G(2R_0 + t_0,\0 ;\delta)$$ 
 where $\tilde x_j \to \0$ and $\tilde T_j \to t_0 + 2R_0$. Due to \eqref{singsupp uj012 is empty} we have that for all $j\in\nn$ sufficiently large and $a>0$ sufficiently small,
\begin{eqnarray}
\label{wf of uj012 uj2 uj4}
\WF(u^j_{012} u^j_3u^j_4) \subset \WF \left( \chi_{a;j}(t,x)     \langle D\rangle^{-k} (\chi_{a;j} \delta_{{\bf T}_j})\right) \subset N^*{\bf T}_j.
\end{eqnarray}
Recall that ${\bf T}_j = \R \times \{\tilde x_j\}$ where $\tilde x_j\neq \0$.

The distribution $u^j_{01234}$ solves the equation
\begin{eqnarray*}
\Box_{g_1} u^j_{01234} &=& C\sum_{\sigma\in S_5} u^j_{\sigma(0)}u^j_{\sigma(1)} u^j_{\sigma(2)\sigma(3)\sigma(4)},\\
& =&C u^{j}_{012} u^j_3 u^j_4 +C\sum_{\sigma\in S_5\backslash S_3} u^j_{\sigma(0)}u^j_{\sigma(1)} u^j_{\sigma(2)\sigma(3)\sigma(4)}
\end{eqnarray*}
with initial condition $u^j_{01234}\mid_{t<-1} = 0$. Recall that we denote by $S_3\subset S_5$ to be the subgroup which maps the set $\{0,1,2\}$ to itself. We have the following 
{\begin{lemma}
\label{regularity of wj}
For $j\in \nn$ sufficiently large, let $u^j_{\rm reg}$ solve
$$\Box_{g_1} u_{\rm reg}^j = C\sum_{\sigma\in S_5\backslash S_3} u^j_{\sigma(0)}u^j_{\sigma(1)} u^j_{\sigma(2)\sigma(3)\sigma(4)}$$
with initial condition $u^j_{\rm reg}\mid_{t<-1} = 0$. Then $(\tilde  T_j + d_{g_2}(\0, \tilde x_j), \0)\notin\singsupp(u^j_{\rm reg})$.
\end{lemma}}
This lemma implies that 
 $(\tilde T_j + d_{g_2}(\0, \tilde x_j), \0)\notin\singsupp(u_{01234}^j)$ 
since the wavefront of the term $u_{012}^j u^j_3 u^j_4$ is given by \eqref{wf of uj012 uj2 uj4} and therefore does not propagate. This contradicts \eqref{singsupp of uj}. The proof is complete.
\end{proof}


We still need to give a proof for Lemma \ref{regularity of wj}. However, it is essentially the same argument as proof of Lemma \ref{not main term}. So we only give a sketch:
\begin{proof}[Proof of Lemma \ref{regularity of wj}]
Use the fact that $\supp(u_3^j)$ and $\supp(u_4^j)$ are contained in $$B_G(\tilde T_j, \tilde x_j; a) \subset\subset B_G(t_0+2R_0,\0;\delta)$$ by \eqref{uj3} and \eqref{uj4}. Note that for $l=0,1,2$, 
$$\singsupp(u_l^j) \subset\pi\circ \FLO_{g_1}^+(\ccl B_{h_0}(\xi_{l})\cap\Char_{g_1})$$
by \eqref{singsuppujl no boomerang}. So \eqref{not if flo of g1 xi0 xi1 xi2} implies that 
$$\singsupp(u^j_l)\cap \left(\supp(u_3^j)\cup \supp(u^j_4)\right) = \emptyset.$$
With these facts about how singular supports are disjoint, we proceed as in the proof of Lemma \ref{not main term}.
\end{proof}

\subsection{Determining the Distance Function}

\begin{proposition}
\label{distances are equal}
Let $R_0<T$ and suppose for $\xi'_0\in S^*_{\0}\U$ the geodesic $\gamma_{g_2}([0,R_0], \xi'_0)$ is a minimizer between end points $\0$ and $y_2:= \gamma_{g_2}(R_0, \xi'_0)  \in M_2$. Set 
\begin{eqnarray}
\label{y1 is in g1 geod}
y_1 := \gamma_{g_1}(R_0, \xi'_0)\in M_1.
\end{eqnarray}
Then 
\begin{eqnarray}
\label{d1=d2}
d_{g_1}(x, y_1) \leq d_{g_2}(x, y_2)
\end{eqnarray}
for all $x\in \U$.
\end{proposition}
Note that thanks to Proposition \ref{unique minimizer from origin}, the statement of Proposition \ref{distances are equal} is actually symmetric with the subscript $1$ and $2$ interchanged. Therefore, it has the following 
\begin{corollary}
\label{distances are equal coro}
For every $y_2\in B_{g_2}(\0; T)$ there exists a unique $y_1\in B_{g_1}(\0;T)$ such that $d_{g_1}(x,y_1) = d_{g_2}(x,y_2)$ for all $x\in \U$. The statement also holds with $1$ and $2$ interchanged.
\end{corollary}
\begin{proof}
Let $\xi'_0\in S^*_\0\U$ and $\gamma_{g_2}([0,R], \xi'_0)$ be a minimizing geodesic segment with end points $\0$ and $y_2$. Set $y_1 := \gamma_{g_1}(R, \xi'_0)$ and by Proposition \ref{unique minimizer from origin}, $\gamma_{g_1}([0,R], \xi'_0)$ is a minimizing segment between $\0$ and $y_1$. By Proposition \ref{distances are equal}, $d_{g_1}(x,y_1) \leq d_{g_2}(x,y_2)$ for all $x\in \U$. Similarly, switching the role of $y_1$ and $y_2$ yields $d_{g_2}(x,y_2)\leq d_{g_1}(x,y_1)$ for all $x\in \U$. Uniqueness is a trivial consequence of the fact that if $y, y'\in M$ and $d_g(y,x) = d_g(y', x)$ for all $x$ in some open set $\U\subset M$, then $y = y'$.
\end{proof}

We now proceed to prove Proposition \ref{distances are equal}. We first make some reductions using the continuity of the distance function. By Proposition \ref{unique minimizer from origin}, for both $j=1,2$, $\gamma_{g_j}([0, R_0], \xi'_0)\subset M_j$ are minimizing segments between $\0$ and $\gamma_{g_j}(R_0,\xi'_0)\in M_j$. Taking $R_0\in (0,T)$ slightly smaller we may assume without loss of generality that both
{\small\begin{eqnarray}
\label{unique min assumption}
\gamma_{g_1}([0, R_0], \xi'_0),\ \gamma_{g_2}([0, R_0], \xi'_0)\ {\mbox {are the unique minimizers and have no conjugate points.}}
\end{eqnarray}}
Another observation is that with $y_1\in M_1$ and $y_2\in M_2$ fixed, it suffices to prove \eqref{d1=d2} for a dense subset of $x\in \U$. As such we only need to prove \eqref{d1=d2} for $x\in \U$ which is joined to $y_1\in M_1$ and $y_2\in M_2$ by unique minimizing geodesics in $M_1$ and $M_2$ which do not contain conjugate points. Furthermore we assume that 
\begin{eqnarray}
\label{x not in the way}
x\notin \gamma_{g_1}([-\delta_0, \delta_0], \xi'_0) = \gamma_{g_2}([-\delta_0, \delta_0], \xi'_0).
\end{eqnarray}


To this end let $x_1\in \U$ satisfy the above assumptions. 
Let $\gamma_{g_2}([0,R_1],\xi'_1)$ be the unique minimizing segment between $x_1$ and $y_2$ in $M_2$ that does not contain conjugate points. The distance $d_{g_2}(x_1, y_2)$ is given by
\begin{eqnarray}
\label{R1}
R_1 = d_{g_2}(x_1,y_2).
\end{eqnarray}
Again, by density, we may without loss of generality assume that 
\begin{eqnarray}
\label{geod 1 and 2 are not tangent}
\dot\gamma_{g_2}(R_0, \xi'_0) \neq - \dot \gamma_{g_2}(R_1, \xi'_1).
\end{eqnarray}
By the fact that $\gamma_{g_2}([0,R_0],\xi'_0)$ is the unique minimizer between end points, \eqref{geod 1 and 2 are not tangent} implies
\begin{eqnarray}
\label{origin not in g2 geod}
\0\notin \gamma_{g_2}((0,2R_0], \xi'_0) \cup \gamma_{g_2}([0, R_0+R_1], \xi'_1).
\end{eqnarray}
Note that since both geodesic segments $\gamma_{g_2}([0,R_0], \xi'_0)$ and $\gamma_{g_2}([0,R_1], \xi'_1)$ are unique minimizers, the two geodesic segments only intersect at the point $y_2$.

We set $\eta'_l := -\dot \gamma_{g_2}(R_l, \xi'_l)^\flat$ for $l=0,1$ and use Lemma \ref{full span lem} to choose $\eta'_2\in S^*_{y_2}M_2\cap \spn\{\eta_0', \eta_1'\}$ arbitrarily close to $\eta'_0$ so that 
\begin{eqnarray}
\label{span downstairs is 2}
{\rm Dim}\left(\spn\{-dt - \eta'_{0}, -dt - \eta'_{1}, -dt - \eta'_{2}\}\right) = 3,\ {\rm Dim}(\spn\{\eta_0', \eta_1', \eta_2'\} )= 2.
\end{eqnarray}
Set $R_2 := R_0$. By Lemma \ref{uniform unique minimizer} if $\eta'_2$ is chosen close to $\eta_0'$, we can conclude that $\gamma_{g_2}([0,R_2],\eta'_2)$ is a unique minimizing segment containing no conjugate points. Set 
$$\xi'_2 := - \dot \gamma_{g_2}(R_2, \eta'_2)^\flat$$
so that
$$y_2 = \gamma_{g_2}(R_0, \xi'_2) =\gamma_{g_2}(R_0,\xi'_0) = \gamma_{g_2}(R_1, \xi'_1).$$
If $\eta'_2$ is chosen sufficiently close to $\eta'_0$, condition \eqref{geod 1 and 2 are not tangent} allows us to assert that 
\begin{eqnarray}
\label{geod k and l are not tangent}
\dot\gamma_{g_2}(R_l, \xi'_0) \neq - \dot \gamma_{g_2}(R_k, \xi'_1)
\end{eqnarray}
for $l,k\in \{0,1,2\}$.

Due to \eqref{origin not in g2 geod}, if $\eta'_2$ is chosen sufficiently close to $\eta'_0$ we have
\begin{eqnarray}
\label{origin not in g2 geod 2}
\0\notin \gamma_{g_2}((0,2R_0], \xi'_0) \cup \gamma_{g_2}([0, R_0+R_1], \xi'_1) \cup \gamma_{g_2}([0,2R_0], \xi'_2).
\end{eqnarray}
By Proposition \ref{no boomerang} the same holds for the $g_1$ geodesics:
\begin{eqnarray}
\label{origin not in g1 geod}
\0\notin \gamma_{g_1}((0,2R_0], \xi'_0) \cup \gamma_{g_1}([0, R_0+R_1], \xi'_1)\cup \gamma_{g_1}([0,2R_0], \xi'_2).
\end{eqnarray}
Set $t_0 = t_2 = 0$ and $t_1 = R_0-R_1$ so that $t_l + R_l = R_0$ for all $l=0,1,2$. Note that by \eqref{R1} this means
\begin{eqnarray}
\label{distance x1 to y0}
d_{g_2}(x_1, y_2) = R_0 - t_1.
\end{eqnarray}
Also, by the strict triangle inequality, $d_{g_2}(x_0, x_1)>|t_1|$ and $d_{g_2}(x_2, x_1) >|t_1|$. So
\begin{eqnarray}
\label{source not in future causal}
(t_l,x_l) \notin I_{g_2}^+(t_k,x_k),\ k\neq l.
\end{eqnarray}

Define the lightlike covectors
\begin{eqnarray}
\label{setting xil}
\xi_l = -dt+ \xi'_l \in T^*_{(t_l,x_l)}\M_2
\end{eqnarray}
for $l=0,1,2$ so that
\begin{eqnarray}
\label{R0y2 in flowout}
(R_0, y_2) \in\bigcap_{l=0}^2\pi\circ \FLO_{g_2}^+(t_l,x_l, \xi_l)
\end{eqnarray}

 By \eqref{origin not in g1 geod} there exists an $h_0>0$ and $\delta>0$ such that 
\begin{eqnarray}
\label{FLO of g1 avoids 2R0}
B_G(2R_0,\0; \delta)\cap\left(\bigcup_{l=0}^2 \pi\circ\FLO_{g_1}^+(\ccl\B_{h_0}(\xi_l) \cap\Char_{g_1})\right) = \emptyset
\end{eqnarray}
Furthermore we can choose the $\delta>0$ in \eqref{FLO of g1 avoids 2R0} to satisfy
\begin{lemma}
\label{unique triple intersection}
i)There exists a $\delta>0$ and $h_0>0$ such that if $\xi\in \B_{h_0}(\xi_0)\cap \Char_{g_1}$ then the lightlike geodesic segment
$$I_{g_1}^-(B_G(2R_0,\0; \delta))\cap \pi\circ\FLO_{g_1}^+(0,\0,\xi)$$
is the unique lightlike geodesic segment joining any two points on it.\\
ii) $I_{g_1}^-(2R_0, \0)\cap \pi\circ\FLO_{g_1}^+(0,\0, \xi_0) \subset [0,R_0]\times M_1$.
\end{lemma}
\begin{proof}
i) By assumption $\gamma_{g_1}([0,R_0]; \xi'_0)$ is the unique minimizer between end points and does not have conjugate points. By Lemma \ref{uniform unique minimizer} there exists an open set $U'\subset S^*_\0\U$ and $\delta' >0$ such that if $\xi'\in U'$ then $\gamma_{g_1}([0, R_0+\delta'], \xi')$ is the unique minimizer between end points without conjugate points. A consequence of this is that 
\begin{eqnarray}
\label{doesn't come back in time}
d_{g_1}\left(\gamma_{g_1}(R_0 +\delta'/2 + t, \xi'), \0 \right) > R_0 + \delta'/2 -t
\end{eqnarray} 
for all $t\geq 0$ and $\xi'\in U'$.
Furthermore for all $\xi'\in U'$, the lightlike geodesic segment
$$\{\pi\circ e_+(s, 0,\0,-dt + \xi')\mid s\in [0,R_0+\delta']\}$$
is the unique lightlike geodesic joining any two points on it. Due to \eqref{doesn't come back in time}, for $\xi'\in U'$
$$I^-_{g_1}(2R_0+\delta'/2, \0)\cap \FLO_{g_1}^+(0,\0,-dt +\xi') \subset \{(t,x)\in \M_1\mid t\in (0,R_0+ \delta'/2)\}.$$

So if we choose $h_0>0$ to satisfy 
$$\left(\B_{h_0}(\xi_0) \cap\Char_{g_1}\right)/ \R^+  \subset \{-dt + \xi'\mid \xi'\in U'\},$$
then for any $\xi\in \B_{h_0}(\xi_0) \cap\Char_{g_1}$, the geodesic segment 
$$I_{g_1}^-(2R_0 + \delta'/2, \0) \cap\pi\circ \FLO_{g_1}^+(0,\0,\xi)$$
is the only lightlike geodesic  segment joining any two points on it.
Therefore, the lemma is verified if $\delta>0$ is chosen small enough so that $B_G(2R_0, \0; \delta) \subset\subset I_{g_1}^-(2R_0 + \delta'/2, \0)$.\\

ii) This can be seen by taking $\xi'= \xi'_0$ in \eqref{doesn't come back in time} and observing that 
$$I_{g_1}^-(2R_0,\0) = \{(2R_0-t,x) \in \M_1\mid d_{g_1}(x,\0)\leq t,\ t\geq 0\}.$$
\end{proof}

The conditions \eqref{span downstairs is 2}, \eqref{R0y2 in flowout}, \eqref{source not in future causal},\eqref{geod k and l are not tangent} allows us to evoke Proposition \ref{fivefold interaction} to conclude
\begin{lemma}
\label{the setup}
We have that
\begin{enumerate}
\item \label{first 3 covectors} For each $l=0,1,2$ there exists a sequence $\xi_{l;j}\in T^*_{(t_l,x_l)}\M_2\cap \Char_{g_2}$ converging to $\xi_l$.
\item \label{first 3 sources} For each $l=0,1,2$ and $j\in\nn$ large there is an $h_j\in (0,h_0)$ with 
$$\B_{h_j}(\xi_{l;j})\subset\subset \B_{h_0}(\xi_{l})$$
so that if $h\in(0,h_j)$ then there exists sources $f_{l;j}(\cdot;h)$ of the form \eqref{source} such that 
$$\xi_{l;j}\in\WF(f_{l;j}(\cdot;h))\subset\ccl( \B_h(\xi_{l;j})\cup\B_h(-\xi_{l;j})),\ \sigma(f_{l;j}(\cdot;h))(t_l, x_l, \xi_{l;j})\neq 0.$$
\item \label{conormal point} A sequence $(\tilde T_j, \tilde x_j)\to (2R_0, \0)$, $\tilde T_j >2R_0$, and $\tilde x_j\neq \0$.
\item \label{choice of a} For any $j\in\nn$ large and $h\in (0,h_j)$ as above, and any $a\in (0,h)$ sufficiently small, we have sources $f_{3;j}(\cdot ;a)$ and $f_{4;j}(\cdot ;a)$ of the form \eqref{def of f4 f5} which are supported in $B_{G}(\tilde T_j, \tilde x_j;a)$.
\end{enumerate}
If $(t,x)\mapsto v^j(t,x; h,a)$ are solutions of \eqref{cubic wave} in $\M_2$ with source 
 $$f^j(\cdot; h,a) := \sum_{l=0}^2 \epsilon_l f_{l;j}(\cdot; h) + \sum_{l=3}^4\epsilon_l f_{l;j}(\cdot ;a)$$ 
 then for each $j\in\nn$ large
\begin{eqnarray}
\label{singsupp of vj01234, parameters}
\tilde T_j + d_{g_2}(\tilde x_j, \0)\in \singsupp(v^j_{01234}(\cdot,\0;h,a)).
\end{eqnarray}
for all $0<a<h<h_j$ small.
\end{lemma}
Note that we now explicitly write out the dependence of the sources $f_{l;j}$ and solution $v^j$ on the parameters $0<a<h<h_j$.



Let $u^j(\cdot;h,a)$ solve \eqref{cubic wave} with metric $g_1$ and the same source as $v^j(\cdot;h,a)$. The fivefold interaction of the nonlinear wave $u^j(\cdot;h,a)$ is 
\begin{eqnarray}
\label{fivefold u}\nonumber
\Box_{g_1} u^j_{01234}(\cdot;h, a) = \sum_{\sigma\in S_5} u^j_{\sigma(0)\sigma(1)\sigma(2)}(\cdot;h, a) u^j_{\sigma(3)}(\cdot;h, a) u^j_{\sigma(4)}(\cdot;h, a),\ u^j_{01234}(\cdot;h, a)\mid_{t<-1} = 0\\
\end{eqnarray}
We write the solution $u^j_{01234}(\cdot;h, a)$ of \eqref{fivefold u} as $u^j_{\rm reg} + u^j_{\rm sing}(\cdot ; h,a)$ where 
\begin{eqnarray}
\label{ureg eq}
\Box_{g_1} u^j_{\rm reg} = \sum_{\sigma\in S_5\backslash S_3} u^j_{\sigma(0)\sigma(1)\sigma(2)}(\cdot;h, a) u^j_{\sigma(3)}(\cdot;h, a) u^j_{\sigma(4)}(\cdot;h, a),\ u^j_{\rm reg}\mid_{t<-1}  = 0
\end{eqnarray}
and 
\begin{eqnarray}
\label{using eq}
\Box_{g_1} u^j_{\rm sing}(\cdot;h,a) = 6 u^j_{012}(\cdot;h,a)u^j_3(\cdot;h,a)u_4^j(\cdot;h,a),\  u^j_{\rm sing}(\cdot;h,a) \mid_{t<-1} = 0.
\end{eqnarray}
Here we denote $S_3\subset S_5$ as elements of the permutation group which maps the set of three letters $\{0,1,2\}$ to itself. Repeating the argument of Lemma \ref{not main term} we get that 
{\begin{lemma}
\label{not main term u}
If $j\in \nn$ is chosen large enough and $a\in (0,h)$ is chosen small enough so that $B_G(\tilde T_j, \tilde x_j; a) \subset \subset B_G(2R_0, \0; \delta)$ and $B_G(\tilde T_j, \tilde x_j; a) \cap (\R\times\{\0\}) = \emptyset$, then
$$(\tilde T_j + d_g(\tilde x_j,\0), \0) \notin \singsupp(u_{\rm reg}).$$
\end{lemma}}
\begin{proof}
This is similar to proof of Lemma \ref{not main term} so we will only give a brief sketch.
To simplify notation we define 
$$f_{\rm source} :=  \sum_{\sigma\in S_5\backslash S_3} u^j_{\sigma(0)\sigma(1)\sigma(2)}(\cdot;h, a) u^j_{\sigma(3)}(\cdot;h, a) u^j_{\sigma(4)}(\cdot;h, a).$$
Note that by \eqref{v3 and v4}
$$u^j_4 = \chi_{a;j} \in C^\infty_c(B_G(\tilde T_j, \tilde x_j; a)),\ u^j_3 =  \chi_{a;j}\langle D\rangle^{-N} \delta_{{\bf T}_j}\in I(\M, N^*{\bf T}_j).$$
Meanwhile by Proposition \ref{propagation}, for $l=0,1,2$
$$u^j_l \in I(\M_1, T^*_{(t_l,x_l)}\M_1, \FLO_{g_1}^+(\ccl(\B_{h_0}(\xi_l) \cup \B_{h_0}(-\xi_l))\cap\Char_{g_1}).$$
This combined with flowout condition \eqref{FLO of g1 avoids 2R0} ensures that for $l = 0,1,2$,
$$\singsupp(u^j_l)\cap \left(\supp(u_3^j) \cup \supp(u_4^j)\right) = \emptyset.$$
Using these facts we can proceed as in Lemma \ref{not main term} to conclude that 
$$\FLO_{g_1}^-(\tilde T_j + d_g(\tilde x_j,\0), \0, \xi) \cap \WF(f_{\rm source}) = \emptyset$$ 
for all $\xi\in T^*_{(\tilde T_j + d_g(\tilde x_j,\0), \0)}\U\cap\Char_{g_1}$. The lemma then follows from  Thm 23.2.9 of \cite{Hormander:2007aa}. 
\end{proof}
With these facts we can now give the
\begin{proof}[Proof of Proposition \ref{distances are equal}]
From \eqref{singsupp of vj01234, parameters} we can use Lemma \ref{differentiate then restrict} to deduce
$$\tilde T_j + d_{g_2}(\tilde x_j,\0)\in \singsupp\left(\partial^5_{\epsilon_0\dots\epsilon_4}\left(v^j(\cdot, \0;h,a)\right)\mid_{\epsilon_0=\cdots= \epsilon_4=0}\right).$$

By the fact that the source-to-solution maps agree, $t\mapsto u^j(t,\0 ; h,a)$ is the same as $t\mapsto v^j(t,\0; h,a)$. So we have that 
$$\tilde T_j + d_{g_2}(\tilde x_j,\0)\in \singsupp\left(\partial^5_{\epsilon_0\dots\epsilon_4}\left(u^j(\cdot, \0;h,a)\right)\mid_{\epsilon_0=\cdots= \epsilon_4=0}\right).$$
Lemma \ref{differentiate then restrict} then implies $\tilde T_j + d_g(\tilde x_j,\0)\in \singsupp(u^j_{01234}(\cdot, \0;h,a))$. By Lemma \ref{not main term u} we have that
\begin{eqnarray}
\label{using is singular}
 (\tilde T_j+ d_g(\tilde x_j,\0), \0)\in\singsupp(u^j_{\rm sing}(\cdot;h,a)).
 \end{eqnarray}
Note that for all $a>0$ small, the distribution $u^j_{\rm sing}(\cdot;h,a)$ solves \eqref{using eq} and the source term of \eqref{using eq} is supported in $B_{G}(\tilde T_j, \tilde x_j; a)$ since $\supp(u_4^j(\cdot;a))\subset B_{G}(\tilde T_j, \tilde x_j; a)$ by \eqref{v3 and v4}. Note that by Condition \eqref{conormal point} of Lemma \ref{the setup}, $x_j\neq \0$. The solution $u^j_4(\cdot;a)$ is smooth and the wavefront set of $u^j_3(\cdot;a)$ is spacelike. So in order for \eqref{using is singular} to hold, 
$$B_{G}((\tilde T_j, \tilde x_j),a)\cap  \singsupp(u^j_{012}(\cdot;h))\neq\emptyset $$
for all $a>0$. We remark that $u^j_{012}(\cdot;h)$ does not depend on the parameter $a>0$ since $f_{l;j}(\cdot; h)$ depends only on the parameter $h>0$ for $l=0,1,2$. So taking intersection over all $a>0$, we conclude that, for each fixed $j\in\nn$ large, if $h\in (0,h_j)$ is sufficiently small,
\begin{eqnarray}
\label{u012 is singular}
(\tilde T_j, \tilde x_j) \in\singsupp(u^j_{012}(\cdot;h)).
\end{eqnarray}
Fix $j\in \nn$ large and for each $h>0$ sufficiently small $u^j_{012}(\cdot;h)$ solves
$$\Box_{g_0} u^j_{012}(\cdot;h) = -6u^j_{0}(\cdot;h)u^j_{1}(\cdot;h)u^j_{2}(\cdot;h),\ u^j_{012}(\cdot;h)\mid_{t<-1} = 0$$
with each of $u^j_l(\cdot; h)$ solving linear inhomogeneous wave equation with source $f_{l;j}(\cdot;h)$ so that $\singsupp(u^j_l(\cdot; h)) \subset \pi\circ\FLO_{g_1}^+(\WF(f_{l;j}(\cdot;h))\cap \Char_{g_1})$ by Proposition \ref{propagation}. 

Take $j\in\nn$ large enough so that $(\tilde T_j, \tilde x_j)\in B_G((2R_0, \0),\delta)$ where $\delta>0$ is chosen so that \eqref{FLO of g1 avoids 2R0} and the conclusion of Lemma \ref{unique triple intersection} holds. By Lemma \ref{only threefold propagate} and \eqref{FLO of g1 avoids 2R0}, in order for \eqref{u012 is singular} to hold we must have 
\begin{eqnarray}
\label{pastcone intersect flowouts}
I_{g_1}^-((\tilde T_j, \tilde x_j))\cap\bigcap_{l=0}^2 \pi\circ\FLO_{g_1}^+(\WF(f_{l;j}(\cdot;h))\cap \Char_{g_1}) \neq \emptyset
\end{eqnarray}
for all $h\in (0,h_j)$.

By Condition \eqref{first 3 sources} of Lemma \ref{the setup}, $\WF(f_{l;j}(\cdot;h)) \subset \ccl\B_h(\pm\xi_{l;j})$ for $l=0,1,2$ so \eqref{pastcone intersect flowouts} becomes
\begin{eqnarray}
\label{pastcone intersect flowouts 2}
I_{g_1}^-((\tilde T_j, \tilde x_j))\cap\bigcap_{l=0}^2 \pi\circ\FLO_{g_1}^+(\ccl\B_h(\xi_{l;j})\cap \Char_{g_1}) \neq \emptyset
\end{eqnarray}
for all $h>0$. 
Furthermore we must have that
{\small\begin{eqnarray}
\label{from intersection to point}
\nonumber
&&\forall h\in (0,h_j), \exists \eta_h \in T_{(t_h, x_h)}^*\M_1/\R^+\cap \Char_{g_1},\ (t_h,x_h) \in \bigcap_{l=0}^2\pi\circ\FLO_{g_1}^+(\ccl\B_h(\xi_{l;j})\cap \Char_{g_1})\\&& {\rm s.t.}\ (\tilde T_j, \tilde x_j)\in \pi\circ\FLO^+_{g_1} (t_h,x_h,\eta_h)
\end{eqnarray} }
for otherwise $(\tilde T_j, \tilde x_j)\notin \singsupp(u_{012}^j(\cdot;h))$ by Thm 23.2.9 of \cite{Hormander:2007aa}.

Since this holds for all $h>0$, compactness then dictates that there exists $(\hat R_j, \hat y_j)\in \M_1$ such that
\begin{eqnarray}
\label{triple intersection of g0}
(\hat R_j, \hat y_j) \in I_{g_1}^-((\tilde T_j, \tilde x_j))\cap\bigcap_{l=0}^2 \pi\circ\FLO_{g_1}^+(t_l,x_l,\xi_{l;j}).
\end{eqnarray}

By Lemma \ref{unique triple intersection}, $(\hat R_j, \hat y_j)$ is unique element in \eqref{triple intersection of g0}. 
Therefore we have that 
\begin{eqnarray}
\label{unique element in triple intersection}
(\hat R_j, \hat y_j) 
=I_{g_1}^-((\tilde T_j, \tilde x_j))\cap \bigcap_{l=0}^2 \pi\circ\FLO_{g_1}^+(t_l,x_l,\xi_{l;j}).
\end{eqnarray}
And by \eqref{from intersection to point}, there is a covector $\hat\eta_j\in T^*_{(\hat R_j, \hat y_j)}\M_1/\R^+\cap \Char_{g_1}$ such that 
\begin{eqnarray}
\label{arriving at the point}
(\tilde T_j, \tilde x_j) \in \pi \circ \FLO_{g_1}^+ (\hat R_j, \hat x_j, \hat \eta_j).
\end{eqnarray}
Taking a converging subsequence as $j\to \infty$ in \eqref{unique element in triple intersection} and \eqref{arriving at the point} we have, by Conditions \eqref{first 3 covectors} and \eqref{conormal point} of Lemma \ref{the setup},
\begin{eqnarray}
\label{limiting flowout}
(\hat R_0, \hat y_0)  = I_{g_1}^-(2R_0,\0) \cap \bigcap_{l=0}^2 \FLO_{g_1}^+(t_l, x_l,\xi_l),
\end{eqnarray}
\begin{eqnarray}
\label{limiting flowout'}
 (2R_0, \0) \in \pi \circ \FLO_{g_1}^+(\hat R_0, \hat y_0,\hat \eta)
\end{eqnarray}
for some $(\hat R_0, \hat y_0) \in \M_1$ and $\hat \eta \in T^*_{(\hat R_0, \hat y_0)}\M_1\cap\Char_{g_1}$. Since $\xi_0\in T^*_{(0,\0)}\M_1\cap \Char_{g_1}$ (see \eqref{setting xil}), this says that
\begin{eqnarray}
\label{both points are in g1 flowout}
 (\hat R_0, \hat y_0) \in \FLO_{g_1}^+(0,\0,\xi_0)
\end{eqnarray}
which in turn implies that
\begin{eqnarray}
\label{hat y0 on xi0}
\hat y_0 = \gamma_{g_1}(\hat R_0, \xi'_0).
\end{eqnarray}
 Evoke part ii) of Lemma \ref{unique triple intersection} combined with \eqref{both points are in g1 flowout} we conclude
\begin{eqnarray}
\label{hat R0 is smaller}
\hat R_0 \leq R_0. 
\end{eqnarray}
By \eqref{unique min assumption}, the geodesic segment $\gamma_{g_1}([0,R_0], \xi'_0)$ is the unique minimizer between the end points $\0$ and $y_1$. So \eqref{hat R0 is smaller} and \eqref{hat y0 on xi0} combined to give
\begin{eqnarray}
\label{distance between hat y0 and y1}
d_{g_1}(\hat y_0, y_1) = R_0 - \hat R_0.
\end{eqnarray}

By \eqref{setting xil}, $\xi_1\in T^*_{(t_1, x_1)}\M_1\cap \Char_{g_1}$. So since $(\hat R_0, \hat y_0) \in \FLO_{g_1}^+(t_1, x_1, \xi_1)$ by \eqref{limiting flowout}, we can conclude that 
\begin{eqnarray}
\label{hat y0 on xi1}
\hat y_0 = \gamma_{g_1}(\hat R_0 - t_1, \xi_1'),\ x_1 = \gamma_{g_1}(0,\xi'_1).
\end{eqnarray}
Combining \eqref{hat y0 on xi1} and \eqref{distance between hat y0 and y1} we get
$$d_{g_1}(x_1, y_1)  \leq d_{g_1}(x_1, \hat y_0) + d_{g_1}(\hat y_0, y_1) \leq (\hat R_0-t_1)  + (R_0- \hat R_0) = R_0 -t_1.$$
So by \eqref{distance x1 to y0} we get that
$$d_{g_1}(x_1, y_1) \leq d_{g_2}(x_1, y_2).$$

\end{proof}


\subsection{Recover Riemannian Structure from Distance Data}
We now complete the proof of Theorem \ref{metric recovery} with the following proposition. Let $(M_i,g_i)$, $i=1,2$, be complete Riemannian manifolds. For some $x_i\in M_i$ and $R>0$, we define
\begin{align*}
 V_i:=\big\{\xi'\in T_{x_i}M_i\ \big|\ d_{g_i}(x_i,\exp_{x_i}(\xi'))=\|\xi'\|_{g_i}<T\big\}.
\end{align*}
Notice that $V_i$ is precisely the set of those tangent vectors $\xi'\in T_{x_i}M_i$ whose associated geodesic $\gamma(\cdot,\xi'):[0,1]\to M_i$, $\gamma(t,\xi')=\exp_{x_i}(t\xi')$ is length-minimizing. 

\begin{proposition}
\label{Riemannian lemma}
Assume that, for some open neighborhoods $\U_i\subset B_{g_i}(x_i,T)$ of $x_i$, there exists an isometry $\psi:(\U_1,g_1)\to (\U_2,g_2)$ such that
\begin{itemize}
\item[(i)] $d\psi(x_1)V_1=V_2$,

\item[(ii)] $d_{g_1}(y,\exp_{x_1}(\xi'))=d_{g_2}(\psi(y),\exp_{x_2}(d\psi(x_1)\xi'))$ for all $y\in U_1$ and $\xi'\in V_1$.
\end{itemize}
Then, there is an extension of $\psi$ to an isometry 
\[\psi:(B_{g_1}(x_1,T),g_1)\to (B_{g_2}(x_2,T),g_2).\]
\end{proposition}
\begin{proof}
Since the Riemannian manifolds $(M_i,g_i)$ are complete, every point of $B_{g_i}(x_i,T)$ can be connected to $x_i$ by means of a length-minimizing geodesic. Namely, \[\exp_{x_i}(V_i)=B_{g_i}(x_i,T).\]
We claim that there exists a continuous bijection $\phi:B_{g_1}(x_1,T)\to B_{g_2}(x_2,T)$ such that
\begin{align*}
 \phi(\exp_{x_1}(\xi'))=\exp_{x_2}(d\psi(x_1)\xi'),\qquad \forall \xi'\in V_1
\end{align*}
Indeed, since the isometry $\psi$ maps geodesics to geodesics, such a $\phi$ is clearly well defined on a Riemannian ball $B_{g_1}(x_1,\epsilon)\subset \U_1$, and indeed
$ \phi|_{B_{g_1}(x_1,\epsilon)}=\psi|_{B_{g_1}(x_1,\epsilon)}$.
Assume that there exist two distinct vectors $\xi',\eta'\in V_1$ such that 
$y_1:=\exp_{x_1}(\xi')=\exp_{x_1}(\eta')$. All we need to show in order to have a well defined continuous bijection $\phi$ is that 
\begin{align}
\label{e:bijection_holds}
\exp_{x_2}(d\psi(x_1)\xi')=\exp_{x_2}(d\psi(x_1)\eta').
\end{align}
We set $z_1:=\exp_{x_1}(\delta\eta')\in U_1$ for some $\delta>0$ small enough, $z_2:=\exp_{x_2}(d\psi(x_1)\delta\eta')$, and $y_2:=\exp_{x_2}(d\psi(x_1)\xi')$. We have
\begin{align*}
 d_{g_2}(x_2,y_2) & = \|\xi'\|_{g_1} = d_{g_1}(x_1,y_1) = d_{g_1}(x_1,z_1) + d_{g_1}(z_1,y_1)
 = d_{g_2}(x_2,z_2) + d_{g_2}(z_2,y_2),
\end{align*}
where the last equality follows by assumption (ii). This implies that $z_2$ lies on a length-minimizing geodesic segment joining $x_2$ and $y_2$, and therefore~\eqref{e:bijection_holds} holds.

We now prove that $\phi$ is a diffeomorphism. Fix an arbitrary point $y_1\in B_{g_1}(x_1,T)$, and consider the unit-speed length-minimizing geodesic $\gamma_1:[0,\ell]\to M_1$ joining $x_1=\gamma_1(0)$ and $y_1=\gamma_1(\ell)$. We fix $\delta>0$ small enough so that $z_1:=\gamma_1(\delta)\in B_{g_1}(x_1,\epsilon)$. Notice that $\gamma_1|_{[\delta,\ell]}$ is the unique length-minimizing geodesic segment joining $z_1$ and $y_1$, and does not contain conjugate points. Analogously, $\gamma_2:=\phi\circ\gamma_1|_{[\delta,\ell]}$ is the unique length-minimizing geodesic segment joining $z_2:=\phi(z_1)$ and $y_2:=\phi(y_1)$, and does not have conjugate points. Therefore, the Riemannian distances $d_{g_i}$ are smooth in an open neighborhood $Z_i\times Y_i$ of $(z_i,y_i)$. We take $Z_i$ to be small enough so that $Z_1\subset B_{g_1}(x_1,\epsilon)$ and $\phi(Z_1)= Z_2$. For each $w\in Z_i$, the derivative $\partial_{y_i} d_{g_i}(w,y_i)\in S_{y_i}^*M$ is precisely the $g_i$-dual of the tangent vector to the unique length-minimizing geodesic segment joining $w$ and $y_i$. Therefore, for a generic triple of points $w_1,w_2,w_3\in Z_1$ sufficiently close to $z_1$, the triple 
\[\partial_{y_1} d_{g_1}(w_1,y_1),\partial_{y_1} d_{g_1}(w_2,y_1),\partial_{y_1} d_{g_1}(w_3,y_1)\] 
is a basis of $T_{y_1}^*M_1$, and the triple 
\[\partial_{y_2} d_{g_2}(\phi(w_1),y_2),\partial_{y_2} d_{g_2}(\phi(w_2),y_2),\partial_{y_2} d_{g_2}(\phi(w_3),y_2)\]
is a basis of $T_{y_2}^*M_2$. This implies that, up to shrinking the open neighborhoods $Y_i$ of $y_i$, the maps
\begin{align*}
\kappa_1:Y_1\to\R^3,\qquad\kappa_1(y)&=(d_{g_1}(w_1,y),d_{g_1}(w_2,y),d_{g_1}(w_3,y)),\\
\kappa_2:Y_2\to\R^3,\qquad\kappa_2(y)&=(d_{g_2}(\phi(w_1),y),d_{g_2}(\phi(w_2),y),d_{g_2}(\phi(w_3),y))
\end{align*}
are diffeomorphisms onto their images. By assumption~(ii), we have $\kappa_1=\kappa_2\circ\phi|_{Y_1}$, and therefore $\phi|_{Y_1}=\kappa_2^{-1}\circ\kappa_1$. This proves that $\phi$ is a local diffeomorphism at $y_1$. Moreover, for all $z\in Z_1$, by differentiating the equality $d_{g_1}(z,y_1)=d_{g_2}(\phi(z),\phi(y_1))$ with respect of $y_1$, we obtain
\begin{align*}
 \underbrace{\partial_{y_1} d_{g_1}(z,y_1)}_{\in S^*_{y_1}M_1}=\partial_{y_2}d_{g_2}(\phi(z),y_2) d\phi(z_1)= d\phi(z_1)^*\underbrace{\partial_{y_2}d_{g_2}(\phi(z),y_2)}_{\in S^*_{y_2}M_2}.
\end{align*}
Since this holds for all $z$ in the open set $Z_1$, we obtain non-empty open sets $W_i\subset S^*_{y_i}M_i$ such that $d\phi(y_1)^* W_2=W_1$. This implies $d\phi(y_1)^* S^*_{y_2}M_2=S^*_{y_1}M_1$, and therefore $(\phi^* g_2)|_{y_1}=g_1|_{y_1}$. Since $y_1$ is arbitrary, we conclude that $\phi$ is a Riemannian isometry as claimed.
\end{proof}

\bibliography{_biblio}
\bibliographystyle{amsalpha}

\end{document}